\RequirePackage{fix-cm}
\documentclass[smallextended]{svjour3}       
\smartqed  
\usepackage[title]{appendix}
\usepackage{cite}
\usepackage{amsmath}
\numberwithin{equation}{section}
\usepackage{float}
\usepackage{graphicx}
\usepackage{clrscode}
\usepackage{subfigure}
\usepackage{amsfonts}
\usepackage{mathrsfs}
\usepackage{latexsym}
\usepackage{amssymb}
\usepackage{color}
\usepackage{caption}
\usepackage{hyperref}
\usepackage{array}
\usepackage{algorithm}
\usepackage{algpseudocode}
\usepackage{amsmath}
\usepackage{graphics}
\usepackage{epsfig}
\usepackage{geometry}
\usepackage{booktabs}

\usepackage{lineno}
\definecolor{b}{rgb}{0,0,0}

\definecolor{r}{rgb}{0,0,0}

\definecolor{g}{rgb}{0,0,0}

\begin{document}
\newcommand{\bul}{\bu^\lambda}
\newcommand{\bum}{\bu^{\mu}}
\newcommand{\buleta}{\bul(\bet)}
\newcommand{\bldn}{\bld{n}}
\newcommand{\bxzero}{\bld{x}_0}
\newcommand{\bxone}{\bld{x}_1}
\newcommand{\bldv}{\bld{v}}
\newcommand{\bldgam}{\mbox{\boldmath$\gamma$}}
\newcommand{\bgam}{\mbox{\boldmath$\gamma$}}
\newcommand{\bal}{\mbox{\boldmath$\alpha$}}
\newcommand{\bet}{\bld{\eta}}
\newcommand{\bldeta}{\bld{\eta}}
\newcommand{\bnu}{\bld{\nu}}
\newcommand{\bom}{\bld{\omega}}
\newcommand{\A}{ {\alpha} }
\newcommand{\uA}{ {\underline{A}} }
\newcommand{\Ab}{ {\overline{A}} }
\newcommand{\ab}{\bar{\alpha}}
\newcommand{\ah}{\hat a(0)}
\newcommand{\ahat}{\hat a}
\newcommand{\Ah}{ \hat{\alpha} }
\newcommand{\ahp}{\hat a\prime(0)}
\newcommand{\ahpp}{\hat a{\prime\prime}(0)}
\newcommand{\Ai}{{\rm Ai}}
\newcommand{\AjN}{A_j^{(N)}}
\newcommand{\AkN}{A_k^{(N)}}
\newcommand{\AkphN}{A_{k-\hf}^{(N)}}
\newcommand{\al}{\alpha}
\newcommand{\ang}[1]{\langle #1 \rangle }
\newcommand{\AN}[1]{\stackrel{N}{A\,}_{#1}\!}
\newcommand{\ANt}[1]{\stackrel{{\tilde N}}{A\,}_{#1}\!}
\newcommand{\atz}{{\tilde a}_0}
\newcommand{\B}{ {\beta} }
\newcommand{\ba}{\begin{array}}
\newcommand{\bA}{{A}}
\newcommand{\bAN}{\stackrel{N}{\bld{A}\,}\!\!}
\newcommand{\bANt}{\stackrel{{\tilde N}}{\bld{A}\,}\!\!}
\newcommand{\bApN}{\stackrel{N}{\bld{A}}\!\!{}^T\!}
\newcommand{\bb}{\bar{\beta}}
\newcommand{\bBM}{\stackrel{M}{\bld{B}\,}\!\!\!}
\newcommand{\bc}{{\bld c}}
\newcommand{\bchi}{{ \mbox{\boldmath $\chi$} }}
\newcommand{\bd}{\bld{d}}
\newcommand{\bD}{\bld{D}}
\newcommand{\be}{\beta}
\newcommand{\bE}{\bld{E}}
\newcommand{\beq}[1]{ 
        \begin{equation} \label{#1} }
\newcommand{\beql}[1]{
\begin{equation} \label{#1}}

\newcommand{\beqa}{\begin{eqnarray}}
\newcommand{\blde}{{\bld e}}
\newcommand{\beps}{{\bld \epsilon}}
\newcommand{\bldf}{{\bld f}}
\newcommand{\BF}[1]{ \par \indent \begin{fgr} \refstepcounter{fgr}
        \label{#1} {\bf \thefgr.} }
\newcommand{\bh}{{\hat \beta}}
\newcommand{\bld}[1]{{ \mbox{\boldmath $#1$} }}
\newcommand{\bnab}{\mbox{\boldmath$\nabla$}}
\newcommand{\bcdot}{{\bld{\cdot}}}
\newcommand{\barn}{{\bar n}}
\newcommand{\barbn}{{\bar \bn}}
\newcommand{\barbx}{{\bar \bx}}
\newcommand{\bhtal}{\hat{\bld{\alpha}}}
\newcommand{\bldscr}[1]{{ \mbox{\boldmath ${\scriptstyle #1}$} }}
\newcommand{\bF}{\bld{F}}
\newcommand{\bG}{\bld{G}}
\newcommand{\bH}{\bld{H}}
\newcommand{\bI}{\bld{I}}
\newcommand{\bk}{\bld{k}}
\newcommand{\bl}{\bld{l}}
\newcommand{\BM}[1]{\stackrel{M}{B\,}_{#1}\!}
\newcommand{\BMt}[1]{\stackrel{{\tilde M}}{B\,}_{#1}\!}
\newcommand{\bn}{{\bld n}}
\newcommand{\bp}{\bld{ p}}
\newcommand{\bP}{{\bld P}}
\newcommand{\bptl}{\bld{\partial}}
\newcommand{\bq}{\bld{q}}
\newcommand{\bldr}{\bld{r}}
\newcommand{\Brm}{{\rm B}}
\newcommand{\bs}{\bld{s}}
\newcommand{\bS}{\bld{S}}
\newcommand{\bsigma}{\bld{\sigma}}
\newcommand{\bsig}{\bld{\sigma}}
\newcommand{\bt}{\bld{t}}
\newcommand{\bT}{\bld{T}}
\newcommand{\btau}{\bld{\tau}}
\newcommand{\btdal}{\tilde{\bld{\alpha}}}
\newcommand{\bu}{\bld{u}}
\newcommand{\bU}{\bld{U}}
\newcommand{\bup}{\bld{\upsilon}}
\newcommand{\bv}{\bld{v}}
\newcommand{\bV}{\bld{V}}
\newcommand{\bw}{{\bld w}}
\newcommand{\bW}{{\bld W}}
\newcommand{\bx}{\bld{x}}
\newcommand{\bxi}{\bld{\xi}}
\newcommand{\bxidotx}{\bxi\cdotb\bx}
\newcommand{\bX}{\bld{X}}
\newcommand{\by}{\bld{y}}
\newcommand{\bz}{\bld{z}}
\newcommand{\bY}{\bld{Y}}
\newcommand{\bzeta}{\bld{\zeta}}
\newcommand{\bzero}{{\bf 0}}
\newcommand{\cdotb}{{\bld{\cdot}}}
\newcommand{\Chapter}[1]{\setcounter{equation}{0} \chapter{#1}}
\newcommand{\cth}{{\cos\theta}}
\newcommand{\Db}{ \bar{\delta} }
\newcommand{\dd}{\,{\rm d}}
\newcommand{\deldel}[2]{ \fr{\partial #1}{\partial #2} }
\newcommand{\del}[1]{ {\partial_#1} }
\newcommand{\de}{\delta}
\newcommand{\dee}[1]{ \,{\rm d}#1 }
\newcommand{\drm}[1]{ \,{\rm d}#1 }
\newcommand{\dlt}{\delta}
\newcommand{\Dh}{ \hat{\delta} }
\newcommand{\Dlt}{\Delta}
\newcommand{\dol}[1]{\overline{\overline{#1}}}
\newcommand{\e}[1]{\bld{e}_{#1}}
\newcommand{\edot}[1]{\dot{\bld{e}}_{#1}}
\newcommand{\ue}{ {\underline{e}} }
\newcommand{\ea}{\end{array}}
\newcommand{\eeq}{ \end{equation} }
\newcommand{\eeqa}{ \end{eqnarray} }
\newcommand{\EF}{ \end{fgr} }
\newcommand{\eps}{ {\epsilon} }
\newcommand{\eq}[2]{
        \vspace{24pt}
        \marginpar
        [\hfill{\footnotesize #1}]
        {\footnotesize #1\hfill} \vspace{-24pt}
        \begin{equation} \label{#1} 
        #2
        \eeq}
\newcommand{\erm}{{\rm e}}
\newcommand{\erfc}{{\rm erfc}}
\newcommand{\eiop}{{\rm e}^{\im\om\psi}}
\newcommand{\eiot}{{\rm e}^{{\rm i}\omega t}}
\newcommand{\F}[1]{ {\cal F}\{#1\} }
\newcommand{\fb}{\bld{f}}
\newcommand{\uf}{\ul{f}}
\newcommand{\fr}[2]{{ \displaystyle \frac{#1}{#2} }}
\newcommand{\frr}[2]{{\textstyle \frac{#1}{#2} }}
\newcommand{\frscr}[2]{{\scriptstyle \frac{#1}{#2} }}
\newcommand{\frtxt}[2]{{\textstyle \frac{#1}{#2} }}
\newcommand{\g}[1]{ {\gamma_{#1}} }
\newcommand{\ga}{ {\gamma_{\bar{\alpha}}} }
\newcommand{\gam}{\gamma}
\newcommand{\gams}{\gamma^2}
\newcommand{\gb}{ {\gamma_{\bar{\beta}}} }
\newcommand{\Gb}{ \bar{\gamma} }
\newcommand{\gef}{\hat{\gamma}^{\rm eff}}
\newcommand{\geff}{  { \gamma_{{\rm eff}} }  }
\newcommand{\Gh}{ \hat{\gamma} }
\newcommand{\hbxi}{ \hat{\bxi} }
\newcommand{\hbxidotx}{\hbxi\cdotb\bx}
\newcommand{\hem}{\mbox{\hspace{1em}}}
\newcommand{\hf}{{ \textstyle \frac{1}{2} }}
\newcommand{\Hf}{\frac{1}{2}}
\newcommand{\hff}{{ \scriptstyle \frac{1}{2} }}
\newcommand{\hftxt}{{ \textstyle \frac{1}{2} }}
\newcommand{\hfdis}{{ \displaystyle \frac{1}{2} }}
\newcommand{\hfscr}{{ \scriptstyle \frac{1}{2} }}
\newcommand{\rh}{{\rm h}}
\newcommand{\Hrm}{{\rm H}}
\newcommand{\hw}{\hat w}
\newcommand{\hxi}{\hat \xi}
\newcommand{\im}{{\rm i}}
\newcommand{\Imag}[1]{ \Im\{{#1}\} }
\newcommand{\imm}{{\cal im}}
\newcommand{\INT}{{\displaystyle \int}}
\newcommand{\itf}{$\mbox{in the form}$}
\newcommand{\io}{{{\rm i}\omega}}
\newcommand{\dint}{ \displaystyle{\int} }
\newcommand{\Int}[2]{{\displaystyle \int_{#1}^{#2}}}
\newcommand{\LA}{ {\lambda} }
\newcommand{\la}{ {\langle} }
\newcommand{\Lb}{ \bar{\lambda} }
\newcommand{\Lh}{ \hat{\lambda} }
\newcommand{\Lim}{{\displaystyle \lim_{N\to\infty}}}
\newcommand{\lm}{\lim_{N\to\infty}}
\newcommand{\M}{ {\mu} }
\newcommand{\Mb}{ \bar{\mu} }
\newcommand{\Mh}{ \hat{\mu} }
\newcommand{\n}{\bld{n}}
\newcommand{\N}{\nu}
\newcommand{\nab}{{\nabla}}
\newcommand{\nabs}{{\nabla^2}}
\newcommand{\NM}{^{(\tilde{N}\hat{M})}}
\newcommand{\npvf}{ \vspace*{\fill}\newpage}
\newcommand{\Nth}{ {\frac{1}{N}} }
\newcommand{\om}{\omega}
\newcommand{\omb}{\bar{\omega}}
\newcommand{\oms}{\omega^2}
\newcommand{\pq}{ (p^2+q^2) }
\newcommand{\Prm}{{\rm P}}
\newcommand{\Pro}[1]{{\displaystyle \prod_{#1}}}
\newcommand{\Prod}[2]{{\displaystyle \prod_{#1}^{#2}}}
\newcommand{\PS}{ {\psi} }
\newcommand{\psib}{ \bar{\psi} }
\newcommand{\qt}{ {\frac{1}{4}} }
\newcommand{\qttxt}{ \textstyle{\frac{1}{4}} }
\newcommand{\Qt}{ \frac{1}{4} }
\newcommand{\re}{{\rm re}}
\newcommand{\ra}{{\rangle}}
\newcommand{\rH}{{\rm H}}
\newcommand{\rJ}{{\rm J}}
\newcommand{\Real}[1]{\Re{\{ {#1} \}} }
\newcommand{\real}[1]{\Re{\{ {#1} \}} }
\newcommand{\rf}[1]{(\ref{#1})}
\newcommand{\refrf}[1]{[\ref{#1}]}
\newcommand{\rF}[1]{\ref{#1}}
\renewcommand{\r}{\rho}
\newcommand{\rhob}{\bar{\rho}}
\newcommand{\rhoh}{\hat{\rho}}
\newcommand{\rkN}{r_k^{(N)}}
\newcommand{\rkhN}{r_{k-\hf}^{(N)}}
\newcommand{\rkphN}{r_{k+\hf}^{(N)}}
\newcommand{\rmbox}[1]{{\mbox{\rm #1}}}
\newcommand{\rotn}{{\mbox{\boldmath $\Omega$}}}
\newcommand{\s}{\sigma}
\newcommand{\us}{ \underline{s}}
\newcommand{\Section}[1]{\setcounter{equation}{0} \section{#1}}
\newcommand{\sgn}{{\rm sgn}}
\newcommand{\Sm}[2]{{\displaystyle \sum_{#1}^{#2}}}
\newcommand{\Srm}{{\rm S}}
\newcommand{\ssf}[1]{{\sf #1}}
\newcommand{\sth}{{\sin\theta}}
\newcommand{\Strain}{{\mbox{\boldmath $\varepsilon$}}}
\newcommand{\strain}{{\varepsilon}}
\newcommand{\Stress}{{\mbox{\boldmath $\tau$}}}
\newcommand{\stress}{\tau}
\newcommand{\PKStress}{\bld{T}}
\newcommand{\PKstress}{T}
\newcommand{\ptl}[1]{ {\partial_#1} }
\newcommand{\Sum}[1]{{\displaystyle \sum_{#1}}}
\newcommand{\T}{\bld{T}}
\newcommand{\taub}{\bar{\tau}}
\newcommand{\taupm}{\tau\prime}
\newcommand{\tb}{\bld{t}}
\renewcommand{\theequation}{\arabic{section}.\arabic{equation}}
\newcommand{\ti}{\theta_i}
\newcommand{\tu}{{\tilde u}}
\newcommand{\tmu}{{\tilde \mu}}
\newcommand{\tnu}{{\tilde \nu}}
\newcommand{\Twmw}{$\mbox{Thus we may write}$}
\newcommand{\txN}{{\tilde x}^{(N)}}
\newcommand{\ul}[1]{\und{#1}}
\newcommand{\und}[1]{\underline{#1}}
\newcommand{\uu}{\ul{u}}
\newcommand{\uxi}{\bld{\xi}}
\newcommand{\ux}{ \bld{x}}
\newcommand{\vb}{\bar{v}}
\newcommand{\vf}{\vfill}
\newcommand{\vsf}{\vspace*{\fill}}
\newcommand{\vh}{\hat{v}}
\newcommand{\Vh}{\hat{V}}
\newcommand{\wmw}{$\mbox{we may write}$}
\newcommand{\wrt}{$\mbox{with respect to}$}
\newcommand{\xb}{\bld{x}}
\newcommand{\x}{\xi}
\newcommand{\X}{\bld{X}}
\newcommand{\xki}{\xi_{k,i}}
\newcommand{\xN}[1]{x_{#1}^{(N)}}
\newcommand{\xNxi}[1]{{x_{#1}^{(N)} (\X)}}
\newcommand{\y}{\bld{y}}
\newcommand{\Y}{\bld{Y}}
\newcommand{\yb}{\bld{y}}
\newcommand{\Z}{ {\zeta} }
\newcommand{\Zb}{ {\overline{\zeta}} }
\newcommand{\ZkmN}{ {\zeta _{k-1}^{(N)}} }
\newcommand{\ZkN}{ {\zeta _k^{(N)}} }
\newcommand{\z}[1]{z_{#1}}

\newcommand{\steq}[1]{\stackrel{\scriptscriptstyle{\mathrm{#1}}}{=}}
\newcommand{\stlp}[1]{\stackrel{\scriptscriptstyle{\mathrm{#1}}}{(}}
\newcommand{\strp}[1]{\stackrel{\scriptscriptstyle{\mathrm{#1}}}{)}}
\newcommand{\stk}[2]{\hspace{1pt}\stackrel{\scriptscriptstyle{\mathrm{#1}}}{#2}\hspace{1pt}}
\newcommand{\stkrse}[2]{\!\raisebox{.3ex}{\mbox{$\stk{#1}{#2}$}}\!}
\newcommand{\szero}{\bld{\sf 0}}
\newcommand{\sA}{\bld{\sf A}}
\newcommand{\sB}{\bld{\sf B}}
\newcommand{\sC}{\bld{\sf C}}
\newcommand{\sE}{\bld{\sf E}}
\newcommand{\sG}{\bld{\sf G}}
\newcommand{\sI}{\bld{\sf I}}
\newcommand{\sM}{\bld{\sf M}}
\newcommand{\sR}{\bld{\sf R}}
\newcommand{\sS}{\bld{\sf S}}
\newcommand{\sU}{\bld{\sf U}}

\newcommand{\sub}[1]{_{,#1}}

\newcommand{\cA}{{\cal A}}
\newcommand{\calC}{{\cal C}}
\newcommand{\calE}{{\cal E}}
\newcommand{\calF}{{\cal F}}
\newcommand{\calP}{{\cal P}}
\newcommand{\calR}{{\cal R}}
\newcommand{\calS}{{\cal S}}
\newcommand{\calT}{{\cal T}}


\newcommand{\FFig}[4]{
\newpage\centerline{\bf #3}
\begin{figure}[#1]
\centerline{\psfig{height=#2in,file=BurridgeFigure#3.pdf}}
\caption{ \label{BurridgeFigure#3} {\footnotesize #4} \newline  BurridgeFigure#3} 
\end{figure}
}

\newcommand{\cross}{\hspace{-2pt}\bld{\times}\hspace{-1pt}}                              %
\newcommand{\irm}{{\rm i}}                                                                                    %
\newcommand{\aps}{  {
\rotatebox{-90}{ \hspace{-10pt}$\backslash$\hspace{-8pt}$\supset$ }\hspace{-1pt}
}  }

\title{Hadamard integrators for wave equations in time and frequency domain: Eulerian formulations via butterfly algorithms
\thanks{Cheng's research was supported by NSFC 11971121, 12241103, the Science and Technology Commission of Shanghai Municipality 23JC1400501 and the Sino-German Mobility Programme (M-0187) by the Sino-German Center for Research Promotion. S. Leung is supported by the Hong Kong RGC under grant 16302223. Qian's research is partially supported by NSF 2012046, 2152011, and 2309534.}}


\titlerunning{Eulerian Hadamard integrator}        

\author{Yuxiao Wei \and Jin Cheng \and Shingyu Leung   \and Robert Burridge \and Jianliang Qian 
}

\authorrunning{Wei, Cheng, Leung, Burridge, and Qian} 

\institute{Yuxiao Wei  \\ \email{19110180029@fudan.edu.cn}  
\and 
\\ Jin Cheng \\ \email{jcheng@fudan.edu.cn}
     \and
\\Shingyu Leung \\ \email{masyleung@ust.hk}
\and 
\\ Robert Burridge \\ \email{burridge137@gmail.com}
\and 
\\Jianliang Qian \\ \email{jqian@msu.edu}
\\
\\
Yuxiao Wei \and Jin Cheng  
\at  School of Mathematical Sciences, Fudan University, Shanghai 200433, China.
\and
Shingyu Leung 
\at  Department of Mathematics, The Hong Kong University of Science and Technology, Clear Water Bay, Hong Kong. 
\and 
Robert Burridge
\at Department of Mathematics and Statistics, University of New Mexico, Albuquerque, NM 87131, USA. 
\and 
Jianliang Qian 
\at Department of Mathematics and Department of CMSE, Michigan State University, East Lansing, MI 48824, USA
}

\date{Received: date / Accepted: date}

\maketitle

\begin{abstract}
Starting from Kirchhoff-Huygens representation and Duhamel's principle of time-domain wave equations, we propose novel butterfly-compressed Hadamard integrators for self-adjoint wave equations in both time and frequency domain in an inhomogeneous medium. First, we incorporate the leading term of Hadamard's ansatz into the Kirchhoff-Huygens representation to develop a short-time valid propagator. Second, using Fourier transform in time, we derive the corresponding Eulerian short-time propagator in the frequency domain; on top of this propagator, we further develop a time-frequency-time (TFT) method for the Cauchy problem of time-domain wave equations. Third, we further propose a time-frequency-time-frequency (TFTF) method for the corresponding point-source Helmholtz equation, which provides Green's functions of the Helmholtz equation for all angular frequencies within a given frequency band. Fourth, to implement the TFT and TFTF methods efficiently, we introduce butterfly algorithms to compress oscillatory integral kernels at different frequencies. As a result, the proposed methods can construct wave field beyond caustics implicitly and advance spatially overturning waves in time naturally with quasi-optimal computational complexity and memory usage. Furthermore, once constructed the Hadamard integrators can be employed to solve both time-domain wave equations with various initial conditions and frequency-domain wave equations with different point sources. Numerical examples for two-dimensional wave equations illustrate the accuracy and efficiency of the proposed methods.
\end{abstract}

\section{Introduction}
\label{intro}
We consider the Cauchy problem for the self-adjoint wave equation in $m$-dimensional space $\mathbb{R}^m$,
\begin{equation}\label{1.1}
  \rho u_{tt}-\nab\cdot(\nu \nabla u)=f,\; \boldsymbol{x}\in \mathbb{R}^m,\; t>0
\end{equation}
with initial conditions
\begin{equation}\label{1.2}
  u(0,\boldsymbol{x})=u^1(\boldsymbol{x}),\; u_t(0,\boldsymbol{x})=u^2(\boldsymbol{x}),
\end{equation}
where $t$ is time, the subscripts ${}_t$ and ${}_{tt}$ represent the first and second time derivatives, respectively, position $\boldsymbol{x}=[x_1,x_2,\cdots,x_m]^T$, the gradient operator $\nab =[\partial_{x_1},\partial_{x_2},\cdots,\partial_{x_m}]^T$, both variables $\rho$ and $\nu$ are analytic and positive functions of position $\boldsymbol{x}$, characterizing certain physical parameters of the medium, $f(t,\bx)$ is the source term, and $u^1(\boldsymbol{x})$ and $u^2(\boldsymbol{x})$ are compactly supported, highly oscillatory smooth functions. 

Taking 
\begin{equation}\label{1.3.0}
    f(t,\bx)=\delta(t)\delta(\bx-\bz)\quad\mbox{and}\quad u^1=u^2=0\quad
\end{equation}
in equations \eqref{1.1} and \eqref{1.2} and applying Fourier transform in time to equation \eqref{1.1} accordingly, we obtain the frequency-domain wave equation, the Helmholtz equation, 
\begin{equation}\label{1.5}
  \nab\cdot(\nu\nab \hat{u})+\omega^2\rho \hat{u}=-\delta(\bx-\bz),
\end{equation}
with the Sommerfeld radiation condition imposed at infinity, where $\omega$ is the angular frequency and $\bz$ is the source.

When the initial conditions in \rf{1.2} are highly oscillatory or the angular frequency in \rf{1.5} is large, the wave field will be highly oscillatory; however, direct numerical methods, such as finite-difference or finite-element methods, for such problems may suffer from dispersion or pollution errors \cite{baygoltur85,babsau00}, so that such methods require an enormous computational grid to resolve these oscillations and are thus very costly in practice. Therefore, alternative methods, such as geometrical-optics-based asymptotic methods, have been sought to resolve these highly oscillatory wave phenomena. 

\color{b}
The essential difficulty in applying geometric optics is that it cannot handle caustics easily \cite{lax57,lud66,masfed81,babbul72,kellew95,ben03,engrun03,luoqiabur14a,luqiabur16,weicheburqia24}. 
One approach is to replace the geometrical-optics ansatz with Gaussian beam (GB) summations \cite{tanqiaral07,tanengtsa09,qiayin10a,qiayin10b,baolaiqia14,qiason21}. Eulerian formulations of GB summations have been developed for Helmholtz equations in the high-frequency regime \cite{leuqiabur07, leuqia09, leuzha09a}  and time-dependent Schrödinger equations in the semi-classical regime \cite{leuqia10}. The main challenge of this approach is that the computations must be done in the high-dimensional phase space. Another approach is to incorporate the geometrical-optics ansatz into the Huygens secondary-source principle, resulting in fast Huygens sweeping (FHS) methods, and these methods have been designed to solve time-dependent Schrödinger equations \cite{leuqiaser13, holeuqia21}, Helmholtz equations \cite{luoqiabur14a, luqiabur16}, frequency-domain Maxwell's equations \cite{qialuyualuobur16, luqiabur18}, time-harmonic vectorial Maxwell's equations \cite{kwaleuwanqia17}, and frequency-domain elastic wave equations \cite{qiasonlubur23}, all in the presence of caustics. FHS methods can propagate wave fields through appropriately partitioned spatial layers by marching in that preferred spatial direction in a layer-by-layer fashion, but these methods only work in media where geodesics satisfy the sub-horizontal condition \cite{symqia03slowm}.

Although caustics occur with a high probability for wave propagation in inhomogeneous media \cite{whi84}, we note that we are still able to use the geometrical-optics type method mainly because of the following fact \cite{avikel63,symqia03slowm}: in an isotropic medium, the point-source eikonal equation has a locally smooth solution near the source point except the source point itself; this implies that caustics will not develop right away on the expanding wavefront away from the source. If we compute the geometrical-optics ansatz of frequency-domain Green's function of \eqref{1.5} in these caustic-free neighborhoods of point sources, for example, the WKBJ ansatz or Babich's ansatz, the locally valid Green's functions can be incorporated into the Huygens secondary-source principle to update the global wave field, resulting in fast Huygens sweeping (FHS) methods for Helmholtz equations \cite{luoqiabur14a, luqiabur16}; alternatively, incorporating the Hadamard's ansatz \cite{had23,CouHil64} for time-domain Green's function of \eqref{1.1} into the Huygens secondary-source principle results in the Hadamard integrator\cite{weicheburqia24} for the wave equation.  We denote by $\bar{T}(\bx_0)$ the time when the first caustic occurs for rays issued from a source point $\bx_0$. For $\Delta T<\min_{\bx_0\in\Omega}\bar{T}(\bx_0)$, where $\bx_0$ is located in a certain bounded domain $\Omega$, the Hadamard's ansatz is valid in this domain $\Omega$ so that we can propagate wave forward in time with the time step $\Delta T$ using the Huygens secondary-source principle.

 Although the frequency-domain Green's function is highly oscillatory, the kernel matrix generated on a uniform Eulerian grid can be constructed and applied with a quasi-linear computational cost using butterfly decomposition. Thus, the FHS method \cite{luoqiabur14a, luqiabur16} can construct high-frequency wave fields with quasi-linear computational complexity and memory usage. However, as a frequency-domain method, the FHS method propagates the wave field layer by layer along a preferred spatial direction, and the propagation along the artificial spatial direction requires geodesics to satisfy the sub-horizontal assumption, implying that the FHS method is applicable only to a limited set of specific media.

Although the time-domain Green's function is non-oscillatory in time direction, it has singularity on irregular wavefronts in time. The Hadamard integrator developed in \cite{weicheburqia24} propagates the wave fields along the \textit{natural} time direction so that the sub-horizontal assumption of geodesics in the FHS methods \cite{luoqiabur14a, luqiabur16} is no longer required in this new setting. However, the singularity on irregular wavefronts prevents the Lagrangian Hadamard integrator \cite{weicheburqia24} from achieving quasi-linear computational complexity.

Then we immediately run into a question: is it possible to develop an Eulerian Hadamard integrator and use butterfly decomposition techniques similar to those in the FHS methods to achieve quasi-linear computational cost?
The answer is yes! Our strategy is based on the following two observations: (1) the time-domain Green's function of the wave equation is the inverse  Fourier transform in time of the corresponding frequency-domain Green's function of the Helmholtz equation, and as shown in \cite{qiasonlubur21}, each term in Hadamard's ansatz for the time-domain Green's function is indeed the inverse Fourier transform in time of the corresponding term in Babich's ansatz \cite{bab65,luqiabur16,luqiabur18} for the frequency-domain Green's function; (2) moreover, while the time-domain Green's function has singularity concentrated on irregular curved wavefronts which is unfriendly to numerical quadrature posed on regular mesh points, the frequency-domain Green's function only has a singularity at the source point so that it is amenable to numerical quadrature defined on regular mesh points. These two observations motivate us to take a detour to apply the time-domain Green's function efficiently by first carrying out relevant applications in the frequency domain and then coming back to the time domain: first, take the Fourier transform in time to go to the frequency domain; second, apply the frequency-domain Green's function to corresponding arguments; third, take the inverse Fourier transform in frequency to get back to time domain. This strategy leads to the time-frequency-time (TFT) method for the Cauchy problem of the wave equation. 

Different from the Lagrangian Hadamard integrator developed in \cite{weicheburqia24}, the TFT method provides resolution-satisfying time sampling of wave solutions, implying that the Hadamard integrator is also a time-domain method for solving the Helmholtz equation \eqref{1.5}. After carefully handling the singularity of the delta function, we can use the TFT method to obtain the time-domain global Green's function numerically. By performing the Fourier transform in time, we obtain the frequency-domain Green's function, which is the solution to \eqref{1.5}. We dub this strategy the time-frequency-time-frequency (TFTF) method.

To accelerate matrix-vector multiplications in the TFT and TFTF methods, we first construct low-rank separable representations of the Hadamard ingredients, including phase and Hadamard coefficients in Hadamard's ansatz, by solving the corresponding governing equations and utilizing the Chebyshev spectral expansions accordingly. Next we consider an algebraic compression tool called butterfly \cite{micboa96,liyaneilkenyin15,lizha17,candemyin09}, a multilevel numerical linear algebra algorithm well-suited for representing highly oscillatory operators, to compress the Hadamard integrator. We have applied the method in \cite{luoqiabur14a,luqiabur16,liusonburqia23} to deal with integrals of Green's functions of the Helmholtz equation. to deal with integral kernels with singularities, we combine the hierarchical matrix representation and the butterfly algorithm by designing a suitable strategy for selecting proxy rows in the interpolation decomposition of butterfly construction, which allows the butterfly compressed Hadamard integrator to enjoy the quasi-linear complexity in the overall CPU time and memory usage with respect to the number of unknowns; such quasi-optimal complexity is valid for both the wave equation in the time-space domain and the Helmholtz equation in the frequency-space domain. 

It is important to note that the Hadamard integrator is independent of initial conditions. Once constructed, it can be used to solve both time-domain wave equations with various initial conditions and frequency-domain wave equations with different point sources. Moreover, once these butterfly representations are constructed, they can be applied recursively when executing the Hadamard integrator; this unique feature allows the Hadamard integrator to be efficiently parallelized.

This paper develops an Eulerian Hadamard integrator. The approach inherits the advantages of the Lagrangian Hadamard integrator in \cite{weicheburqia24} by marching in time direction; moreover, since the Eulerian integrator relies on the butterfly algorithm to speed up the matrix-vector multiplication arising from the integration process, the computational complexity of the Eulerian integrator for the wave equation in both time and frequency domain can be reduced to quasi-linear. Moreover, our Eulerian Hadamard integrator is also able to solve the point-source Helmholtz equation in inhomogeneous media by treating caustics implicitly and handling overturning waves naturally, thus settling the issue of how to remove the sub-horizontal condition in the fast Huygens sweeping methods \cite{luoqiabur14a,luqiabur16,qialuyualuobur16,luqiabur18} to accommodate overturning waves.  

 The rest of the paper is organized as follows. We introduce in Section \ref{khform} the Kirchhoff-Huygens representation formula which utilizes Green's functions to propagate waves. Connecting Hadamard's ansatz and Babich's ansatz using Fourier transform, we propose in Section \ref{sec3} a novel Eulerian Hadamard-Kirchhoff-Huygens (HKH) propagator for the Cauchy problem of the time-dependent wave equation and the point-source problem of the frequency-dependent wave equation (the Helmholtz equation), and we further obtain the Eulerian Hadamard integrators by applying it recursively. We then develop in Section \ref{sec4} numerical strategies for implementing the Eulerian Hadamard integrators. To accelerate the constructions and applications of the integrators, in Section \ref{sec5}, we construct interpolative decomposition butterfly (IDBF) and hierarchically off-diagonal butterfly (HODBF) representations of the integral kernels, resulting in quasi-optimal computational complexity and memory usage. Section \ref{sec6} presents two-dimensional (2-D) results to demonstrate the accuracy, efficiency and convergence of the proposed methods. We conclude the paper with some comments in Section \ref{sec7}.
  \color{black}
\section{Kirchhoff-Huygens representation formula}
 \label{khform}
 We are interested in solving the Cauchy problem for the self-adjoint wave equation \eqref{1.1} 
with initial conditions \eqref{1.2}. 
Let $S(t)$ be the solution operator corresponding to $f=0$ and $u^1=0$; namely, for a function $\psi$, $u(t, \bx)=(S(t) \psi)(\bx)$ solves \rf{1.1} with the source term $f=0$ and the initial conditions $u^1=0$ and $u^2=\psi$. Then Duhamel's principle implies that the solution for wave equation \rf{1.1} with initial conditions \rf{1.2} is
\begin{equation}
u(t,\bx_0)=\partial_t(S(t) u^1)(\bx_0)+(S(t) u^2)(\bx_0)+\int_0^t\left(S(t-s) f(s,\cdot)\right)(\bx_0) d s.
\end{equation}
Consider the time-domain Green's function $G(t, \bx_0; \bx)$ which satisfies equation \rf{1.1} with the source term and initial conditions specified in \rf{1.3.0}. Then, for a compactly supported function $\psi$ such that the initial condition $u^2=\psi$, by \cite{weicheburqia24} the solution operator $S$ can be written as 
\begin{equation}\label{2.2.1}
    (S(t)\psi)(\bx_0)=\int_V \rho(\bx) G(t,\bx_0;\bx) \psi(\bx) d \bx,
\end{equation}
where $V$ is a region of space which contains the support of $\psi$ and does not change in time. Since the boundary integrals are not considered, we need to set a sufficiently large computational domain to ensure that the waves do not reach the boundary \cite{weicheburqia24}.

Then we have 
\begin{equation}
\label{2.5.0}
\begin{aligned}
u(t,\bx_0)&=\int_V \rho(\bx)[G_t(t,\bx_0;\bx) u^1(\bx)+  G(t,\bx_0;\bx) u^2(\bx)] d \bx+\int_0^t \int_V \rho(\bx) G(t-s,\bx_0;\bx) f(s,\bx) d\bx d s.
\end{aligned}
\end{equation}
Differentiating \rf{2.5.0} with respect to $t$, we obtain
\begin{equation}\label{2.6.0}
\begin{aligned}
    u_t(t,\bx_0)&=\int_V \rho(\bx)[G_{tt}(t,\bx_0;\bx) u^1(\bx)+G_t(t,\bx_0;\bx) u^2(\bx)] d \bx+\int_0^t \int_V \rho(\bx)G_t(t-s,\bx_0;\bx) f(s,\bx) d\bx d s.
\end{aligned}
\end{equation}

We refer to \rf{2.5.0} and \rf{2.6.0} as the Kirchhoff-Huygens representation formula. When the medium is homogeneous, Green's function for the wave equation is known so that the formula has been used frequently in practice. When the medium is inhomogeneous and the source term $f=0$, in \cite{weicheburqia24} we have used Hadamard's ansatz to compute the needed Green's function for the formula, leading to a short-time valid propagator in a Lagrangian formulation, dubbed the Hadamard-Kirchhoff-Huygens (HKH) propagator, which is able to propagate a highly oscillatory wave field for a short period of time; furthermore, a Lagrangian Hadamard integrator is obtained by recursively applying the HKH propagator in time. In this article, we continue to use Hadamard's ansatz to construct Green's functions in the Kirchhoff-Huygens formula, but we compute the HKH propagator using an Eulerian formulation, resulting in an Eulerian Hadamard integrator.

To simplify the presentation, in the current work we only consider the two-dimensional case, where $m=2$, and we leave the three-dimensional case as a future work. Consequently, we consider the following two cases:
\begin{itemize} 
\item Case 1. The pure initial-value problem for equation \eqref{1.1}, where $f=0$, and $u^1$ and $u^2$ are compactly supported, highly oscillatory smooth functions. 
\item Case 2. The point-source problem for equation \eqref{1.5}, where the frequency-domain Green's function of the Helmholtz equation naturally corresponds to the time-domain Green's function of the wave equation \eqref{1.1} with $f=\delta(t)\delta(\bx-\bx_0)$ and $u^1=u^2=0$.  
\end{itemize}

In Case 1 for pure initial-value problems, we employ the same setup as in \cite{weicheburqia24}, aiming to achieve a quasi-optimal computational complexity using Eulerian formulations. 

In Case 2 for point-source problems of the Helmholtz equation, since the time-domain Green's function corresponding to $f=\delta(t)\delta(\bx-\bx_0)$ and $u^1=u^2=0$ in the Cauchy problem \eqref{1.1}-\eqref{1.2} is naturally linked to the frequency-domain Green's function of the Helmholtz equation, we extend the Hadamard integrator to point-source wave equations in frequency domain so as to solve the point-source Helmholtz equation. Because the Hadamard integrator marches in the time direction, it can naturally handle overturning waves in frequency domain, which is in sharp contrast to our previous frequency-domain asymptotic methods in \cite{luoqiabur14a,luqiabur16}, in that the methods in \cite{luoqiabur14a,luqiabur16} are not able to treat overturning waves. 

On the other hand, for a generic source term $f(x,t)\not\equiv 0$ with compact support in both space and time, the Kirchhoff-Huygens representation \rf{2.5.0}-\rf{2.6.0} contains a volume integral of $m+1$ dimension, which is expensive to evaluate; therefore, we leave this as a future work.

\section{Eulerian Hadamard integrators}
\label{sec3}
\subsection{Hadamard's ansatz for time-domain Green's functions}
Consider Green's function of the self-adjoint wave equation 
\begin{equation}
\label{2.1}
  \rho u_{tt}-\nab\cdot(\nu \nab u)=\delta(t)\delta(\bx-\bx_0),\;\; \boldsymbol{x}\in \mathbb{R}^m,\;\; t>0
\end{equation}
with initial conditions 
\begin{equation}
\label{2.2}
  u(0,\boldsymbol{x})=0,\; u_t(0,\boldsymbol{x})=0.
\end{equation}
Let
\begin{equation}\label{2.3}
  c=\sqrt{\frac{\nu}{\rho}}, \quad n=\frac{1}{c}=\sqrt{\frac{\rho}{\nu}} .
\end{equation}
Then Hadamard's ansatz \cite{weicheburqia24} is written as 
\begin{equation}\label{2.4}
  u(t, \boldsymbol{x}_0; \bx)=\sum_{s=0}^{\infty} v_s(\boldsymbol{x}_0; \bx) f_{+}^{s-\frac{m-1}{2}}\left[t^2-\tau^2(\boldsymbol{x}_0; \bx)\right],
\end{equation}
where the Gelfand-Shilov generalized functions $f_{+}^{\lambda}(\zeta)$ are given in \cite{gelshi64}. They are defined for $\lambda>-1$ as follows:
\begin{equation}\label{2.5}
  f_{+}^\lambda(\zeta)=\frac{\zeta_{+}^\lambda}{\lambda !},
\end{equation}
where
\begin{equation}\label{2.6}
  \zeta_{+}^\lambda=\left\{\begin{array}{cc}
0, & \text { for } \zeta<0, \\
\zeta^\lambda, & \text { otherwise},
\end{array}\right.
\end{equation}
and by analytic continuation for other values of $\lambda$. \color{b} We have the important relationship \cite{gelshi64}
\begin{equation}\label{2.8.f}
  (f_{+}^{\lambda})^{\prime}=f_{+}^{\lambda-1}.
\end{equation}
\color{black} 
In particular, $\lambda=0$ corresponds to the Heaviside function, yielding a jump discontinuity at $\zeta=0$,  while by \eqref{2.8.f}, $\lambda=-1$ corresponds to the $\delta$-function.  
Here $\tau(\boldsymbol{x}_0; \bx)$ is the phase function, also known as traveltime, satisfying the eikonal equation \rf{3.7.0}, 
\begin{equation}\label{3.7.0}
    |\nab \tau |=n, \; \tau(\bx_0; \bx)|_{\bx=\bx_0}=0.
\end{equation}
And the coefficient functions $v_s(\bx_0; \bx)$ in \rf{2.4} satisfy the recurrent system
\beq{SASWE3.01}
4\,\rho\,\tau\,\fr{\drm v_{s}}{\drm \tau}\,+\,v_{s}\,[\nab\bcdot(\nu\,\nab\,\tau^2)\,+\,2\,(2s-m)\,\rho] = \nab\bcdot(\nu\nab v_{s-1}), \quad s=0,1,2,\cdots,
\eeq
and $v_{-1}\equiv 0$. By leveraging properties of the Gelfand-Shilov function as demonstrated in \cite{weicheburqia24}, we have the initial condition for $v_0$ at $\bx_0$:
\begin{equation}\label{3.9.0}
    v_0(\bx_0;\bx)|_{\bx=\bx_0}=\frac{n^m(\bx_0)}{2 \rho(\bx_0) \pi^{\frac{m-1}{2}}}=\frac{n^{m-2}(\bx_0)}{2 \nu(\bx_0) \pi^{\frac{m-1}{2}}}.
\end{equation}

 We introduce some definitions and properties of generalized functions \cite{gelshi64} here. We first denote $K$ the space of test functions which consists of all real functions $\phi(\bx)$ with continuous derivatives of all orders and with bounded support. We next define a generalized function as any linear continuous functional on $K$. These functionals (generalized functions) which can be given by an equation such as 
\begin{equation}
    (\Psi,\phi) = \int_{R^n} \Psi(\bx)\phi(\bx) d\bx \quad\mbox{ for }\forall \phi\in K
\end{equation}
with $\Psi$ an absolutely integrable function in every bounded region of $R^n$ shall be called regular, and all others (including the delta function) will be called singular. It is worth noting that $f_{+}^{\lambda}(\zeta)$ is regular when $\lambda >-1$ and singular when $\lambda \leq -1$. 

\subsection{Babich's ansatz for frequency-domain Green's functions}
In order to develop Eulerian formulations, we introduce Babich's ansatz \cite{bab65,luqiabur16,qialuyualuobur16,luqiabur18} for the frequency-domain Green's function $\hat{u}=\hat{G}(\omega,\bx_0;\bx)$ of the Helmholtz equation \rf{1.5}, which is written as  
\begin{equation}\label{3.10.0}
    \hat{u}(\omega,\bx_0; \bx)=\hat{G}(\omega,\bx_0;\bx)=\sum_{s=0}^{\infty}v_s(\bx_0; \bx)f_{s+1-\frac{m}{2}}(\tau(\bx_0;\bx),\omega),
\end{equation}
where $\tau$ still represents the traveltime, $v_s$ represent the coefficient functions, 
\begin{equation}
    f_p(\tau, \omega)=i \frac{\sqrt{\pi}}{2} e^{i p \pi}\left(\frac{2 \tau}{\omega}\right)^p H_p^{(1)}(\omega \tau)=i \frac{\sqrt{\pi}}{2}\left(\frac{2 \tau}{\omega}\right)^p H_{-p}^{(1)}(\omega \tau),
\end{equation}
and $H_p^{(1)}$ is the p-th Hankel function of the first kind. In \cite{luqiabur16}, using properties of the basis $f_{p}(\tau,\omega)$, we have obtained that $\tau$ satisfies the eikonal equation \rf{3.7.0} and $v_s$ satisfy the recurrent system \rf{SASWE3.01} with the initial condition \rf{3.9.0}. That is, the time-independent phase and coefficient functions in Hadamard's ansatz are the same as the frequency-independent phase and coefficient functions in Babich's ansatz. Furthermore, \cite{luqiabur18} shows that the basis functions of Hadamard's ansatz and Babich's ansatz are linked to each other: 
\begin{equation}
    \int_\tau^{\infty} e^{i \omega t} f_{+}^{\left(\nu-\frac{1}{2}\right)}\left(t^2-\tau^2\left(\bx_0; \bx\right)\right) d t=\frac{1}{2} i \sqrt{\pi}\left(\frac{2 \tau}{\omega}\right)^\nu e^{i \omega \nu} H_\nu^{(1)}(\omega \tau) \equiv f_\nu(\omega, \tau).
\end{equation}

Because of the reciprocity, we have $\hat{G}(\omega,\bx_0;\bx)$=$\hat{G}(\omega,\bx;\bx_0)$, which will be used frequently in the following algorithmic development. 
Letting $G_0(t,\bx_0;\bx)$ and $\hat{G}_0(\omega,\bx_0;\bx)$ be the leading term of Hadamard's ansatz \rf{2.4} and Babich's ansatz \rf{3.10.0}, respectively, for a point source located at $\bx_0$, we have the following crucial relations,  
\begin{equation}\label{3.13.0}
    G_0(t,\bx_0;\bx)=\frac{1}{2\pi} \int_{-\infty}^{\infty} \hat{G}_0(\omega,\bx_0;\bx) e^{-i\omega t} d\omega=\frac{1}{\pi}\Re\left[ \int_{0}^{\infty} \hat{G}_0(\omega,\bx_0;\bx) e^{-i\omega t} d\omega \right],
\end{equation}
where $\Re$ represents the real part,
and 
\begin{equation}\label{3.14.0}
    \hat{G}_0(\omega,\bx_0;\bx)=\int_{0}^{\infty} G_0(t,\bx_0;\bx) e^{i\omega t }d t, 
\end{equation}
where $G_0(t,\bx_0;\bx)|_{t<0}\equiv 0$.

The above observations inspire us to replace the time-domain Green's function with a truncated finite-bandwidth Fourier summation of the frequency-domain Green's function, where the singularities of the time-domain Green's function are naturally transferred to the singularities at source points of the frequency-domain Green's function, resulting in the time-frequency-time (TFT) method for Case 1 and the time-frequency-time-frequency (TFTF) method for Case 2, as alluded to in Section \ref{khform}.

\subsection{TFT Eulerian Hadamard integrator}
Now we start with developing the TFT method for Case 1, where $f\equiv 0$, and $u^1$ and $u^2$ are compactly supported, highly oscillatory smooth functions, and we will compute the integrals of $G_0$, $\dot{G}_0$ and $\ddot{G}_0$, respectively, where, by the integration or integral of a Green's function, we mean the integral of the product of a Green's function and a sufficiently well-behaved function. 

First, we apply the time-domain Green's function to $\phi(\bx)$, a test function in the space $K$ as defined above, and use \rf{3.13.0} to obtain 
\begin{equation}\label{3.19.2}
\begin{aligned}
\int_{V_{\boldsymbol{x}_0}}G_0(t,\bx_0;\bx)\phi(\bx) d\boldsymbol{x}&=\int_{V_{\boldsymbol{x}_0}}\frac{1}{\pi}\int_{-\infty}^{\infty} e^{-i\omega t} \hat{G}_0(\omega,\bx_0;\bx) d \omega\phi(\bx) d\boldsymbol{x}\\
&=\frac{1}{2\pi}\int_{-\infty}^{\infty}\int_{V_{\boldsymbol{x}_0}}  \hat{G}_0(\omega,\bx_0;\bx) \phi(\bx) d\boldsymbol{x} e^{-i\omega t} d \omega\\
&= \frac{1}{\pi}\Re\left[\int_{0}^{\infty}\int_{V_{\boldsymbol{x}_0}}  \hat{G}_0(\omega,\bx_0;\bx) \phi(\bx) d\boldsymbol{x} e^{-i\omega t} d \omega\right].\\
\end{aligned}
\end{equation}   
Here $V_{\boldsymbol{x}_0}$ is an integration region independent of time $t$, satisfying 
\begin{equation}\label{3.8}
      V_{\boldsymbol{x}_0}\supseteqq \{\boldsymbol{x}:\tau< t_0\}  \supseteqq \{\boldsymbol{x}:\tau< \Delta T\},
    \end{equation}
  where $t_0$ is a constant that depends on the medium and satisfies $t<\Delta T<t_0<\bar{T}(\bx_0)$.    

\color{g}
To apply the time derivatives of a Green's function in the frequency domain, we follow \cite{weicheburqia24} to introduce the geodesic (ray) polar transformation for a given source $\bx_0$:
\begin{equation}\label{gec}
    P[\bx_0]: \bx \rightarrow (\tau,\boldsymbol{\omega}),
\end{equation}
where $\tau$ is the traveltime, $\boldsymbol{\omega}\in \mathbb{S}^{m-1}$ is the take-off angle of the ray, and $(\tau,\boldsymbol{\omega})$ are the geodesic polar coordinates. Within any neighborhood of $\bx_0$ not containing any caustics other than $\bx_0$, there is one and only one ray connecting $\bx$ and $\bx_0$, which means $P[\bx_0]$ is well-defined and one-to-one, and $\boldsymbol{x}$ is a smooth function of the point $\boldsymbol{y}=\bx_0+\tau \boldsymbol{\omega}$. To facilitate our following discussions, we also write down the volume element as the following,
\begin{equation}
    \mathrm{d} \bx=\mathrm{d} s \mathrm{d} S=c \mathrm{d} \tau \mathrm{d} S=c \left|\frac{\partial S}{\partial \boldsymbol{\omega}} \right|\mathrm{d} \tau \mathrm{d}\boldsymbol{\omega},
\end{equation}
where $s$ is the arc length along the ray, $ \mathrm{d} S=\left|\frac{\partial S}{\partial \boldsymbol{\omega}}\right| \mathrm{d} \boldsymbol{\omega}$ is the element of area cut out on the wave front $\tau=$ const. by rays emanating from the solid angle element $\mathrm{d} \omega$ at the source. According to the appendix of \cite{weicheburqia24}, we have
\begin{equation}\label{2.68}
  \left|\frac{\partial S}{\partial \omega}\right|=\frac{\tau^{m-1}}{4 \rho_0 c_0^m \pi^{m-1} \rho c v_0^2},
\end{equation}
where $\rho_0=\rho(\bx_0)$ and $c_0=c(\bx_0).$

We now apply $\dot{G}(t,\bx_0;\bx)$ to a test function $\phi\in K$ and use the geodesic polar transformation and its inverse,
\begin{equation}\label{3.32.0}
    \begin{aligned}
        \int \dot{G}(t,\bx_0;\bx) \phi(\bx) d\bx &= 2t\int f^{-\frac{1}{2}-1}_+[t^2-\tau^2]v_0(\bx)\phi(\bx) d\bx\\
         &=2 t \int\left(-\frac{\partial}{2 \tau \partial \tau}\right) f_{+}^{-\frac{1}{2}}\left[t^2-\tau^2\right] v_0 \phi  d \bx \\
        &=2 t \int_0^{2\pi} \int_0^t \left(-\frac{\partial}{2 \tau \partial \tau}\right) f_{+}^{-\frac{1}{2}}\left[t^2-\tau^2\right] v_0 \phi \left|\frac{\partial S}{\partial \boldsymbol{\omega}}\right| c \mathrm{d} \tau \mathrm{d} \theta\\
        &=t\int_{0}^{2\pi}\left[f_+^{-\frac{1}{2}}[t^2-\tau^2]\frac{\phi}{4 \rho_0 c_0^2 \pi \rho  v_0}\right]_{\tau=t}^0 \mathrm{~d} \theta \\
&\quad +t \int_{0}^{2\pi} \int_0^tf_+^{-\frac{1}{2}}[t^2-\tau^2]\left(\frac{\partial}{\partial \tau}\right)\left[\frac{\phi}{4 \rho_0 c_0^2 \pi \rho  v_0}\right] \mathrm{d} \tau \mathrm{d} \theta\\
&=\frac{\phi(\bx_0)}{\rho(\bx_0)} +t \int_{0}^{2\pi} \int_0^t G_0(t,\bx_0;\bx) \frac{1}{v_0} \left(\frac{\partial}{\partial \tau}\right)\left[\frac{\phi}{4 \rho_0 c_0^2 \pi \rho  v_0}\right] \mathrm{d} \tau \mathrm{d} \theta\\
&=\frac{\phi(\bx_0)}{\rho(\bx_0)} +t \int_{V_{\boldsymbol{x}_0}} G_0(t,\bx_0;\bx) \frac{1}{v_0} \left(\frac{\partial}{\partial \tau}\right)\left[\frac{\phi}{4 \rho_0 c_0^2 \pi \rho  v_0}\right] \frac{1}{c|\frac{\partial S}{\partial \boldsymbol{\omega}}|}\mathrm{d} \bx\\
&=\frac{\phi(\bx_0)}{\rho(\bx_0)} +t \int_{V_{\boldsymbol{x}_0}} G_0(t,\bx_0;\bx) \left(\frac{\partial}{\partial \tau}\right)\left[\frac{\phi}{\rho  v_0}\right] \frac{\rho}{\tau}\mathrm{d} \bx,
    \end{aligned}
\end{equation} 
where $\theta$ is the take-off angle $\boldsymbol{\omega}$ when $m=2$, and we have used \rf{2.68}, \rf{3.9.0} and $f^{\lambda}_+[0]=0$. Further using \rf{3.19.2}, we have 
\begin{equation}\label{3.2.1}
  \begin{aligned}
 &\int \dot{G}_0(t, \bx_0;\bx) \phi(\bx) \mathrm{d} \boldsymbol{x}
= \frac{\phi(\bx_0)}{\rho(\bx_0)} +\frac{t}{\pi} \Re\left[\int_{0}^{\infty}e^{-i\omega t}\int_{V_{\boldsymbol{x}_0}} \hat{G}_0(\omega,\bx_0;\bx) \left(\frac{\partial }{\partial \tau}\right)\left[\frac{\phi}{\rho v_0}\right]\frac{\rho v_0}{\tau}\mathrm{d} \boldsymbol{x}\mathrm{d}\omega\right].\\
\end{aligned}
\end{equation}
Here the directional derivative
\begin{equation}
    \frac{\partial }{\partial \tau} =c^2\nab \tau \cdot \nab,
\end{equation}
is obtained by utilizing the method of characteristics to solve the eikonal equation \rf{3.7.0}.

As for $\ddot{G}_0$, for $t>0$, we have  
\begin{equation}
    \rho \ddot{G}_0(t,\bx_0;\bx)-\nab \cdot(\nu \nab G_0(t,\bx_0;\bx))=0.
\end{equation}
Applying $\ddot{G}_0(t,\bx_0;\bx)$ to a test function $\phi(\bx)$ and using integration by parts and \rf{3.19.2}, we have
\begin{equation}
\label{3.2.2}
\begin{aligned}
    \int_{V_{\boldsymbol{x}_0}} \ddot{G}_0(t,\bx_0;\bx)\phi(\boldsymbol{x})\mathrm{d} \boldsymbol{x}&= \int_{V_{\boldsymbol{x}_0}} \nab\cdot(\nu \nab G_0(t,\bx_0;\bx)) \frac{\phi}{\rho} d\bx\\
    &=\int G_0(t,\bx_0;\bx)\nab\cdot(\nu\nab \frac{\phi}{\rho}) \mathrm{d}\boldsymbol{x}\\
    &= \frac{1}{\pi}\Re\left[\int_{0}^{\infty}\int_{V_{\boldsymbol{x}_0}}  \hat{G}_0(\omega,\bx_0;\bx) \nab\cdot(\nu\nab \frac{\phi}{\rho}) d\boldsymbol{x} e^{-i\omega t} d \omega\right].
\end{aligned}
\end{equation}
Now, taking $\phi(\bx)$ as $\rho(\bx) u^1(\bx)$ and $\rho(\bx) u^2(\bx)$, respectively, in the Kirchhoff-Huygens representation formulas \eqref{2.5.0}-\eqref{2.6.0}, we obtain the following Eulerian formulations,
\begin{equation}\label{3.44.0}
    \begin{aligned}
        u(t,\bx_0)&= u^1(\bx_0)+\frac{1}{\pi}\Re\left[\int_{0}^{\infty}\int_{V_{\boldsymbol{x}_0}}  \hat{G}_0(\omega,\bx_0;\bx)\rho u^2(\boldsymbol{x})  d\boldsymbol{x}e^{-i\omega t} d \omega\right]\\
&\quad +\frac{t}{\pi} \Re\left[\int_{0}^{\infty}\int_{V_{\boldsymbol{x}_0}}\hat{G}_0(\omega,\bx_0;\bx)(\frac{\partial }{\partial \tau})(\frac{u^1}{ v_0})\frac{\rho v_0}{\tau} \mathrm{d} \boldsymbol{x}e^{-i\omega t}\mathrm{d}\omega\right],\\
    \end{aligned}
\end{equation}
and 
\begin{equation}\label{3.45.0}
    \begin{aligned}
        u_t(t,\bx_0)&= u^2(\bx_0)+\frac{t}{\pi} \Re\left[\int_{0}^{\infty}\int_{V_{\boldsymbol{x}_0}}\hat{G}_0(\omega,\bx_0;\bx)(\frac{\partial }{\partial \tau})(\frac{u^2}{ v_0})\frac{\rho v_0}{\tau}\mathrm{d} \boldsymbol{x}e^{-i\omega t}\mathrm{d}\omega\right]\\
&\quad +\frac{1}{\pi} \Re\left[\int_{0}^{\infty}\int_{V_{\boldsymbol{x}_0}} \hat{G}_0(\omega, \bx_0;\bx)\nab\cdot(\nu\nab u^1)\mathrm{d} \boldsymbol{x}e^{-i\omega t}\mathrm{d}\omega\right].\\
    \end{aligned}
\end{equation}
\color{black}

Since the frequency-domain Green's function is only singular at the source point which is amenable to Eulerian evaluation on uniform meshes, we refer to \rf{3.44.0}-\rf{3.45.0} as the Eulerian Hadamard-Kirchhoff-Huygens (HKH)-TFTF propagator, which is used to propagate the wave field from $\tau=0$ to $\tau=t$, where $0<t\leq \Delta T$. Although the HKH-TFTF  propagator is only valid for a short-time period in a caustic-free neighborhood, recursively applying this propagator in time yields the {\it TFT Eulerian  Hadamard integrator} to solve time-dependent wave equations globally in time, where caustics are treated implicitly and spatially overturning waves are handled naturally in {\it time domain}.

Here we remark in passing that we use the term ``propagator'' to indicate a short-time valid solution operator, while we use the term ``integrator'' to indicate a globally valid solution operator in either time-space or frequency-space domain.

\subsection{TFTF Eulerian Hadamard integrator}

Now we develop the TFTF method for Case 2, point-source problems for the Helmholtz equation, by first solving the corresponding time-domain wave equations and then applying the Fourier transform in time to the resulting time-domain solutions. 

According to Duhamel's principle, we can rewrite the corresponding time-domain wave equation of the point-source Helmholtz equation \rf{1.5} in the same form as Case 1, 
\beq{2.03}
\rho(\bx) u_{tt}(t,\bx) - \nab\cdot[\nu(\bx)\nab u(t,\bx)] = 0\;\;\mbox{ for }\;\; t>0\,,
\eeq
with initial conditions
\beq{2.04}
u(0,\bx) = 0\; \mbox{ and }\; u_t(0,\bx) = \fr{1}{\rho(\bz)}\,\dlt(\bx-\bz) \,.
\eeq
The challenge arises from the fact that the wave solution in Case 2 is the time-domain Green's function, a generalized function. In each time step, we need to numerically \color{b} apply \color{black} generalized functions to other generalized functions. Specifically, with the initial conditions at $T=0$ containing the singular generalized function $\delta$, we have 
\begin{equation}
    u(t,\bx_0)=\int \rho(\bx) G(t,\bx_0;\bx) \frac{\delta(\bx-\bz)}{\rho(\bz)} d\bx=G(t,\bx_0;\bz) 
\end{equation}
and 
\begin{equation}
    u_t(t,\bx_0)=\int \rho(\bx) G_t(t,\bx_0;\bx) \frac{\delta(\bx-\bz)}{\rho(\bz)} d\bx=G_t(t,\bx_0;\bz)
\end{equation}
in the \textit{distributional sense}. Since $G_t(t,\bx_0;\bz)$ is \textit{singular}, no locally integrable function is available that can be used to update $u_t(t, \bx_0)$ through {\it direct} integration. To overcome this challenge, instead of initializing at $T=0$ using the singular delta function, we initialize the wave field at $T=\Delta T$ with a smooth approximation, which is obtained by truncating the inverse Fourier transform of the frequency-domain Green's function. 

According to \rf{2.5.0} and \rf{2.6.0}, we approximate  the initial wave fields at time $t\in (0, \Delta T]$ as follows,
\begin{equation}
\label{3.2.4}
\begin{aligned}
    u(t,\boldsymbol{x})&\approx G(t,\bz;\bx)\approx \frac{1}{\pi}\Re\left[\int_0^{B}e^{-i\omega t}\hat{G}_0(\omega,\bz;\bx)d\omega\right],\\
\end{aligned}
\end{equation}
and 
\begin{equation}
\label{3.2.5}
\begin{aligned}
    u_t(t,\boldsymbol{x})&\approx \frac{1}{\pi}\Re\left[\int_0^{B}(-i\omega)e^{-i\omega t}\hat{G}_0(\omega,\bz;\bx)d\omega\right],\\
\end{aligned}
\end{equation}
where $B$ is an artificial bandwidth, leading to truncated, frequency-domain smooth approximations to the wave fields. Although here we have dropped the frequency information beyond the bandwidth $B$ of the time-domain solution, the Shannon-Nyquist Sampling theorem indicates that the remaining frequency content is sufficient for reconstructing the band-limited frequency-domain solution with angular frequency $\omega<\frac{B}{2}$ from the truncated wave solutions.

Now recursively applying the HKH-TFTF propagator starting from the initial time $t=\Delta T$ yields the global time-domain Green's function. Further we can utilize the Fourier transform to get the frequency-domain solution $\hat{u}(\omega,\bx)$ 
\begin{equation}\label{tftf}
\hat{u}(\omega,\bx)=\int_0^{\infty}e^{i\omega t}u(t,\bx)dt,
\end{equation}
leading to the {\it TFTF Eulerian  Hadamard integrator} to solve frequency-domain point-source wave equations (point-source Helmholtz equations) globally in time and thus globally in space, where the latter (``globally in space'') is implied by the eikonal-defined traveltime function.

Consequently, the TFTF Eulerian Hadamard integrator is able to treat caustics implicitly and handle spatially overturning waves naturally in {\it frequency domain}.


 
\section{Numerics for Eulerian Hadamard integrators}
\label{sec4}
We present numerics for implementing TFT and TFTF Eulerian Hadamard integrators. We will design Eulerian methods to numerically discretize the HKH\color{b}-TFTF \color{black} propagator so as to obtain the corresponding Eulerian Hadamard integrators by recursively applying the propagator in time. For highly oscillatory wave fields, we maintain a fixed number of points per wavelength (PPW) to uniformly discretize the computational domain into regular grid points. 

To start with, we briefly discuss higher-order sweeping schemes for the eikonal and transport equations, yielding the squared-phase function  $\tau^2$ and the Hadamard coefficient $v_0$, respectively. Subsequently, the low-rank representations based on Chebyshev interpolation are introduced to compress those smooth ingredients used in the propagator. Finally, we present preliminary algorithms for implementing TFT and TFTF Eulerain Hadamard integrators.

\subsection{Hadamard ingredients}
In order to form the HKH-TFTF propagator \eqref{3.44.0}-\eqref{3.45.0}, we need to compute some ingredients corresponding to $\tau$ and $v_0$, including
\begin{equation}\label{2.76}
  \tau, \; v_0,\; \nab\tau,\; \nab \tau \cdot \nab v_0.
\end{equation}
Since $\tau^2$ rather than $\tau$ is smooth at source points, we obtain smooth representations for the following quantities:
\begin{equation}\label{2.77}
  \tau^2,\; v_0,\; \nab\tau^2,\; \nab \tau^2 \cdot \nab v_0,
\end{equation}
which can be transformed back to \rf{2.76} via dividing by $\sqrt{\tau^2}$. We refer to the ingredients in \rf{2.77} as Hadamard ingredients.

$\tau^2$ and $v_0$ are obtained by solving the eikonal equation \rf{3.7.0} first and the transport equation \rf{SASWE3.01} afterwards. $\nab \tau^2$ is obtained by a proper finite-difference method. As for $\nab \tau^2 \cdot \nab v_0$, instead of using finite differencing which may reduce accuracy, we again use the transport equation \rf{SASWE3.01} to obtain
 \begin{equation}\label{2.79}
   \nab \tau^2 \cdot \nab v_0=\frac{2m\rho v_0-\nab\cdot\left[\nu \nab \tau^2\right]v_0}{2\nu}.
 \end{equation}
 
 \subsection{Numerical schemes}
 The leading-order term of the Hadamard's ansatz is defined by two functions, the eikonal $\tau$ satisfying the eikonal equation \rf{3.7.0} and the Hadamard coefficient $v_0$ satisfying the transport equation \rf{SASWE3.01}. Since we have assumed that the Hadamard's ansatz is valid locally around the source point, we need access to these two functions in order to construct the ansatz. Since the eikonal equation as a first-order nonlinear partial differential equation does not have analytical solutions in general, we have to use a robust, high-order numerical scheme to numerically solve this equation; moreover, the eikonal equation equipped with a point-source condition is even more tricky to deal with due to the upwind singularity at the source point \cite{qiasym02adapt}. To make the situation even more complicated, the transport equation \rf{SASWE3.01} for the Hadamard coefficient $v_0$ is weakly coupled with the eikonal equation \rf{3.7.0} in that the coefficients of the former equation depend on the solution of the latter.

Fortunately, this set of weakly coupled equations with point-source conditions has been solved to high-order accuracy by using Lax-Friedrichs weighted essentially non-oscillatory (LxF-WENO) sweeping schemes as demonstrated in \cite{qiayualiuluobur16}. The high-order schemes in \cite{qiayualiuluobur16} have adopted essential ideas from many sources including  \cite{liuoshcha94,oshshu91,jiapen00,qiasym02adapt,kaooshqia04,zha05,zharechov05,fomluozha09,luoqiabur14,luqiabur16} and have been used in many applications. Consequently, we will adopt these schemes to our setting as well and we omit details here.

\subsection{Multivariate Low-rank representation} 
The Hadamard ingredients \rf{2.76} are time-independent. According to \cite{qiayualiuluobur16}, if variable functions $\rho$ and $\nu$ are analytic, then $\tau^2$ and $v_0$ are analytic in the
region of space containing a point source but no other caustics. Consider the source region $\Omega_S$ and the receiver region $\Omega_R$, with $\Omega_S \subset \Omega_R$, such that rays originating from any $\bx_0 \in \Omega_S$ will not intersect in $\Omega_R$, which means that no caustics will develop in $\Omega_R$; consequently, $\tau^2(\boldsymbol{x}_0;\bx)$ and $v_0(\boldsymbol{x}_0;\bx)$ are analytic in $\Omega_S\times \Omega_R$ so that we can construct Chebyshev-polynomials based low-rank representations of Hadamard ingredients with respect to both $\boldsymbol{x}_0$ and $\boldsymbol{x}$ \cite{luqiabur16,luqiabur16b,luqiabur18,liusonburqia23}. 

Letting $f(x,y,x_0,y_0)=\tau^2,\; v_0,\;\frac{\partial \tau^2}{\partial x},\; \frac{\partial \tau^2}{\partial y},\;\nab\tau^2\cdot \nab v_0$, respectively, we consider the following analytical low-rank representation using Chebyshev interpolation, 
\begin{equation} \label{3.33}f(x,y,x_0,y_0)\approx\sum_{i=1}^{n_i}\sum_{j=1}^{n_j}\sum_{k=1}^{n_k}\sum_{l=1}^{n_l}F(i,j,k,l)\bar{T}_i(x)\bar{T}_j(y)\bar{T}_k(x_0)\bar{T}_l(y_0),
\end{equation}
where $n_i$, $n_j$, $n_k$, and $n_l$ are the orders of Chebyshev interpolation, $F$ is a 4-D tensor of size $n_i\times n_j\times n_k\times n_l$ which contains the spectral coefficients to be determined, and $\bar{T}_i$, $\bar{T}_j$, $\bar{T}_k$, and $\bar{T}_l$ represent the Chebyshev interpolants defined via translating the standard Chebyshev polynomials $T_m$ defined on $[-1,1]$ to the corresponding domain, 
\begin{equation}\label{ray3.34}
    T_m(s)=\cos(m\arccos(s)),\quad\quad s\in [-1,1].
\end{equation}
We can obtain $F$ by applying fast cosine transforms with respect to $x$, $y$, $x_0$, and $y_0$ to the tensor $f(x_i^c,y_j^c,x_{0,k}^c,y_{0,l}^c)$, where $x_i^c$, $y_j^c$, $x_{0,k}^c$, and $y_{0,l}^c$ are $n_i$, $n_j$, $n_k$, and $n_l$-order Chebyshev nodes in $\Omega_R\times \Omega_S$, respectively, which are also obtained by translating the $n$-order Chebyshev nodes $\{s_{m}\}$ in $[-1,1]$,
\begin{equation}
    s_{m}=\cos\left(\frac{2m-1}{2 n}\right),\quad \quad m=1,2, \cdots, n.
\end{equation}
To calculate $f(x_i^c,y_j^c,x_{0,k}^c,y_{0,l}^c)$, we adopt a computational strategy as used in \cite{liusonburqia23}: first choose a region $\Omega_c$ which is slightly larger than $\Omega_R$ and on which the ingredients are still analytic, then compute numerically the ingredients in $\Omega_c$ with sources located at $(x_{0,k}^c,y_{0,l}^c)$, respectively, and finally use cubic spline interpolation of the just computed ingredients on uniform grids to obtain $f(x_i^c,y_j^c,x_{0,k}^c,y_{0,l}^c)$. Here $\Omega_c$ is introduced to ensure the accuracy of the interpolations near the boundary of $\Omega_R$.

After obtaining $F$, we can evaluate $f(x,y,x_0,y_0)$ at any point $(x,y,x_0,y_0)\in \Omega_R\times \Omega_S$ using \rf{3.33}, 
and the evaluation can be further accelerated by partial summation \cite{luoqiabur14a}. So far, we have finished the pre-computation and compression of the Hadamard ingredients in $\Omega_R\times \Omega_S$ so that the integral kernels can be constructed by assembling these ingredients into the Hadamard's or Babich's ansatz.

Different from marching along a spatial direction \cite{luqiabur16}, when we use the Hadamard integrator to march in time direction, we always take $\Omega_S$ as a subset of $\Omega_R$. Actually, for given caustic-free regions $\Omega_S$ and $\Omega_R$, the maximum time step is bounded by the minimum traveltime
\begin{equation}
    \tau_{\min}=\min_{\bx_0\in \partial \Omega_S,\;\bx\in \partial \Omega_R} \tau(\bx_0;\bx).
\end{equation}
We do not need to obtain the strict upper bound. Instead, we can determine the time step by dividing the minimum distance between $\partial \Omega_S$ and $\partial \Omega_R$ by the maximum wave speed $\max c$. Due to possible caustics in an inhomogeneous medium, \color{b} we need to partition the computational domain $\Omega$ into $N_I$ source regions $\{\Omega_{S}^{\ell}, \;\ell=1,2,\cdots N_I\}$ and designate a corresponding receiver region  $\Omega^{\ell}_R$ for each source region $\Omega^{\ell}_S$, such that we can utilize the locally valid short-time HKH propagator to update the wave fields in $\Omega^{\ell}_S$ using the information in $\Omega_R^{\ell}$, separately. Since the waves will not propagate to the boundary of the computational domain, the wave field in the parts of the receiver regions that lie outside the computational domain is naturally set to zero. \color{black}  The partitioning of regions and the time step only depend on the medium.

\color{b}
With the low-rank approximations, we are able to evaluate the Hadamard ingredients with $O(1)$ computational complexity, which is crucial for the numerical algorithm that follows. Additionally, we will demonstrate that the computational cost of constructing these low-rank approximations is independent of the oscillatory wave field. 

Since both $\tau^2$ and $v_0$ are non-oscillatory analytic functions in a caustic-free region, this means that we can solve the eikonal and transport equations with high-order accuracy on a coarser grid (compared to the sampling grid for oscillatory wave fields); on the other hand, the orders of the multivariable Chebyshev interpolations are chosen according to the accuracy requirement in compressing the Hadamard ingredients, which are in turn medium-dependent. For instance, in the numerical examples that we are going to show, we set $n_i=n_j=13$ and $n_k=n_l=11$, which do not change with the oscillatory nature of the wave field.
\color{black}
 \subsection{Numerical discretization of HKH-TFTF propagator}
To numerically implement the HKH-TFTF propagator, we shall first truncate the infinite frequency-domain integration interval $[0,\infty)$ with an artificial frequency-bandwidth $B$ that is problem-specific in the sense that the bandwidth depends on the highly-oscillatory initial conditions in Case 1 and on the angular frequency $\omega$ in Case 2.

In the TFT method, once the oscillatory, smooth initial conditions \rf{1.2} are given, the wave solution is essentially bandlimited in the sense that the frequency content of the corresponding wave solution is determined and decays rapidly outside a certain frequency band. Thus we can select a sufficiently large artificial frequency-domain bandwidth $B$ to truncate the Fourier integral in frequency, where the fast decaying property ensures that truncation does not introduce significant truncation errors.

\color{g}
In the TFTF method, the bandwidth $B=\kappa \omega$ is selected according to the angular frequency $\omega$, where $\kappa$ is an  oversampling parameter. 

Using asymptotic properties of Hankel functions \cite{abrste65}, for $m=2$, we have the following asymptotic approximation for the frequency-domain Green's function, 
\begin{equation}\label{3.16.0}
  \hat{G}_0(\omega,\bx_0;\bx)\sim i\frac{\sqrt{\pi}}{2} v_0\sqrt{\frac{2}{\pi  \omega\tau}}e^{i\omega\tau-\frac{\pi}{4}} \quad\text{as}\quad \omega\tau\rightarrow \infty. 
\end{equation}
Hence, truncation errors may occur when we reconstruct the wave solution in Case 2 by 
\begin{equation}\label{4.8.1}
    G_0(t,\bx_0;\bx)=\Re[\int_0^{\infty}e^{-i\omega t} \hat{G}_0(\omega,\bx_0;\bx) d \omega]\approx \Re[\int_0^{B}e^{-i\omega t} \hat{G}_0(\omega,\bx_0;\bx) d \omega].
\end{equation}
For fixed $t$ and $\bx\in \{\bx:\tau(\bx_0;\bx)\neq t\}$, $e^{-i\omega t} \hat{G}_0(\omega,\bx_0;\bx)$ is oscillatory and decays at $O(\frac{1}{\sqrt{\omega}})$. Therefore, truncation of \rf{4.8.1} to bandwidth $B$ will introduce an error of no more than $O(\frac{1}{\sqrt{B}})$.

After truncating the infinite frequency domain to a bounded interval, we uniformly discretize the frequency interval $[0,B]$ into 
\begin{equation}
    \{\omega_k: \omega_k=(k+\frac{1}{2})\Delta \omega, k=0,1,2,\cdots,N_{\omega}-1 \},
\end{equation}
where $B=N_{\omega} \Delta \omega$. Since the frequency-domain Green's function $\hat{G}_0(\omega,\bx_0;\bx)$ in the two-dimensional case has logarithmic singularity when $\omega\tau \rightarrow 0$, 
\begin{equation}\label{3.15.0}
     \hat{G}_0(\omega,\bx_0;\bx)\sim \frac{\sqrt{\pi}}{2} v_0 \frac{2}{\pi}\ln(\omega\tau) \quad\text{as}\quad\omega\tau\rightarrow 0,
\end{equation}
 we need to treat the integrals of $\hat{G}_0(\omega,\bx_0;\bx)$ over $[0,\Delta \omega]$ with care.

Now we discretize the HKH-TFTF propagator \rf{3.44.0}-\rf{3.45.0} in the frequency domain as follows,
\color{g}
\begin{equation}\label{4.14.1}
    \begin{aligned}
        u(t,\bx_0)&\approx u^1(\bx_0)+\frac{1}{\pi}\Re\left[\int_{0}^{\Delta \omega}\int_{V_{\boldsymbol{x}_0}}  \hat{G}_0(\omega,\bx_0;\bx)\rho u^2(\boldsymbol{x})  d\boldsymbol{x}e^{-i\omega t} d \omega\right]\\
&\quad +\frac{\Delta\omega}{\pi}\Re \left[ \sum_{k=1}^{N_{\omega}-1}e^{-i\omega_k t}\int_{V_{\boldsymbol{x}_0}} \hat{G}_0(\omega_k,\bx_0;\bx)\rho u^2(\boldsymbol{x}) d\boldsymbol{x}\right]\\
&\quad +\frac{t}{\pi} \Re\left[\int_{0}^{\Delta \omega}\int_{V_{\boldsymbol{x}_0}}\hat{G}_0(\omega,\bx_0;\bx)(\frac{\partial }{\partial \tau})(\frac{u^1}{ v_0})\frac{\rho v_0}{\tau} \mathrm{d} \boldsymbol{x}e^{-i\omega t}\mathrm{d}\omega\right]\\
&\quad+\frac{t\Delta\omega}{\pi}\Re\left[\sum_{k=1}^{N_{\omega}-1}e^{-i\omega_k t} \int_{V_{\boldsymbol{x}_0}} \hat{G}_0(\omega_k,\bx_0;\bx)(\frac{\partial }{\partial \tau})(\frac{u^1}{ v_0})\frac{\rho v_0}{\tau}\mathrm{d} \boldsymbol{x} \right],\\
    \end{aligned}
\end{equation}
and 
\begin{equation}\label{4.15.1}
    \begin{aligned}
        u_t(t,\bx_0)&\approx u^2(\bx_0)+\frac{t}{\pi} \Re\left[\int_{0}^{\Delta \omega}\int_{V_{\boldsymbol{x}_0}}\hat{G}_0(\omega,\bx_0;\bx)(\frac{\partial }{\partial \tau})(\frac{u^2}{ v_0})\frac{\rho v_0}{\tau}\mathrm{d} \boldsymbol{x}e^{-i\omega t}\mathrm{d}\omega\right]\\
&\quad+\frac{t\Delta\omega}{\pi}\Re\left[\sum_{k=1}^{N_{\omega}-1}e^{-i\omega_k t} \int_{V_{\boldsymbol{x}_0}} \hat{G}_0(\omega_k,\bx_0;\bx)(\frac{\partial }{\partial \tau})(\frac{u^2}{ v_0})\frac{\rho v_0}{\tau}\mathrm{d} \boldsymbol{x} \right]\\
&\quad +\frac{1}{\pi} \Re\left[\int_{0}^{\Delta \omega}\int_{V_{\boldsymbol{x}_0}} \hat{G}_0(\omega, \bx_0;\bx)\nab\cdot(\nu\nab u^1)\mathrm{d} \boldsymbol{x}e^{-i\omega t}\mathrm{d}\omega\right]\\
&\quad+\frac{\Delta\omega}{\pi}\Re\left[\sum_{k=1}^{N_{\omega}-1}e^{-i\omega_k t} \int_{V_{\boldsymbol{x}_0}} \hat{G}_0(\omega_k, \bx_0;\bx)\nab\cdot(\nu\nab u^1)\mathrm{d} \boldsymbol{x} \right].
    \end{aligned}
\end{equation}

Before we further discretize the propagator in space, we need to deal with the frequency-space integral
\begin{equation}\label{4.12.2}
    \int_{V_{\boldsymbol{x}_0}} \int_0^{\Delta \omega} e^{-i\omega t}\; \hat{G}_0(\omega,\bx_0;\bx)\;d\omega\;\psi(\bx_0;\bx)\; d\bx,
\end{equation}
where $\psi(\bx_0;\bx)$ represents the frequency-independent part of the integrand in the propagator. The difficulty comes from the singularity of the frequency-domain Green's function when $\omega\rightarrow 0$. However, since $\tau$ has a small upper bound that only depends on the medium, we have $\omega\tau \sim 0$ uniformly; we thus use the polynomial approximation of the Hankel function $H_0^1(\omega\tau)$ \cite{abrste65} to approximate the kernel $\hat{G}_0$ so that 
\begin{equation}
       I(t,\bx_0;\bx)\doteq \int_0^{\Delta \omega} e^{-i\omega t} \hat{G}_0(\omega,\bx_0;\bx) d\omega 
\end{equation}
can be evaluated, and the technical details can be found in  Appendix \ref{A1}. Then integral \rf{4.12.2} can be rewritten using the time-dependent kernel $I(t,\bx_0;\bx)$
\begin{equation}
     \int_{V_{\boldsymbol{x}_0}} \int_0^{\Delta \omega} e^{-i\omega t}\; \hat{G}_0(\omega,\bx_0;\bx)\;d\omega\;\psi(\bx_0;\bx)\; d\bx=\int_{V_{\boldsymbol{x}_0}}  I(t,\bx_0;\bx) \psi(\bx_0;\bx) d\bx.
\end{equation}

We now discretize the HKH-TFTF propagator in space. Consider a source region $\Omega_S$ and the corresponding receiver region $\Omega_R$ which are uniformly discretized into $N$ and $M$ grid points, respectively. Then, a key issue is how to evaluate the self-interaction (diagonal) terms of the kernels. In \eqref{4.14.1} and \eqref{4.15.1}, we have the following singular diagonal terms
\begin{eqnarray}
&&   \hat{G}_0(\omega,\bx_0;\bx)|_{\bx=\bx_0}, \label{dig1}  \\
&&     I(t,\bx_0;\bx)|_{\bx=\bx_0}, \label{dig2}\\
&&   \left.\left[\hat{G}_0(\omega,\bx_0;\bx)\frac{\partial }{\tau \partial\tau}\left(\frac{u^k}{v_0}\right)\right]\right|_{\bx=\bx_0},  \quad k=1, 2,\label{dig3}\\
&&   \left.\left[\int_0^{\Delta \omega} e^{-i\omega t}\hat{G}_0(\omega,\bx_0;\bx)\frac{\partial }{\tau \partial\tau}\left(\frac{u^k}{v_0}\right) d\omega\right]\right|_{\bx=\bx_0}, \quad k=1, 2.\label{dig4}
\end{eqnarray}
After uniformly discretizing the computational regions, the diagonal terms can be replaced by integral averages over a cell centered at $\bx_0$, where the non-singular parts of the integrand are treated as constants. Letting $c_j$ be the cell of size $h$ centered at $\bx_0$, the Appendix of \cite{liusonburqia23} provides 
\begin{equation}\label{4.11.1}
\begin{aligned}
     \hat{G}_0(\omega,\bx_0;\bx_0)&\doteq \frac{1}{h^2}\int_{c_{j}} \hat{G}_0(\omega,\bx_0;\bx) d \bx\\
     &\approx\frac{i}{\left(2h n_0 \omega\right)^2}\left[8 \int_0^{\frac{\pi}{4}} \frac{h n_0 \omega}{2 \cos \theta} H_1^{(1)}\left(\frac{h n_0 \omega}{2 \cos \theta}\right) d \theta+4 i\right],
\end{aligned}
\end{equation}
where $n_0=n(\bx_0)$. Based on \rf{4.11.1}, we derive the approximation for
\begin{equation}
     I(t,\bx_0;\bx_0)\doteq \frac{1}{h^2} \int_{c_j} \int_0^{\Delta \omega} e^{-i\omega t} \hat{G}_0(\omega,\bx_0;\bx) d\omega d\bx 
\end{equation}
in Appendix \ref{A2}, where we have used the polynomial approximation of the Hankel function $H_1^{(1)}$. We further show in Appendix \ref{A3} that 
\begin{equation}\label{4.22.2}
\begin{aligned}
     &\quad\left. \left[\hat{G}_0(\omega,\bx_0;\bx)\frac{\partial }{\tau \partial\tau}\left(\frac{u^k}{v_0}\right)\right]\right|_{\bx=\bx_0}\\
   &\doteq \frac{1}{h^2}\int_{c_j} \left[\hat{G}_0(\omega,\bx_0;\bx)\frac{\partial }{\tau \partial\tau}\left(\frac{u^k}{v_0}\right)\right] d\bx=0, \quad k=1, 2.
\end{aligned}
\end{equation}
Multiplying \rf{4.22.2} by $e^{-i\omega t}$, we integrate the resulting formula over $[0,\Delta \omega]$ to compute the last diagonal term, leading to 
\begin{equation}
\begin{aligned}
    &\quad\left.\left[\int_0^{\Delta \omega} e^{-i\omega t}\hat{G}_0(\omega,\bx_0;\bx)\frac{\partial }{\tau \partial\tau}\left(\frac{u^k}{v_0}\right) d\omega\right]\right|_{\bx=\bx_0}\\
    &\doteq \frac{1}{h^2}\int_{c_j} \int_0^{\Delta \omega} e^{-i\omega t}\hat{G}_0(\omega,\bx_0;\bx)  d\omega\frac{\partial }{\tau \partial\tau}\left(\frac{u^k}{v_0}\right) d\bx=0, \quad k=1, 2.
\end{aligned}  
\end{equation}

At this point, we define the following frequency-dependent kernels:
\begin{eqnarray}
&&   U_1^{\omega}=[(U_1^{\omega})_{i,j}]=\left[\hat{G}_0(\omega,\bx_0^i;\bx^j)\right]_{1\leq i\leq N, 1\leq j\leq M},  \\
&&     U_2^{\omega}=[(U_2^{\omega})_{i,j}]=\left[ \frac{ \hat{G}_0(\omega,\bx_0^i;\bx^j)}{2\tau^2(\bx_0^i;\bx^j)}\frac{\partial \tau^2(\bx_0^i;\bx^j)}{\partial x}\right]_{1\leq i\leq N, 1\leq j\leq M},\\
&&  U_3^{\omega}=[(U_3^{\omega})_{i,j}]=\left[ \frac{ \hat{G}_0(\omega,\bx_0^i;\bx^j)}{2\tau^2(\bx_0^i;\bx^j)}\frac{\partial \tau^2(\bx_0^i;\bx^j)}{\partial y}\right]_{1\leq i\leq N, 1\leq j\leq M},\\
&&   U_4^{\omega}=[(U_4^{\omega})_{i,j}]=\left[ \frac{ \hat{G}_0(\omega,\bx_0^i;\bx^j)}{2v_0(\bx_0^i;\bx^j)\tau^2(\bx_0^i;\bx^j)}\nab \tau^2(\bx_0^i;\bx^j)\cdot \nab v_0(\bx_0^i;\bx^j)\right]_{1\leq i\leq N, 1\leq j\leq M},
\end{eqnarray}
and the corresponding time-dependent kernels 
\begin{equation}
    \check{U}^{t}_s=[(\check{U}^{t}_s)_{i,j}]=\left[\int_0^{\Delta \omega} e^{i\omega t}U_s^{\omega}(\bx_0^i;\bx^j) d\omega\right]_{1\leq i\leq N, 1\leq j\leq M},\;\; s=1,2,3,4.\;
\end{equation}

We further define vectors
\begin{equation}
    \begin{aligned}
        & f_1=[\rho(\bx^1)u^2(\bx^1),\quad\cdots,\quad\rho(\bx^M)u^2(\bx^M))]^T,\\
        & f_2=[c(\bx^1)u_x^1(\bx^1),\quad\cdots, \quad c(\bx^M)u_x^1(\bx^M)]^T,\\
        & f_3=[c(\bx^1)u_y^1(\bx^1),\quad\cdots,\quad c(\bx^M)u_y^1(\bx^M)]^T,\\
        & f_4=[\nu(\bx^1)u^1(\bx^1),\quad\cdots,\quad\nu(\bx^M)u^1(\bx^M)]^T,\\
        &  f_5=[\nab\cdot(\nu(\bx^1)\nab u^1(\bx^1),\quad\cdots,\quad\nab\cdot(\nu(\bx^M)\nab u^1(\bx^M)]^T,\\
        &   f_6=[c(\bx^1)u_{x}^2(\bx^1),\quad\cdots,\quad c(\bx^M)u_{x}^2(\bx^M)]^T,\\
        &f_7=[c(\bx^1)u_{y}^2(\bx^1),\quad\cdots,\quad c(\bx^M)u_{y}^2(\bx^M)]^T,\\
        &  f_8=[\nu(\bx^1)u^2(\bx^1),\quad\cdots,\quad\nu(\bx^M)u^2(\bx^M)]^T.
    \end{aligned}
\end{equation}
Then we introduce the following frequency-sampling functions of size $N\times 1$ for $1\leq k \leq N_{\omega}$
\begin{equation}
    \begin{aligned}
        & \hat{u}_1(\omega_k)=U_1^{\omega_k}f_1,\\
        &  \hat{u}_2(\omega_k)=U_2^{\omega_k}f_2+U_3^{\omega_k}f_3-U_4^{\omega_k}f_4,\\
        &  \hat{u}_3(\omega_k)=U_1^{\omega_k}f_5,\\
        &  \hat{u}_4(\omega_k)=U_2^{\omega_k}f_6+U_3^{\omega_k}f_7-U_4^{\omega_k}f_8,
    \end{aligned}
\end{equation}
and the time-sampling functions of size $N\times 1$ for $0<t_l\leq \Delta T$ 
\begin{equation}
    \begin{aligned}
        &\check{u}_1(t_{\ell})=\check{U}_1^{t_{\ell}}f_1,\\
&\check{u}_2(t_{\ell})=\check{U}_2^{t_{\ell}}f_2+\check{U}_3^{t_{\ell}}f_3-\check{U}_4^{t_{\ell}}f_4,\\
&   \check{u}_3(t_{\ell})=\check{U}_1^{t_{\ell}}f_5,\\
&   \check{u}_4(t_{\ell})=\check{U}_2^{t_{\ell}}f_6+\check{U}_3^{t_{\ell}}f_7-\check{U}_4^{t_{\ell}}f_8.
    \end{aligned}
\end{equation}

Once we obtain all the above sampling functions, according to \rf{3.44.0} and \rf{3.45.0}, we update $u(t_{\ell},\Omega_S)$ and $u_t(t_{\ell},\Omega_S)$ for all $0<t_{\ell}\leq \Delta T$  as
\begin{equation}\label{4.46}
      u(t_{\ell},\Omega_S)\approx \frac{1}{\pi}\Re\left[\check{u}_1(t_{\ell})+t_{\ell}\check{u}_2(t_{\ell})+\Delta \omega\sum_{k=1}^{N_{\omega}-1}e^{-i\omega_k t_{\ell}}\left(\hat{u}_1(\omega_k)+t_{\ell}\hat{u}_2(\omega_k)\right)\right]
\end{equation}
and
\begin{equation}\label{4.47}
      u_t(t_{\ell},\Omega_S)\approx \frac{1}{\pi}\Re\left[\check{u}_3(t_{\ell})+t_{\ell}\check{u}_5(t_{\ell})+\Delta \omega\sum_{k=1}^{N_{\omega}-1}e^{-i\omega_k t_{\ell}}\left(\hat{u}_3(\omega_k)+t_{\ell}\hat{u}_4(\omega_k)\right)\right].
\end{equation}
\color{black}

\subsection{Direct computation of the Hadamard integrators}
We first present a direct summation algorithm for the TFT Eulerian Hadamard integrator for Case 1, and based on this, we further present the TFTF Eulerian Hadamard integrator for Case 2.

\textbf{Algorithm 1} (TFT Eulerian Hadamard integrator for Case 1)
\begin{enumerate}
    \item  Uniformly discretize the computational domain $\Omega$ into a wave-resolution-satisfying uniform grid; partition $\Omega$ into $N_I$ non-overlapping source regions $\{\Omega_S^{\ell},\; \ell=1, 2, \cdots, N_I\}$ and designate the corresponding receiver regions to be $\{\Omega_R^{\ell},\; \ell=1, 2, \cdots, N_I\}$, set the initial time $T_0$, the large time step $\Delta T<\tau_{\min}$ for updating the initial condition, the small time step $\Delta t=\frac{\Delta T}{N_T}$ for wave resolution, and an ending time $T_{end}$ to ensure that the waves do not reach the computational boundary; determine the bandwidth $B$ according to the initial conditions; set the frequency step $\Delta \omega = \frac{B}{N_{\omega}}$, where $N_{\omega}=O(B)$, and generate $\omega_k=(k+\frac{1}{2})\Delta \omega$ for $k = 0, 1,2,\cdots N_{\omega}-1$; pre-compute and compress Hadamard ingredients \rf{2.77} in all sub-region pairs; precompute the diagonal terms for the kernel matrices.
    \item Initialize $u(T_0,\bx)=u^1(\bx)$ and  $u_t(T_0,\bx)=u^2(\bx)$; set time $T=T_0$ and the loop variable $p=0$.
    \item For the current time step with $T=T_0+p \Delta T$:
    \begin{enumerate}
        \item compute finite-difference derivatives of $u(T,\bx)$ and $u_t( T,\bx)$ and construct vectors $f_1$, $f_2$, $\cdots,$ $f_8$.
        \item determine a sub-region $\Omega_p$ of $\Omega$ that contains the region of influence of current data by extending each direction of the non-zero region of the solution $u(T,\bx)$ outward by $\max c \Delta T$; designate the source regions intersecting with $\Omega_p$ as the regions to be updated.
     \item  for all the source regions $\Omega^{\ell}_S$ to be updated:    
\begin{enumerate}
    \item\color{g} Interpolate the Hadamard ingredients via partial summation; construct the kernels 
    \begin{equation*}
        U_s^{\omega_k},\; s=1, 2, 3, 4,\; 1\leq k\leq N_\omega-1,
 \quad\text{and}\quad
        \check{U}_s^{n\Delta t},   s=1, 2, 3, 4, 1\leq n \leq N_T;
    \end{equation*}
  compute the matrix-vector multiplications to obtain the corresponding frequency- and time-sampling functions $\{\hat{u}_i(\omega_k)\}$ and $\{\check{u}_i(n\Delta t)\}$:
     \begin{equation*}
       \hat{u}_i(\omega_k),\;  i=1, 2, 3, 4, \; 1\leq k \leq N_{\omega}-1,  \quad\text{and}\quad
 \check{u}_i(n\Delta t),  i=1, 2, 3, 4, \; 1\leq n \leq N_{T}.
   \end{equation*}
    \item Update wavefields $u(T+n\Delta t,\Omega_S^{\ell})$ and $u_t(T+n\Delta t, \Omega_S^{\ell})$ for $n=1, 2 ,\cdots, N_T$ via \rf{4.46} and \rf{4.47}.  \color{black}
    \end{enumerate}
    \item Set to zero the wave fields within the source regions that do not require updating.
    
    \end{enumerate}
    \item Update $T=T_0+(p+1)\Delta T$. If $T<T_{end}$, then $p\leftarrow p+1$ and go to Step 3; else, stop.
\end{enumerate}

Now we consider the computational cost of Algorithm 1. For a given frequency bandwidth $B$, according to the Shannon-Nyquist Sampling theorem, we will have a spatial grid of size $O(B^2)$ in the two-dimensional space. To ensure that wave field is well resolved, we have $\Delta t=O(1/B)$, leading to $O(B^3)$ unknowns to compute. 

As mentioned before, $\Delta T=O(1)$ since it only depends on the medium, so is $N_I$. Thus we will perform Step 3 for $O(1)$ times, and in each iteration, there are $O(B)$ matrix-vector multiplications of size $O(B^2\times B^2)$, resulting in a computational complexity of at least $O(B^5)$. It is extremely expensive and impractical compared to the number of unknowns. To resolve this issue, we will utilize butterfly algorithms and their hierarchical extensions to accelerate matrix-vector multiplication.

\color{g}
As for Case 2, we first use \eqref{3.2.4} to approximate the time-domain solution in the interval $0< t \leq \Delta T$ and then apply \textbf{Algorithm 1} to solve the following wave equation in the interval $\Delta T<t \leq T_{end}$, 
\begin{equation}\label{4.62.0}
  \rho u_{tt}-\nab\cdot(\nu \nabla u)=0,\; \boldsymbol{x}\in \mathbb{R}^m,\; \Delta T<t \leq T_{end},
\end{equation}
with initial conditions
\begin{equation}\label{4.63.0}
    u(\Delta T,\bx)=\frac{1}{\pi}\Re\left[\sum_k e^{-i\omega_k \Delta T}\hat{G}_0(\omega_k,\bz,\bx)  \right],
\end{equation}
\begin{equation}\label{4.64.0}
    u_t(\Delta T,\bx)=\frac{1}{\pi}\Re\left[\sum_k (-i\omega_k) e^{-i\omega_k \Delta T}\hat{G}_0(\omega_k,\bz,\bx)  \right].
\end{equation}
Consequently, the frequency domain solution $\hat{u}(\omega,\Omega)$ is constructed by \rf{tftf}. This leads to the TFTF Eulerian Hadamard integrator.

\color{black}
\section{Fast computation of Hadamard integrators}
\label{sec5}
Here we propose a butterfly-based algorithm for the construction and application of each kernel, which is quasi-optimal in the sense that  the complexity and memory usage are reduced from $O(B^4)$ to $O(B^2\log^2(B))$. 

For large $\omega_k$, the discretized integral operators $U_i^{\omega_k}$ are full rank due to the oscillation of $\hat{G}(\omega_k,\boldsymbol{x}_0;\bx)$, but the judiciously selected submatrices of the discretized operators are low-rank compressible due to the so-called complementary low-rank property \cite{luqiabur16}. 
Consequently, we can use the butterfly algorithm to compress those frequency-dependent kernels. 

Even though the time-dependent kernels do not oscillate (which means that the complementary low-rank property is naturally satisfied), in this work we will use the same butterfly algorithm to compress them. Since the ranks of the submatrices of the time-dependent kernels are smaller than those of the frequency-dependent kernels, this fact naturally helps us reduce the computational complexity and memory usage required for compression and application. 

\color{b}
We have partitioned the computational region $\Omega$ into $N_I$ source regions $\{\Omega_S^{\ell},\; \ell = 1, 2, \cdots, N_I\}$ and designated a corresponding receiver region $\Omega_R^{\ell}$ for each source region $\Omega_S^{\ell}$. As shown in Fig. \ref{figure1}, we further partition each receiver region $\Omega_R^{\ell}$ into nine sub-regions, one of which overlaps with the corresponding source region $\Omega_S^{\ell}$ (the red square), while the remaining eight (the white squares) are adjacent to the corresponding source region.
\color{black}

\begin{figure}[htbp]
    \centering
    \includegraphics[width=0.35\linewidth]{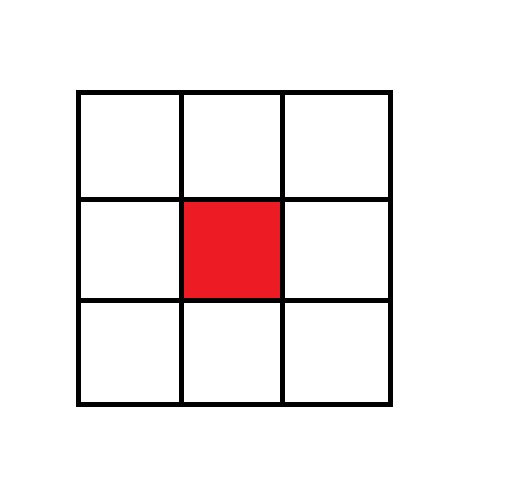}
    \caption{Illustration of $\Omega_S^{\ell}$ and $\Omega_R^{\ell}$: The red square represents the source region $\Omega_S^{\ell}$, and the largest square divided into nine parts represents the corresponding receiver region $\Omega_R^{\ell}$}
    \label{figure1}
\end{figure}

We first describe the interpolative-decomposition butterfly representation for the non-overlapping sub-regions and then introduce the hierarchical structure to handle the overlapping sub-region.

\subsection{Interpolative-decomposition butterfly (IDBF) representation}
For the integral kernel $K(\Omega_1,\Omega_2)$ defined on two non-overlapping regions $\Omega_1$ and $\Omega_2$,
we closely follow \cite{liyaneilkenyin15, lizha17, liuxinguomirerighyli21, liusonburqia23, candemyin09} to introduce the interpolative-decomposition butterfly algorithm. The algorithm first recursively subdivides the geometrical point sets $\Omega_1$ and $\Omega_2$ associated with the rows and columns of these operators into two subsets of approximately equal sizes using a k-dimensional tree clustering algorithm, until the number of elements in each subset is less than a preset value. When $\Omega_1$ and $\Omega_2$ have different sizes, the number of levels of the binary tree will depend on the set of a smaller size. This procedure generates two complete binary trees $\mathcal{T}_{1}$ and $\mathcal{T}_{2}$ of $L$ levels with root level $0$ and leaf level $L$.  In both trees, a non-leaf node $\tau$ at level $l$ has two children $\tau_1$ and $\tau_2$, where $\tau=\tau_1 \cup \tau_2$ and $\tau_1 \cap \tau_2=\emptyset$. For a non-root node $\tau$, its parent is denoted $p_\tau$.

The $L$-level butterfly representation of the kernel $K(\Omega_1,\Omega_2)$ requires the complementary low-rank property: at any level $0 \leq l \leq L$, for any node $\tau$ at level $l$ of $\mathcal{T}_1$ and any node $\nu$ at level $L-l$ of $\mathcal{T}_2$, the subblock $K(\tau, \nu)$ is numerically low-rank with rank $r_{\tau, \nu}$ bounded by a small number $r$ called the butterfly rank. As discussed in \cite{liusonburqia23}, even if there is a singularity in the integral kernel which may cause some subblocks to have non-constant ranks, the butterfly algorithm will still work when the singularity only exists at the interface of adjacent regions.

Thus, for any subblock $K(\tau, \nu)$, where $\tau\in \mathcal{T}_1$ and $\nu\in\mathcal{T}_2$, the complementary low-rank property permits a low-rank representation using interpolative decomposition (ID) as
\begin{equation}\label{bf1}
    K(\tau, \nu) \approx K(\tau, \bar{\nu}) V_{\tau, \nu}, 
\end{equation}
where the skeleton matrix $K(\tau, \bar{\nu})$ contains $r_{\tau, \nu}$ skeleton columns indexed by $\bar{\nu}$, and the interpolation matrix $V_{\tau, \nu}$ has bounded entries. The ID can be computed via, for instance, rank-revealing QR decomposition with a relative tolerance $tol$, which actually determines the accuracy of the overall butterfly representation. Here we briefly describe the so-called column-wise butterfly representation \cite{liyaneilkenyin15}\cite{liusonburqia23}. 

At level $l=0$, the interpolation matrices $V_{\tau, \nu}$ are explicitly formed, and at level $l>0$, they are represented in a nested fashion. To see this, consider a node pair $(\tau, \nu)$ at level $l>0$, and let $\nu_1$ and $\nu_2$, and $p_\tau$ be the children of $\nu$ and parent of $\tau$, respectively. From \rf{bf1}, we have
\begin{equation}\label{bf2}
\begin{aligned}
K(\tau, \nu)&=\left[K(\tau,\nu_1),K(\tau,\nu_2)\right]\\
     &\approx\left[K(\tau,\bar{\nu}_1),K(\tau,\bar{\nu}_2)\right]\left[\begin{array}{ll}
V_{p_\tau, \nu_1} & \\
& V_{p_\tau, \nu_2}
\end{array}\right]\\
&\approx K(\tau,\bar{\nu}) W_{\tau,\nu}\left[\begin{array}{ll}
V_{p_\tau, \nu_1} & \\
& V_{p_\tau, \nu_2}
\end{array}\right].
\end{aligned}
\end{equation}
Here $W_{\tau, \nu}$ and $\bar{\nu}$ are the interpolation matrix and skeleton columns from the ID of $\left[K\left(\tau, \bar{\nu}_1\right), K\left(\tau, \bar{\nu}_2\right)\right]$, respectively. This allows representing $V_{\tau, \nu}$ as
\begin{equation}
    V_{\tau, \nu}=W_{\tau, \nu}\left[\begin{array}{ll}
V_{p_\tau, \nu_1} & \\
& V_{p_\tau, \nu_2}
\end{array}\right],
\end{equation}
where we refer $W_{\tau,\nu}$ as the transfer matrix. We will select $O(r_{\tau,\nu})$ proxy rows $\hat{\tau}\subset \tau$ to approximately compute $V_{\tau,\nu}$ and $W_{\tau,\nu}$ via ID:
\begin{eqnarray}
&&    K(\hat{\tau}, \nu) \approx K(\hat{\tau}, \bar{\nu}) V_{\tau, \nu}, \quad l=0, \\
&&    \left[K\left(\hat{\tau}, \bar{\nu}_1\right), K\left(\hat{\tau}, \bar{\nu}_2\right)\right] \approx K(\hat{\tau}, \bar{\nu}) W_{\tau, \nu}, \quad 0<l<L.
\end{eqnarray}
And we do not select proxy rows at level $L$. 

In the current work, we select the proxy rows as follows.
Consider a subblock $(\tau,\nu)$ at level  $0<l<L$. Let $\bx^i$ be the geometrical centroid in the computational domain corresponding to the row index $i$.  Let $s_i$ be the index set consisting of the indices of neighboring centroids of the centroid $\bx^i$, where $i\in \bar{\nu}_1 \cup \bar{\nu}_2$; 
for example, $s_i$ may be the index set of all indices for centroids that are within a $10*\mathrm{h}$ distance of $\bx^i$, where $h$ is the mesh size. Let $f_\tau$ be the index set consisting of $\chi_1\left|\bar{\nu}_1 \cup \bar{\nu}_2\right|$ indices for uniformly selected centroids near the boundary of the subdomain $\tau$ with an oversampling parameter $\chi_1$, and let $g_{\tau}$ denote the index set consisting of $\chi_2\left|\bar{\nu}_1 \cup \bar{\nu}_2\right|$ indices for randomly selected centroids inside the subdomain $\tau$ with an oversampling parameter $\chi_2$, where $\left|\bar{\nu}_1\right|$ indicates the cardinality of the set $\bar{\nu}_1$. Accordingly, we choose the proxy rows as
\begin{equation}
    \hat{\tau}=\left(\cup_{i \in \bar{\nu}_1 \cup \bar{\nu}_2} s_i\right) \cap \tau \cup f_\tau\cup g_{\tau}.
\end{equation}

Numerically, we take the oversampling parameters $\chi_1=\chi_2=3$. \color{b} Here, the selection criteria for proxy rows and  oversampling parameters are empirical, determined by the properties of the kernel function and numerical accuracy tests. This ensures that the error introduced by the butterfly decomposition is significantly smaller than that from the numerical Hadamard integrator itself.
\color{black}
The proxy rows can greatly reduce the complexity. With all the interpolation and transfer matrices computed, the butterfly representation of $K(\Omega_1, \Omega_2)$ is:
\begin{equation}
    K(\Omega_1, \Omega_2) \approx K^L W^L W^{L-1} \ldots W^1 V^0 .
\end{equation}

Let $\nu_1$, $\nu_2$, $\ldots$, and $\nu_{2^{L-l}}$ denote the nodes at level $L-l$ of $\mathcal{T}_2$, and $\tau_1$, $\tau_2$, $\ldots$, and $\tau_{2^l}$ denote the nodes at level $l$ of $\mathcal{T}_1$. The interpolation factor $V^0$, the transfer factors $W^l$ for $l=1, \ldots, L$, and the skeleton factor $K^L$ are assembled in the following way, 
\begin{eqnarray}
   && V^0=\operatorname{diag}\left(V_{\tau, \nu_1}, \ldots, V_{\tau, \nu_{2^{L}}}\right), \quad\left(\tau, \nu_i\right) \text { at level } l=0,\label{bf8}\\
   && K^L=\operatorname{diag}\left(K\left(\tau_1, \bar{\nu}\right), \ldots, K\left(\tau_{2^L}, \bar{\nu}\right)\right), \quad\left(\tau_i, \nu\right) \text { at level } l=L,\label{bf9}\\
   && W^l=\operatorname{diag}\left(W_{\tau_1}, \ldots, W_{\tau_{2^{l-1}}}\right), \quad l=1, \ldots, L,\label{bf10} \\
   && W_{\tau_i}=\left[\begin{array}{c}
\operatorname{diag}\left(W_{\tau_1^1, \nu_1}, \ldots, W_{\tau_1^1, \nu_{2^{L-l}}}\right) \\
\operatorname{diag}\left(W_{\tau_i^2, \nu_1}, \ldots, W_{\tau_i^2, \nu_{2^{L-l}}}\right)
\end{array}\right], \quad (\tau_i^{\{1,2\}}, \nu_j) \text { at level } l,\label{bf11}
\end{eqnarray}
where $\tau_i^1$ and $\tau_i^2$ denote the children of $\tau_i$. 

\subsection{Hierarchically off-diagonal butterfly (HODBF) representation}
When $\Omega_1=\Omega_2$, the resulting integral kernels $K(\Omega_1,\Omega_1)$ have singular diagonal terms. We use the hierarchically off-diagonal butterfly (HODBF)\cite{liuguomir16}\cite{liusonburqia23} representation to resolve such singularities.

We still use the $L$-level binary tree $\mathcal{T}_1$ for $\Omega_1$ as in the last subsection. For any two siblings $\tau_1$ and $\tau_2$ with parent $\tau$ at level $l$ of $\mathcal{T}_1$, we can directly extract two $(L-l)$-level subtrees $\mathcal{T}_{\tau_1}$ and  $\mathcal{T}_{\tau_2}$ rooted at $\tau_1$ and $\tau_2$ from $\mathcal{T}_1$, respectively, due to the nested structure of the binary tree. Next, we construct the $(L-l)$-level butterfly representation for $K(\tau_1,\tau_2)$, where $\tau_1\cap \tau_2=\emptyset$. Finally, there are $2^{l}$ butterfly representations at each level $l=1,2,\cdots,L$. The subblocks $K(\tau,\tau)$ for $\tau$ at level $L$ are kept as the dense blocks. Such a representation is called the hierarchically off-diagonal butterfly (HODBF) representation.

\subsection{The complexity and memory usage}
We analyze the complexity and memory usage of IDBF and HODBF for $K=U_i^{\omega}$ for $i=1,2,3,4$ and $\check{U}^t_j$ for $ j=1,2,3,4$. Using the Chebyshev interpolation and partial summation, we can evaluate each element in those kernels with $O(1)$ complexity \cite{liusonburqia23}. 

For IDBF of $K(\Omega_1,\Omega_2)\in \mathbb{C}^{O(n)\times O(n)}$, we see that, in \rf{bf8}-\rf{bf11}, $V^0$ and $K^L$ contain $2^L$ diagonal blocks each with $O\left(r_{\tau,\nu}\right)$ non-zeros, and $W^l$ contains $2^L$ blocks $W_{\tau, \nu}$ each with $O\left(r_{\tau ,\nu}^2\right)$ non-zeros. After selecting proxy rows, we will evaluate $O(n)$ IDs of size $O(r_{\tau,\nu})\times O(r_{\tau,\nu})$ for $O(\log(n))$ levels, resulting in $O(n\log(n))$ complexity and memory usage if $O(r_{\tau,\nu})=O(1)$, analogous to the classic result of the  butterfly algorithm.

However, for the kernel with singularity near the interface in our setting, we follow \cite{liusonburqia23} to analyze the sub-blocks near the interface separately. Letting $l_m=\frac{L}{2}$ denote the middle level of the butterfly, we can show that among the $O\left(n\right)$  subblocks $K(\tau, \nu)$ at each level $l$, there are $O\left(2^{\left|l-l_m\right| / 2} n^{1/4}\right)$ subblocks representing interactions between adjacent or close-by geometric subdomains. The non-constant ranks of these subblocks scale as $r_{\tau,\nu}=O(2^{-|l-l_m|/2}n^{1/4})$, which are dominated by the interface degrees of freedom (DOFs) between the two subdomains. Moreover, each non-constant-rank subblock requires $O(r^2_{\tau,\nu})$ memory usage and matrix entry computation and $O(r^3_{\tau,\nu})$ ID cost. Thus, they require 
\begin{equation}
\sum_l r_{\tau, \nu}^2 O\left(2^{\left|l-l_m\right| / 2} n^{1 / 4}\right)=O\left(n^{3 / 4}\right) \leq O\left(n\right) \quad \text { memory usage }
\end{equation}
and 
\begin{equation}
\sum_l r_{\tau, \nu}^3 O\left(2^{\left|l-l_m\right| /2} n^{1 / 4}\right)=O\left(n\right) \quad \text { complexity }.
\end{equation}
Hence, we can still state that the memory usage and complexity of the IDBF for $K\in\mathbb{C}^{O(n)\times O(n)} $ are $O(n\log(n))$. Due to the $O(n\log(n))$ nonzeros, we immediately conclude that the complexity of the application of the IDBF matrix is $O(n\log(n))$.

For HODBF of $K(\Omega_1,\Omega_1)\in \mathbb{C}^{O(N)\times O(N)}$, letting $L=O(\log(N))$ denote the level of the binary tree $\mathcal{T}_1$, we follow the DOFs analysis in \cite{micboa96} to sum up the complexity and memory usage of the $(L-l)$-th level butterfly for all levels $l$ and the $2^L$ dense blocks at level $L$, 
\begin{equation}
    2^L O(1)+\sum_{l=1}^L 2^l O(N/2^l\log(N/2^l))=O(N\log^2(N)).
\end{equation}
That is, the memory usage and complexity of the construction and application of the HODBF for $K(\Omega_1,\Omega_1)\in \mathbb{C}^{O(N)\times O(N)}$ are $O(N\log^2(N))$.


\subsection{Butterfly-compressed Eulerian Hadamard \color{b}integrators \color{black}}
\label{sec5.4}
For Case 1, starting from Algorithm 1 and using IDBF and HODBF, we precompute the low-rank representations of frequency-dependent kernels $\{U^{\omega_k}_i\}$ and corresponding time-dependent kernels $\{\check{U}_i^{t_\ell}\}$ and apply these compressed kernels in each time step, resulting in the butterfly-compressed TFT Hadamard integrator (Algorithm 2) with quasi-linear memory usage and computational complexity.\\

\noindent\textbf{Algorithm 2: Butterfly-compressed TFT Eulerian Hadamard integrator for Case 1}
\begin{enumerate}
    \item  Uniformly discretize the computational domain $\Omega$ into a wave-resolution-satisfying uniform grid; partition $\Omega$ into $N_I$ non-overlapping source regions $\{\Omega_S^{\ell},\; \ell=1, 2, \cdots, N_I\}$ and designate the corresponding receiver regions to be $\{\Omega_R^{\ell},\; \ell=1, 2, \cdots, N_I\}$,  set the initial time $T_0$,  the large time step $\Delta T<\tau_{\min}$ for updating the initial condition, the small time step $\Delta t=\frac{\Delta T}{N_T}$ for wave resolution, and an ending time $T_{end}$ to ensure that the waves do not reach the computational boundary; determine the bandwidth $B$ according to the initial conditions, set the frequency step $\Delta \omega = \frac{B}{N_{\omega}}$, where $N_{\omega}=O(B)$, and generate $\omega_k=(k+\frac{1}{2})\Delta \omega$ for $k = 0, 1,2,\cdots N_{\omega}-1$; precompute and compress Hadamard ingredients \rf{2.77} in all sub-region pairs.
    \item For each source region $\Omega_S^{\ell}$, divide the corresponding receiver region $\Omega_R^{\ell}$ into several subregions that are similar in size to the source region, including a subregion that coincides with the source region and several adjacent subregions;  compute the IDBF and HODBF representations of \color{g}
    \begin{equation*}
        U_s^{\omega_k}, s=1,2,3,4, 1\leq k\leq N_\omega-1 \quad\text{and}\quad \check{U}_s^{n\Delta t}, s=1,2,3,4, 1\leq n \leq N_T.
    \end{equation*}
    \color{black}
    \item Initialize $u(T_0,\bx)=u^1(\bx)$, 
 $u_t(T_0,\bx)=u^2(\bx)$; set time $T=T_0$ and loop variable $p=0$.
    \item For the current time step with $T=T_0+p \Delta T$:
\begin{enumerate}
        \item compute finite-difference derivatives of $u(T,\bx)$ and $u_t( T,\bx)$ and construct vectors $f_1$, $f_2$, $\cdots$, $f_8$.
        \item determine a subregion $\Omega_p$ of $\Omega$ that contains the region of influence of current data by extending each direction of the non-zero region of the solution $u(T,\bx)$ outward by $\max c \Delta T$; designate the source regions intersecting with $\Omega_p$ as the regions to be updated.
     \item  for all the source regions $\Omega^{\ell}_S$ to be updated:     
     \begin{enumerate}
    \item  apply the IDBF matrices and HODBF matrices to obtain the corresponding frequency- and time-sampling functions \color{g} 
         \begin{equation*}
       \hat{u}_i(\omega_k),\; i=1,2,3,4, \; 1\leq k \leq N_{\omega}-1, 
       \quad\text{and}\quad
 \check{u}_i(n\Delta t), i=1,2,3,4, \; 1\leq n \leq N_{T}.
   \end{equation*}
   \color{black}
    \item update wavefields $u(T+n\Delta t,\Omega_S^{\ell})$ and $u_t(T+n\Delta t, \Omega_S^{\ell})$ via \rf{4.46} and \rf{4.47}.  
    \end{enumerate}
    \item set to zero the wave fields within the source regions that do not require updating.
    \end{enumerate}
    \item Update $T=T_0+(p+1)\Delta T.$ If $T<T_{end}$, then $p\leftarrow p+1$ and go to Step 4; else, stop.
\end{enumerate}

Equipped with butterfly representations, we construct the frequency-domain and time-domain kernels in the precomputation step and recursively apply them during time marching. The computational complexity and memory usage of constructing and applying each kernel are reduced from $O(B^4)$ to $O(B^2\log^2(B))$, resulting in the Eulerian Hadamard integrator with overall complexity and memory usage of $O(B^3\log^2(B))$. Although we divide $\Omega$ into $N_I$ parts, $N_I=O(1)$ which only depends on the medium; consequently, such partition does not affect the overall complexity estimate. 

 Further, assume that we solve the wave equation with $m_1$ different initial conditions simultaneously. In this case, we only construct the $O(B)$ kernels once, and the computational cost is 
 \begin{equation*}
      O(B)O(B^2 \log^2(B)) \doteq  O(C_1 B^3 \log^2(B));  
 \end{equation*}
applying all these kernels to $m_1$ initial conditions, the computational cost is $O(m_1 C_2B^3\log^2(B))$; we use the fast Fourier transform to obtain the time-domain solutions on $O(B^3)$ mesh points, with a cost of $O(m_1 B^3\log(B))$. Since we only need to update the initial conditions for $\frac{T_{end}}{\Delta T} = O(1)$ times, the overall computational complexity is 
    \begin{equation*}
        O\left((C_1+\frac{T_{end}}{\Delta T}m_1C_2)B^3\log^2(B)\right).
    \end{equation*}
    For the butterfly algorithm, $C_1$ depends on the accuracy of the representation and $r^3_{\tau,\nu}$ from the IDs. In any case, $C_1$ is far greater than $C_2$; thus, when solving different Cauchy problems in the same medium, the efficiency of the algorithm will be more significant. On the other hand, even if IDBF and HODBF need to be recalculated at each step due to limited computing resources, the overall computational complexity remains quasi-linear. 

    The quasi-linear memory usage makes parallel computing feasible. Since there are $O(B)$ kernels to compress, we do not embed parallel computing into the decomposition of each kernel matrix. Instead, we use the parallel toolbox in Matlab to simultaneously calculate the butterfly representations of $O(B)$ kernel matrices.

    \color{g}
For Case 2, we first use \eqref{3.2.4} to approximate the time-domain solution in the interval $0< t \leq \Delta T$ and then apply \textbf{Algorithm 2} to solve the wave equation \eqref{4.62.0}-\eqref{4.64.0} in the interval $\Delta T < t \leq T_{end}$. Subsequently, we perform the  Fourier transform \rf{tftf} to obtain the frequency-domain solution. This leads to the butterfly-compressed TFTF Eulerian Hadamard integrator. Here, we only perform an additional Fourier transform, which means that the computational complexity and memory usage of the butterfly-compressed TFTF Eulerian Hadamard integrator remain $O(C_1B^3\log^2(B))$. 

Compared to the TFT method, which yields $O(B^3)$ unknowns in the time-space domain, the TFTF method obtains $O(B^3)$ unknowns in the frequency-space domain. Namely, both Eulerian Hadamard integrators achieve a quasi-linear computational complexity.
\color{black}

\section{Numerical examples}
\label{sec6}
This section provides several numerical examples to demonstrate the accuracy and efficiency of the proposed butterfly-compressed Hadamard integrators by solving wave equations in time and frequency domains. All computations were performed on a computer equipped with 512GB of RAM and two 28-core Intel Xeon 6258R processors. 

Because exact solutions for the wave equations are not available in general, we solve the wave equation using  \color{b} the FDTD method \cite{tafhagpik05} for time-domain ``reference'' solutions and a finite-difference based direct solver \cite{joshisuh96} for frequency-domain ``reference'' solutions, where the perfectly-matched layer (PML) absorbing boundary conditions are imposed \cite{grosim10,ber94,engyin11}, \color{black} and we always keep the number of points per wavelength (PPW) greater than 40 to reduce the dispersion error in the two direct methods. For the proposed Hadamard integrator, we set PPW$\sim 5$ with respect to the artificial bandwidth $B$.

\subsection{Sinusoidal model}
We introduce a sinusoidal model with the following setup.
\begin{itemize}
    \item $\rho=\fr{1}{(1+0.2\sin(3\pi(x+0.05))\sin(0.5\pi y))^2}$, $\nu=1$ such that $c=1+0.2\sin(3\pi(x+0.05))\sin(0.5\pi y).$
    \item $B=64\pi$, $h=\frac{1}{200}$, $\Delta\omega=\frac{\pi}{2}$, and $\Delta t=\frac{1}{192}$.
    \item $\Omega=[-1,2]\times[-1,2]$; the sizes of $\Omega_S^{\ell}$ and $\Omega_R^{\ell}$ are $0.2\times 0.2$ and $0.6\times 0.6$,  respectively; $\Delta T=\frac{1}{8}$.
     \item Orders of Chebyshev interpolations of Hadamard ingredients in the four variables $x_0$, $y_0$, $x$, and $y$ are $11$, $11$, $13$, and $13$, respectively.
    \item Tolerance used in the interpolative decomposition: ${\rm tol} = 10^{-9}$.
\end{itemize}
\color{b}
We show the velocity model, some rays, and wavefronts in Fig.~\ref{velocity1}, where the rays and wavefronts are obtained by solving the eikonal equation using the method of characteristics. Caustics occur when rays intersect and wavefronts fold. Note that in the current domain, caustics develop only along the $y$ direction. The frequency-domain fast Huygens sweeping method developed in \cite{luoqiabur14a,luqiabur16} can also be used to compute the point-source wave field in this medium by propagating waves layer by layer along the $y$ direction.

     \begin{figure}[htbp]
     \centering
     \subfigure[]{
     \includegraphics[scale=0.4]{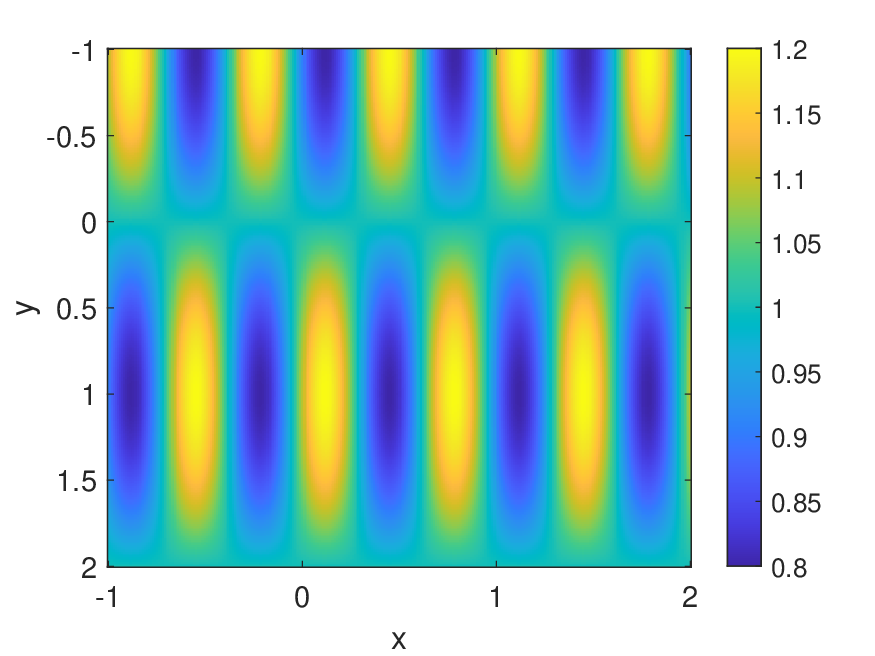}}
     \subfigure[]{
     \includegraphics[scale=0.4]{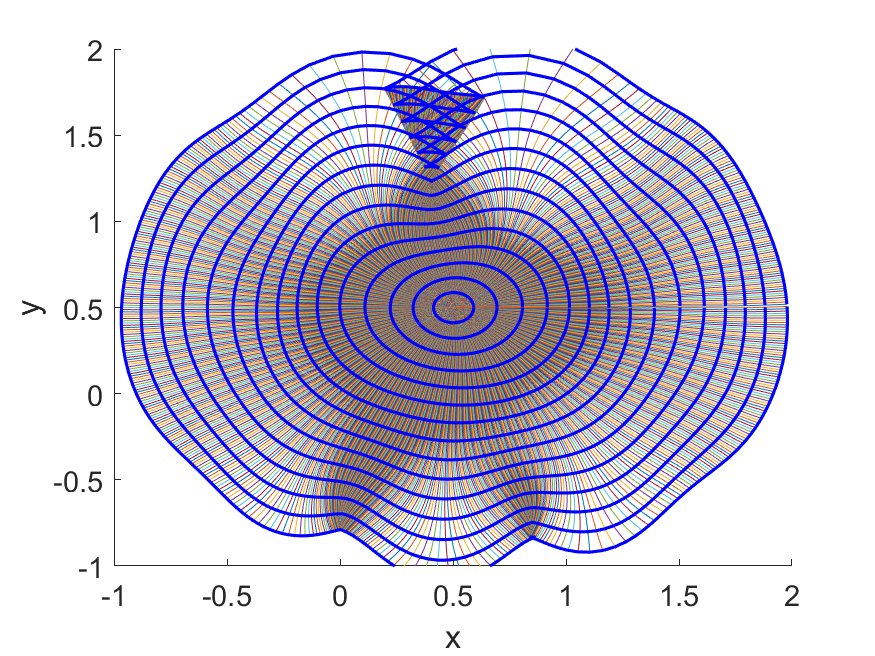}}
     \caption{\color{b} Sinusoidal model. (a) The velocity; (b) Rays and wavefronts with source $\bx_0=[0.5,0.5]$. The thick blue lines represent equal-time wavefronts (traveltime contours) with the contour interval equal to $0.1$, and thin colored lines represent rays with different take-off angles}
     \label{velocity1}
     \end{figure}
\color{black}

\subsubsection{TFT method for Case 1}
For Case 1, we take the initial conditions as
\begin{equation}
    u(0,x,y)=\sin\left(\beta \pi \frac{x+y-1}{\sqrt{2}}\right)\exp\left(-100((x-0.5)^2+(y-0.5)^2)\right),\; u_t(0,x,y)=0,
\end{equation}
where $\beta=8$, $16$ and $32$, respectively, and we utilize the Hadamard integrator to simultaneously solve these three Cauchy problems. Fig. \ref{figure2.0} shows the comparisons between the TFT solutions and the \color{b}reference \color{black} solutions in $[0,1]\times [0,1]$ at $t=0.5$. In  Fig. \ref{figure3.0}, we further compare the solutions along some lines. \color{b} Table \ref{table1} shows the relative $L^2$ and $L^{\infty}$ errors between the TFT solutions and the reference solutions. 
\color{b}

Here, we observe that compared with the case of $\beta=16$, the TFT method shows larger errors for the cases of $\beta=8$ and $\beta=32$. The reason for the case of $\beta=8$ is that the asymptotic error of the Hadamard integrator dominates in the low-frequency band, which is typical of the error behavior for microlocal-analysis based asymptotic methods \cite{liusonburqia23}; on the other hand, as $\beta$ gradually increases to reach or even exceed the Nyquist frequency $B/2$, the time-domain solution reconstructed through 
frequency-domain samplings will no longer be accurate, as indicated by the Shannon-Nyquist Sampling theorem. Similar error performance will be observed in subsequent numerical examples. We will further discuss the error behaviors and convergence results in detail in Section \ref{sec6.4}.

\color{black}
  \begin{figure}[htbp]
     \centering
                     \subfigure[]{
     \includegraphics[scale=0.35]{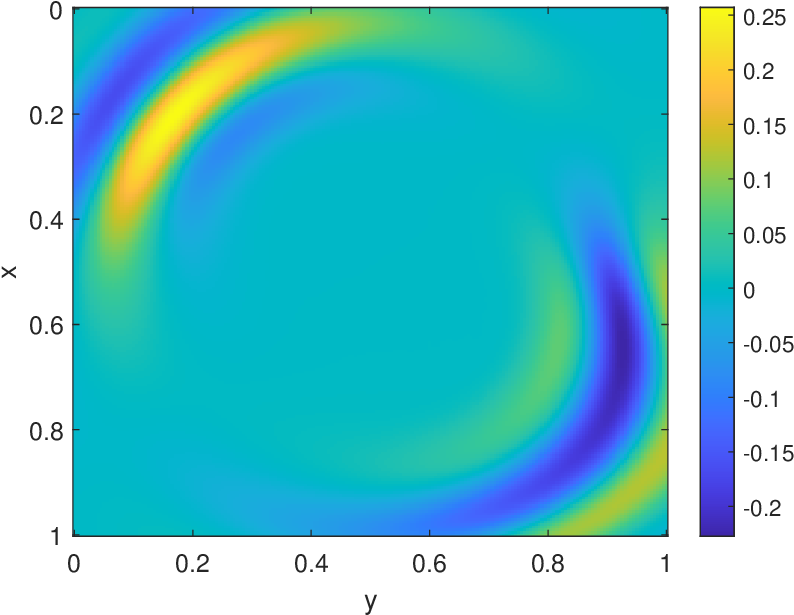}}
                \subfigure[]{
     \includegraphics[scale=0.35]{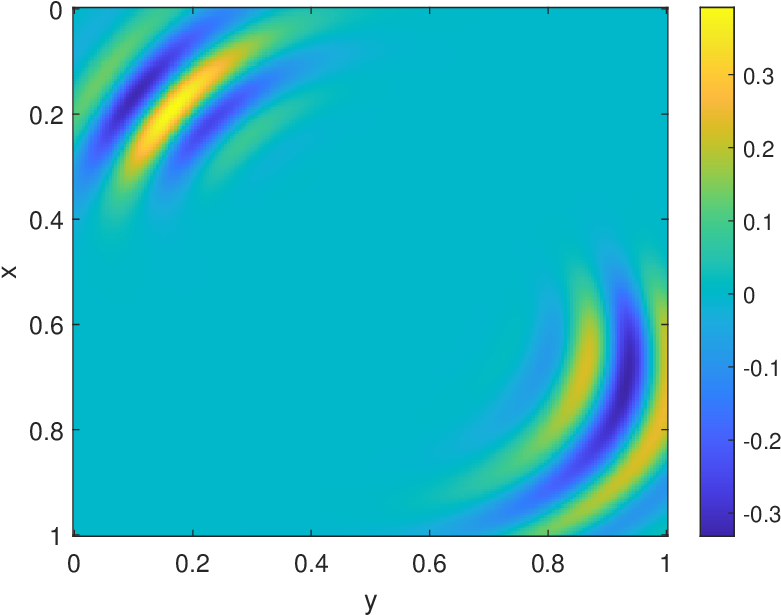}}
                \subfigure[]{
     \includegraphics[scale=0.35]{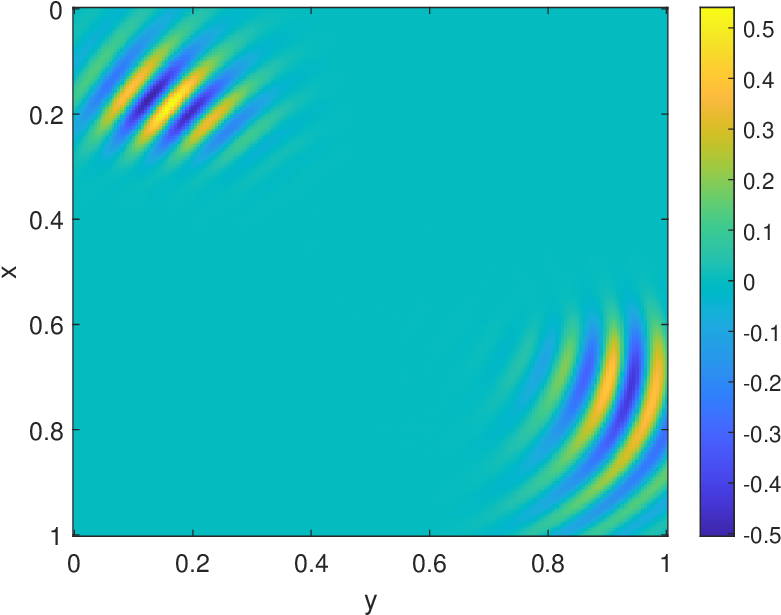}}\\
               \subfigure[]{
     \includegraphics[scale=0.35]{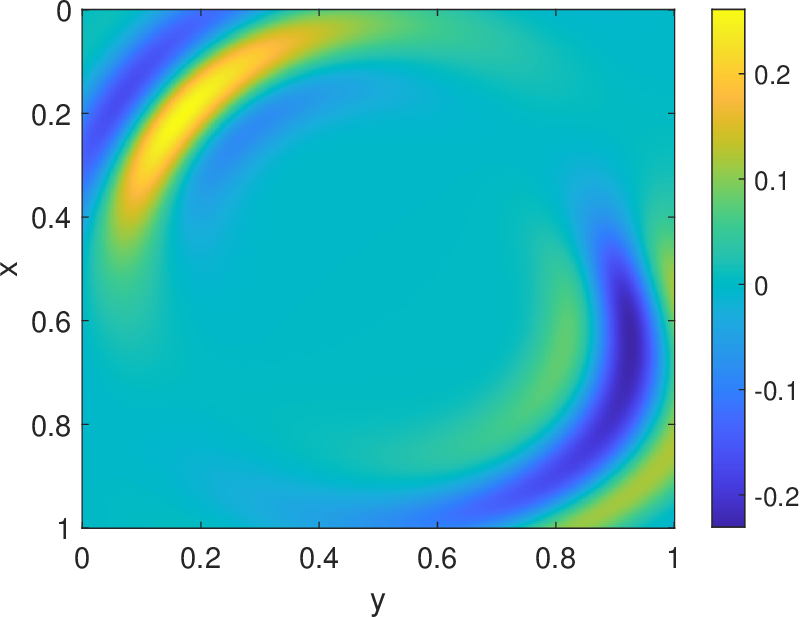}}
                    \subfigure[]{
     \includegraphics[scale=0.35]{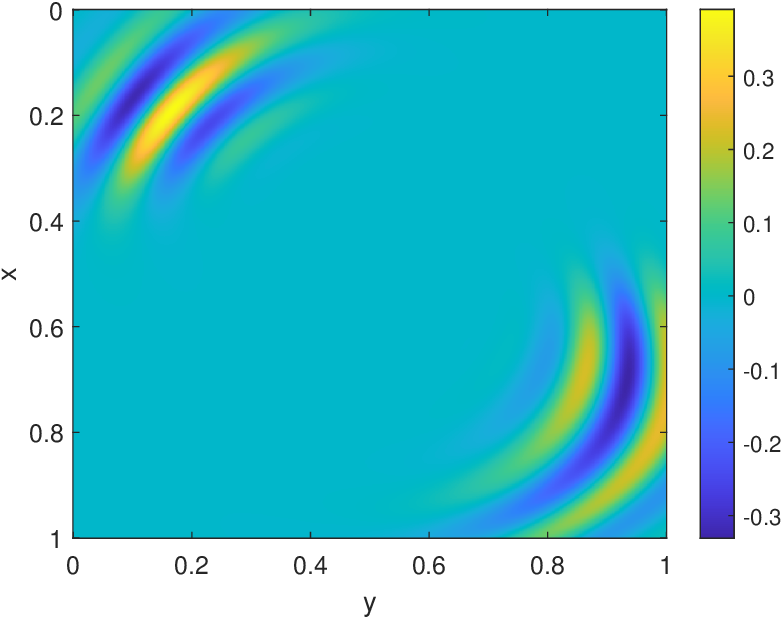}}
                    \subfigure[]{
     \includegraphics[scale=0.35]{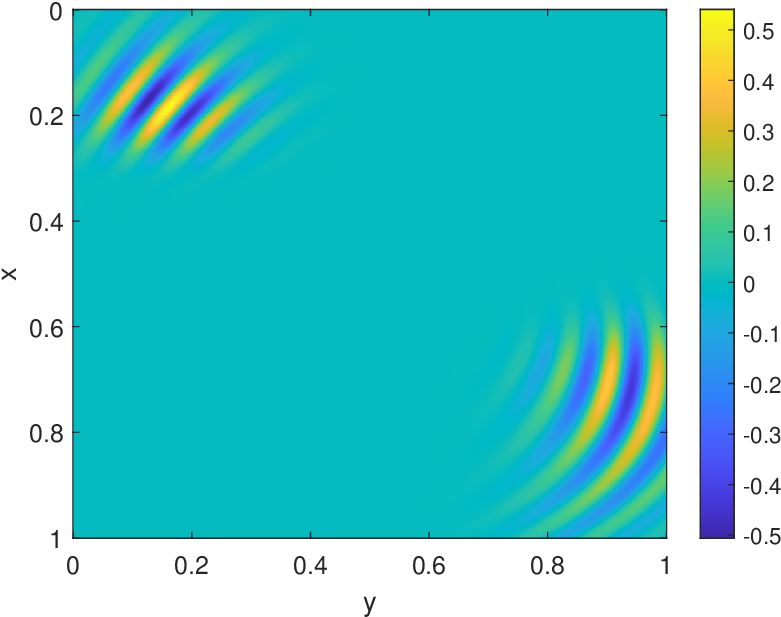}}\\
     \caption{Sinusoidal model. $u(t,\bx)$, $t=0.5$: (a) TFT solution with $\beta=8$; (b) TFT solution with $\beta=16$; (c) TFT solution with $\beta=32$; (d) \color{b}reference \color{black} solution with $\beta=8$;  (e) \color{b}reference \color{black} solution with $\beta=16$;  (f) \color{b}reference \color{black} solution with $\beta=32$ }
     \label{figure2.0}
     \end{figure}   
     
       \begin{figure}[htbp]
     \centering
               \subfigure[]{
     \includegraphics[scale=0.35]{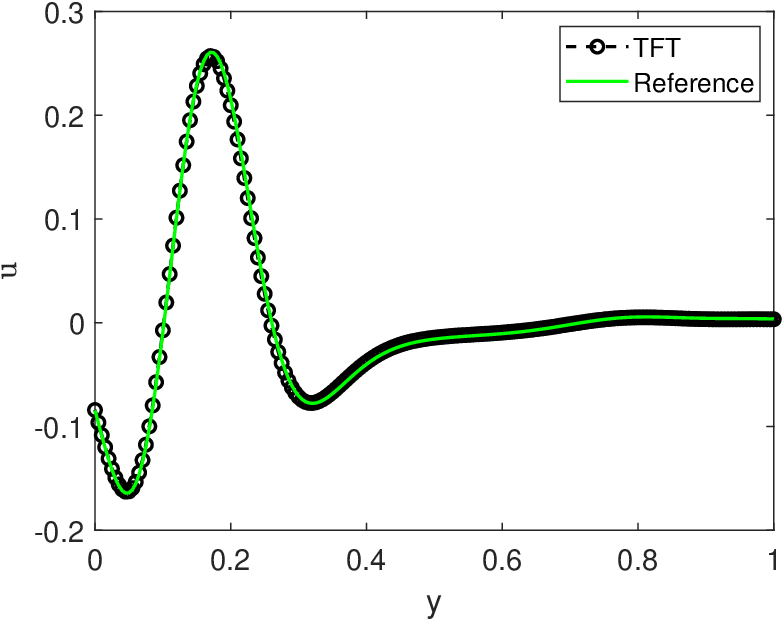}}
                    \subfigure[]{
     \includegraphics[scale=0.35]{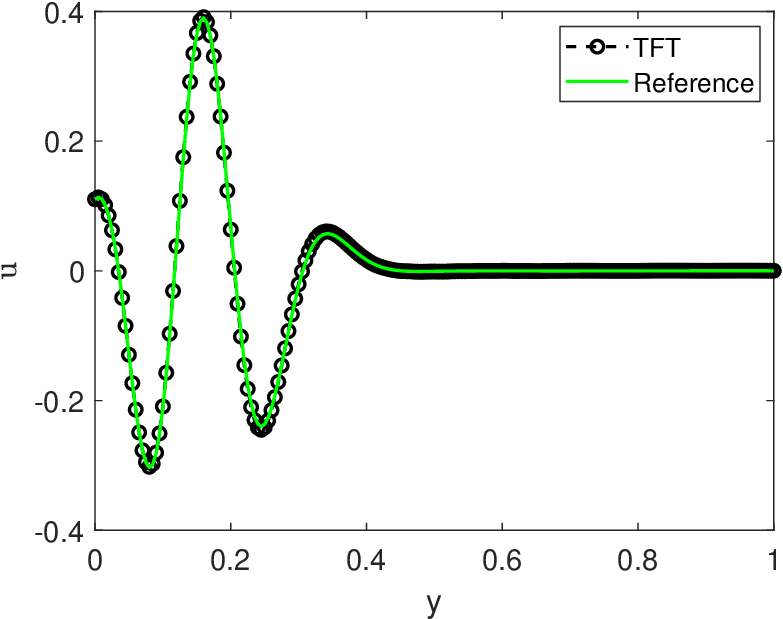}}
                    \subfigure[]{
     \includegraphics[scale=0.35]{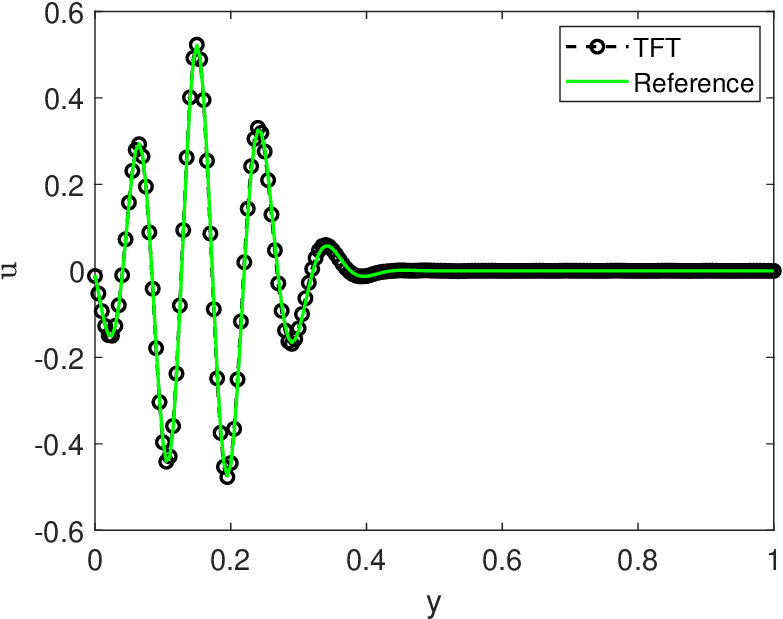}}\\
                \subfigure[]{
     \includegraphics[scale=0.35]{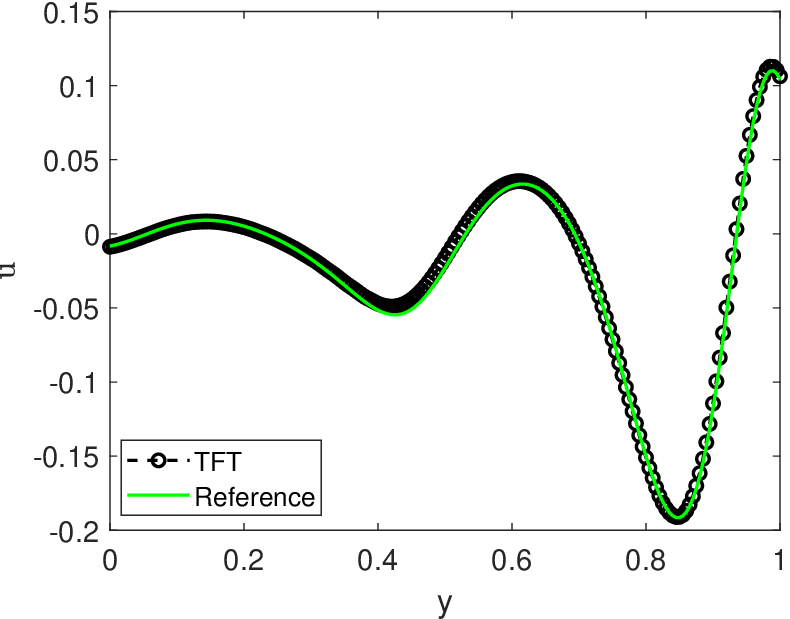}}
                \subfigure[]{
     \includegraphics[scale=0.35]{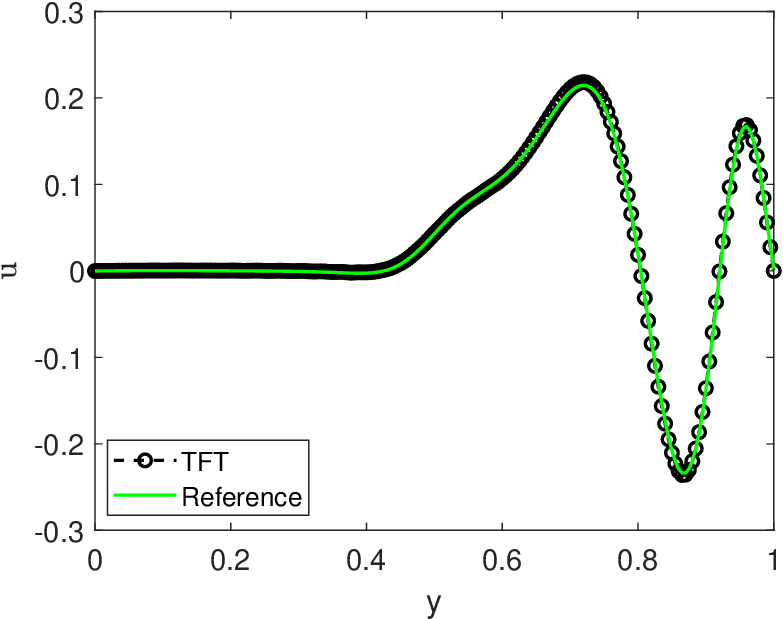}}
                \subfigure[]{
     \includegraphics[scale=0.35]{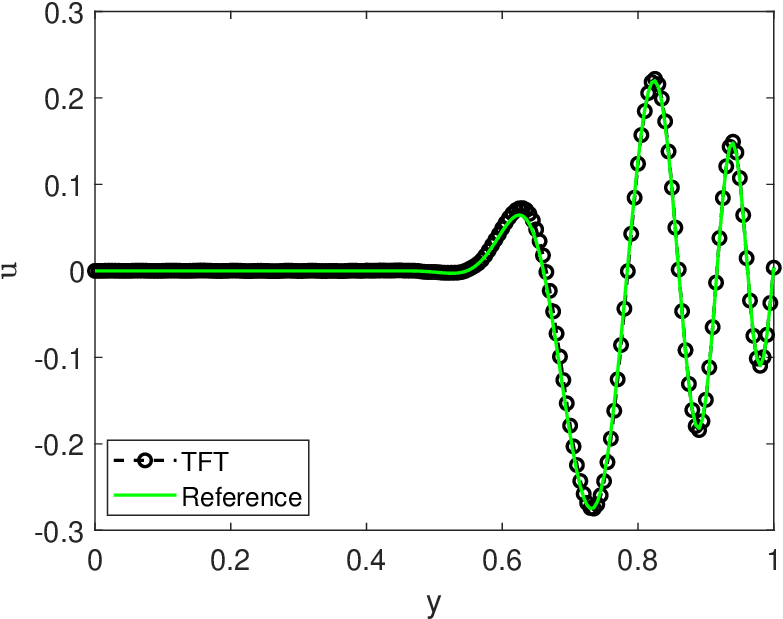}}\\
     \caption{Sinusoidal model. 
     Slices of time-domain wave fields $u(t,\bx)$ when $t=0.5$: (a) a slice at $x=0.2$ with $\beta=8$; (b) a slice at $x=0.2$ with $\beta=16$; (c) a slice at $x=0.2$ with $\beta=24$; (d) a slice at $y=0.85$ with $\beta=8$; (e) a slice at $y=0.85$ with $\beta=16$; (f) a slice at $y=0.85$ with $\beta=24$
     }
     \label{figure3.0}
     \end{figure}
     
\begin{table}[ht]
\centering
\caption{The relative $L^2$ and $L^{\infty}$ errors of TFT solutions for Sinusoidal model at $t=0.5$}\label{table1}
\begin{tabular}{ccccc}
\toprule
$\beta$ & 8 & 16  &  32\\
\midrule
Relative $L^2$ error & $3.22e-2$ & $1.97e-2$  & $3.77e-2$ \\
Relative $L^{\infty}$ error & $3.08e-2$ & $1.88e-2$  & $3.64e-2$\\
\bottomrule
\end{tabular}
\end{table}

\subsubsection{TFTF method for Case 2}
For Case 2, taking $\bz=[0.5,0.5]$ and $T_{end}=1.5$, we obtain the time-domain point-source wave fields $u(t, \bx)$. As shown in Fig. \ref{figure2}, we can clearly observe the occurrence of caustics.

 \begin{figure}[htbp]
     \centering
          \subfigure[]{
     \includegraphics[scale=0.35]{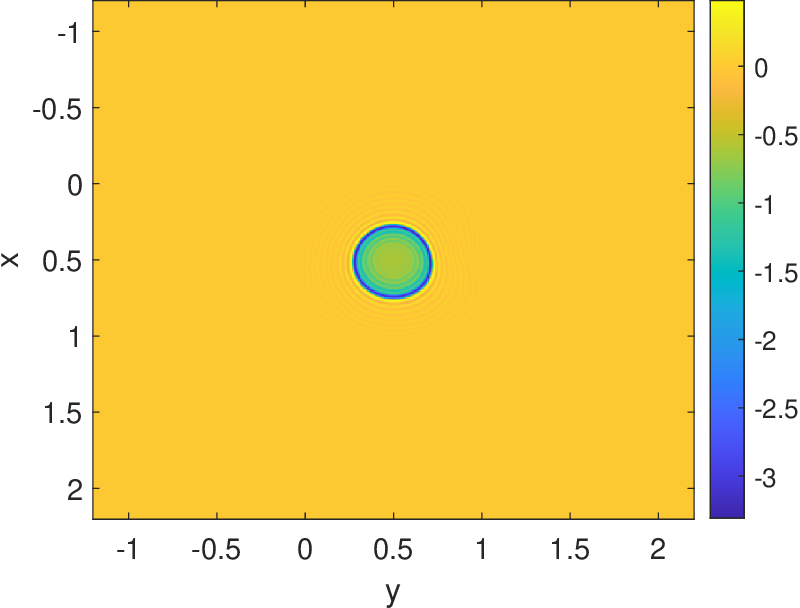}}
                \subfigure[]{
     \includegraphics[scale=0.35]{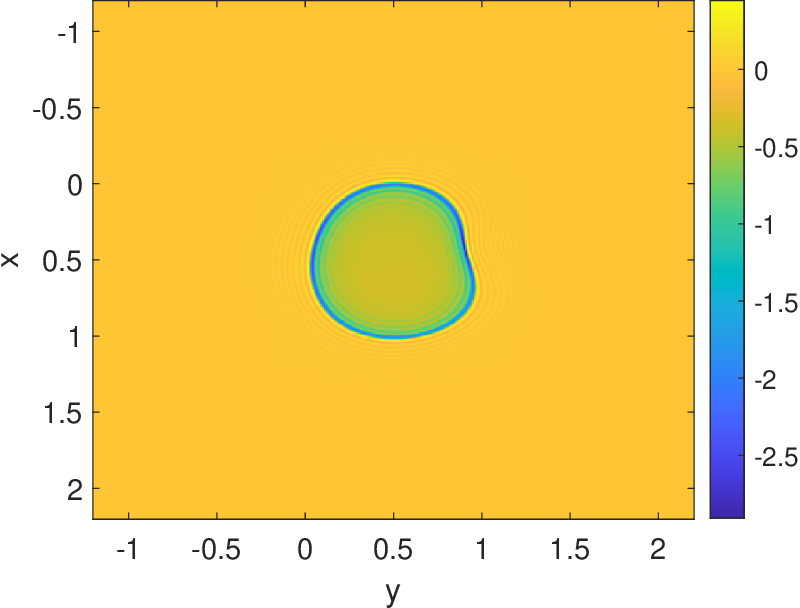}}
     \subfigure[]{
     \includegraphics[scale=0.35]{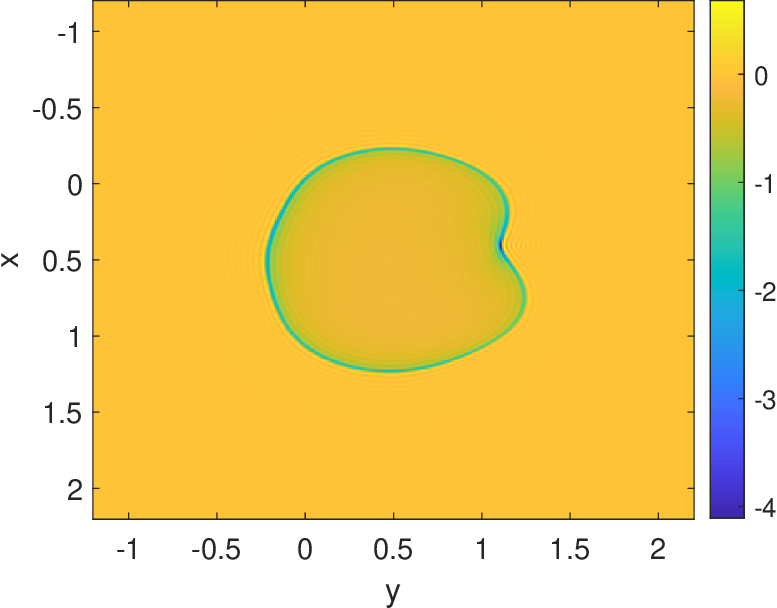}}\\
          \subfigure[]{
     \includegraphics[scale=0.35]{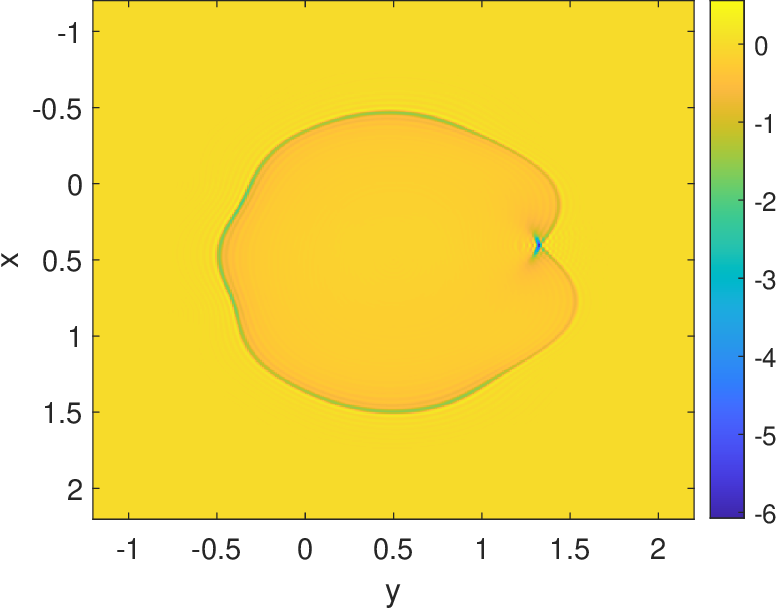}}
                \subfigure[]{
     \includegraphics[scale=0.35]{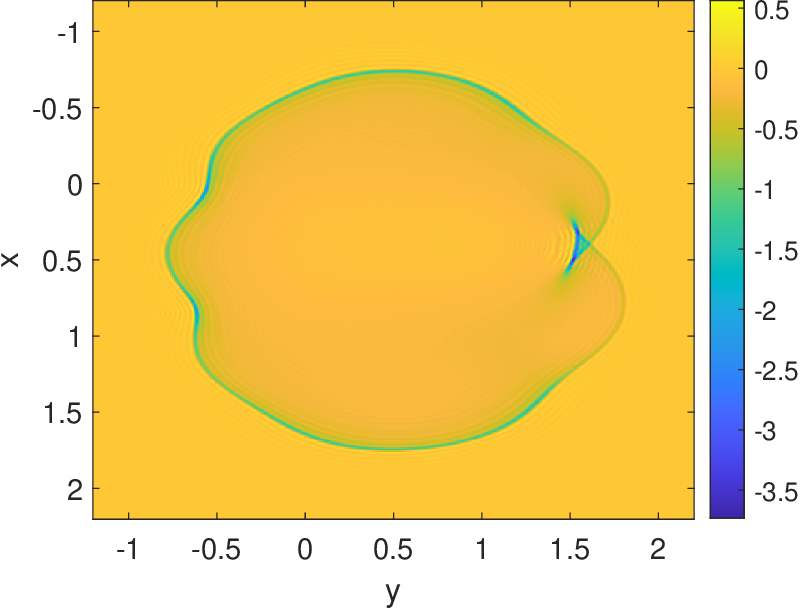}}
     \subfigure[]{
     \includegraphics[scale=0.35]{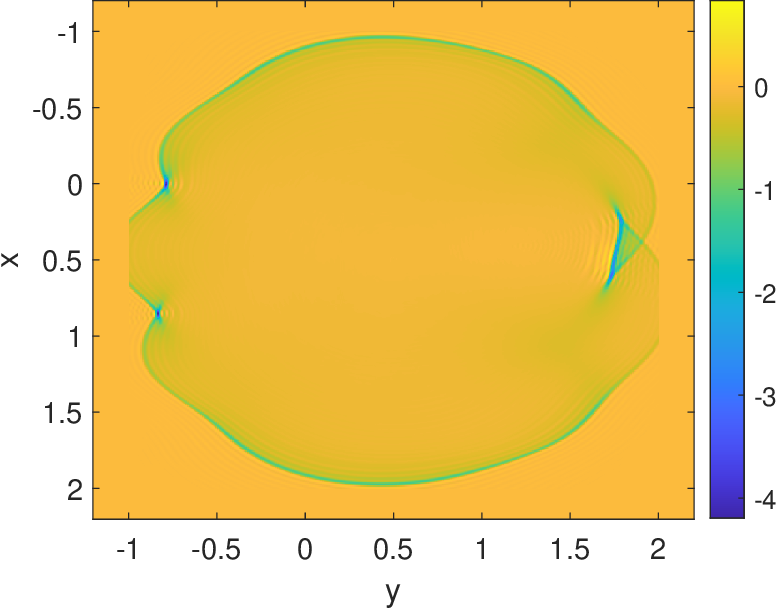}}\\
     \caption{Sinusoidal model. Time-domain point-source wave fields $u(t,\bx)$. (a) $T=0.25$; (b) $T=0.5$; (c) $T=0.75$; (d) $T=1$; (e) $T=1.25$; (f) $T=1.5$ }\label{figure2}
     \vspace{-3.75mm}
     \end{figure}

Further applying Fourier transform in time, we obtain the frequency-domain point-source wave fields $\hat{u}(\omega,\bx)$ with different angular frequencies $\omega$. In Fig. \ref{figure3}, we compare the TFTF solution with the \color{b}reference \color{black} solution for $\omega=8\pi$, $16\pi$ and $32\pi$, respectively, where $\bx\in [0,1.5]\times[0,1.5]$. \color{b} Table \ref{table2} shows the relative $L^2$ and $L^{\infty}$ errors between the TFTF solutions and the reference  solutions.
\color{black}

 \begin{figure}[htbp]
     \centering
          \subfigure[]{
     \includegraphics[scale=0.35]{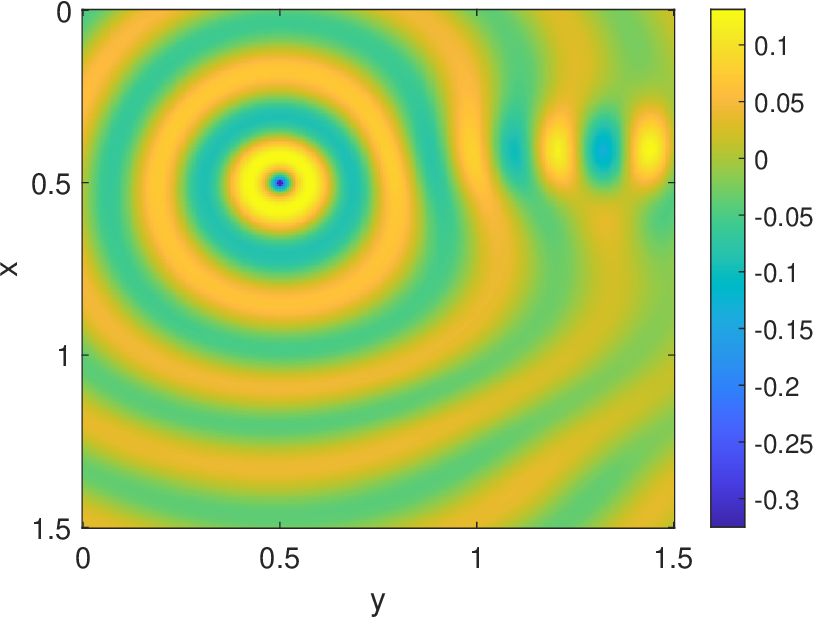}}
                     \subfigure[]{
     \includegraphics[scale=0.35]{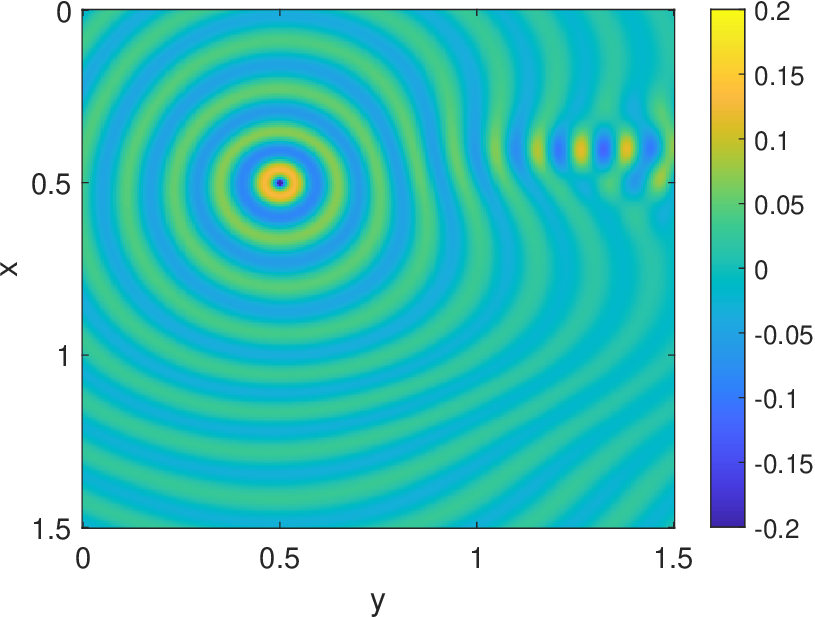}}
               \subfigure[]{
     \includegraphics[scale=0.35]{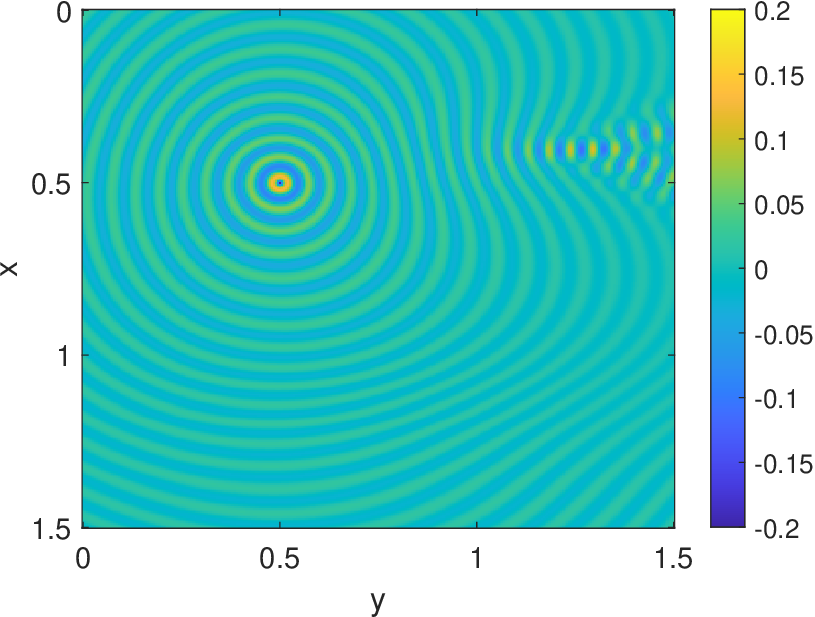}}\\
                    \subfigure[]{
     \includegraphics[scale=0.35]{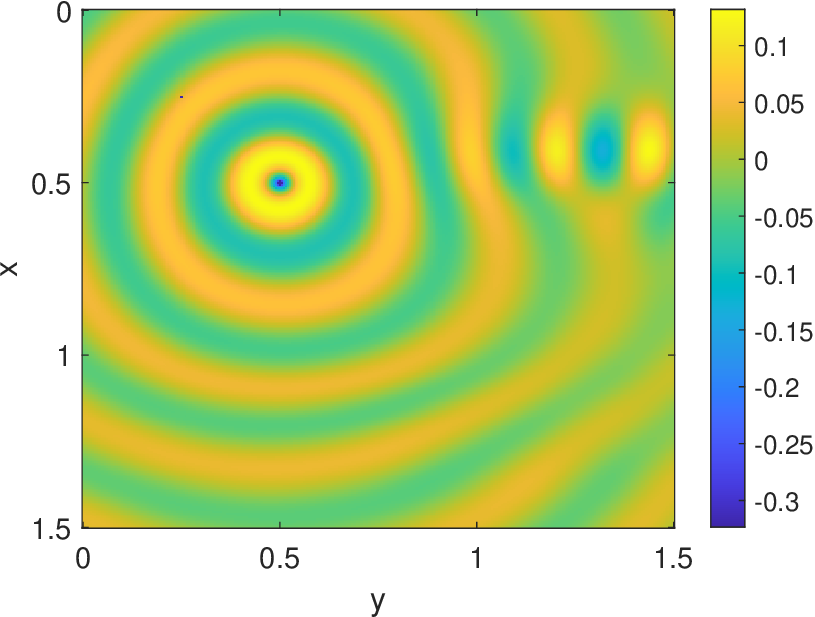}}
                     \subfigure[]{
     \includegraphics[scale=0.35]{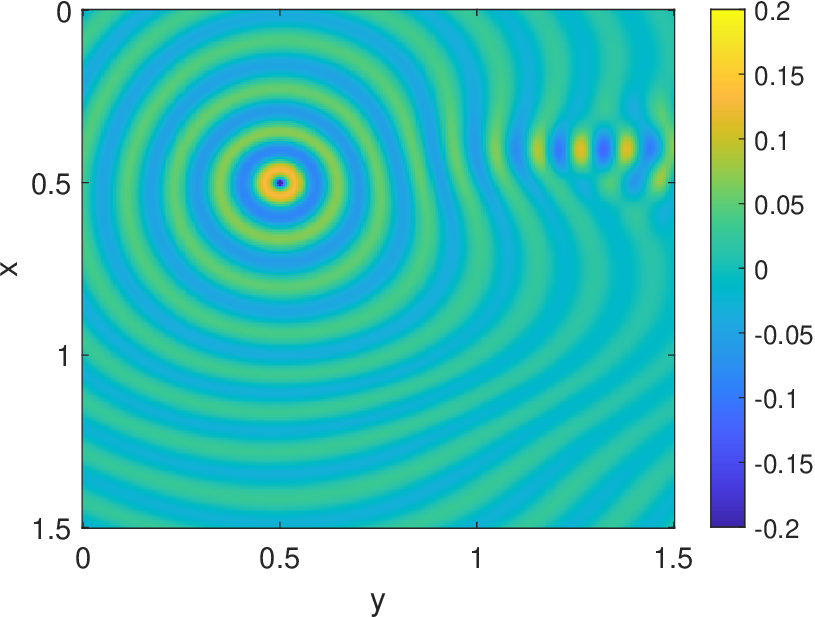}}
          \subfigure[]{
     \includegraphics[scale=0.35]{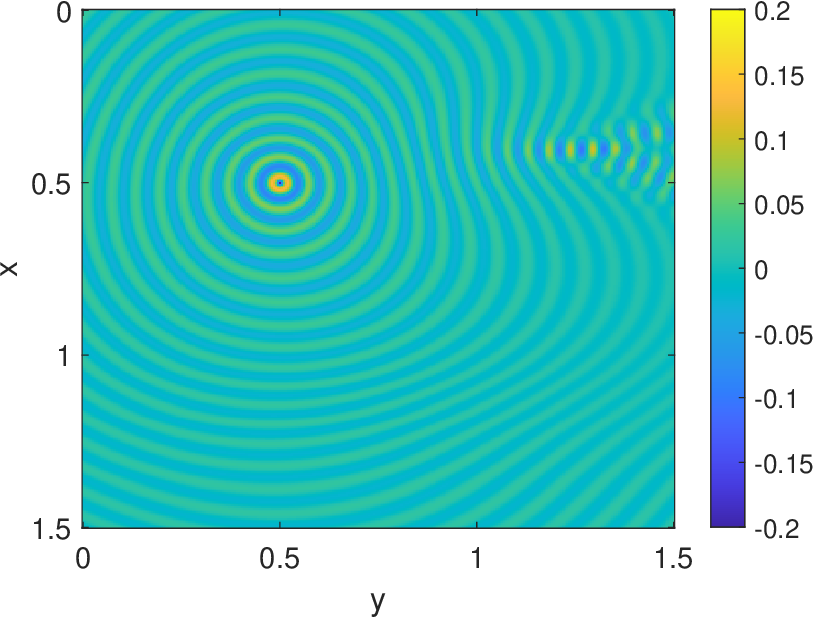}}\\
     \caption{Sinusoidal model. $\hat{u}(\omega,\bx)$. (a) TFTF solution with $\omega=8\pi$; (b) TFTF solution with $\omega=16\pi$; (c) TFTF solution with $\omega=32\pi$; (d) \color{b}reference \color{black} solution with $\omega=8\pi$; (e) \color{b}reference \color{black} solution with $\omega=16\pi$; (f) \color{b}reference \color{black} solution with $\omega=32\pi$ }\label{figure3}
     \vspace{-3.75mm}
     \end{figure}
     \begin{table}[ht]
\centering
\caption{The $L^2$ and $L^{\infty}$ errors of TFTF solutions for Sinusoidal model}\label{table2}
\begin{tabular}{ccccc}
\toprule
$\omega$ & $8\pi$ & $16\pi$  &  $32\pi$\\
\midrule
Relative $L^2$ error & $3.21e-2$ & $2.48e-2$  & $5.88e-2$ \\
Relative $L^{\infty}$ error & $3.81e-2$ & $2.79e-2$  & $3.84e-2$\\
\bottomrule
\end{tabular}
\end{table}
Fig. \ref{figure4} shows the comparisons of the two solutions along lines that traverse through the caustic region for various  $\omega$. Although they match well, it can be observed that the errors at $\omega=8\pi$ and $32\pi$ are slightly larger than those at  $\omega=16\pi$, and we will explain this phenomenon in Section \ref{sec6.4}.

      \begin{figure}[htbp]
     \centering
          \subfigure[]{
     \includegraphics[scale=0.35]{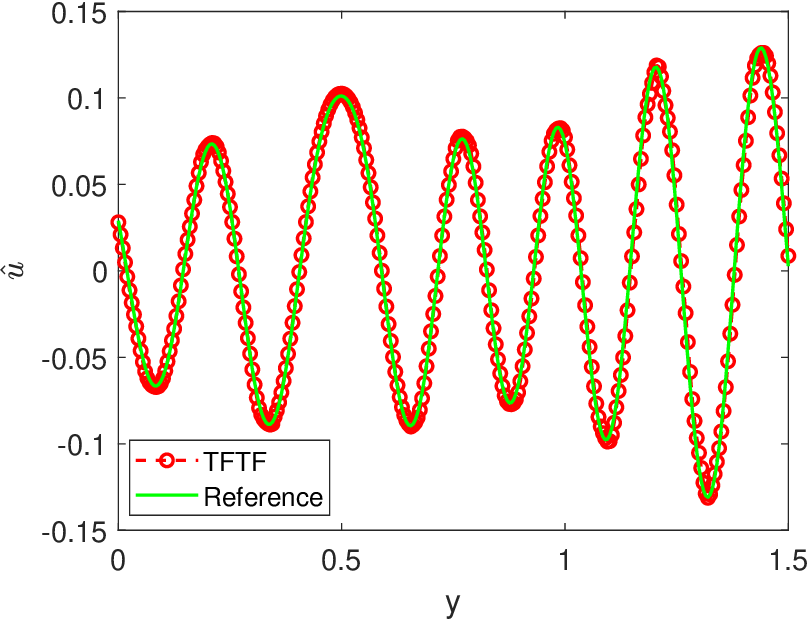}}
               \subfigure[]{
     \includegraphics[scale=0.35]{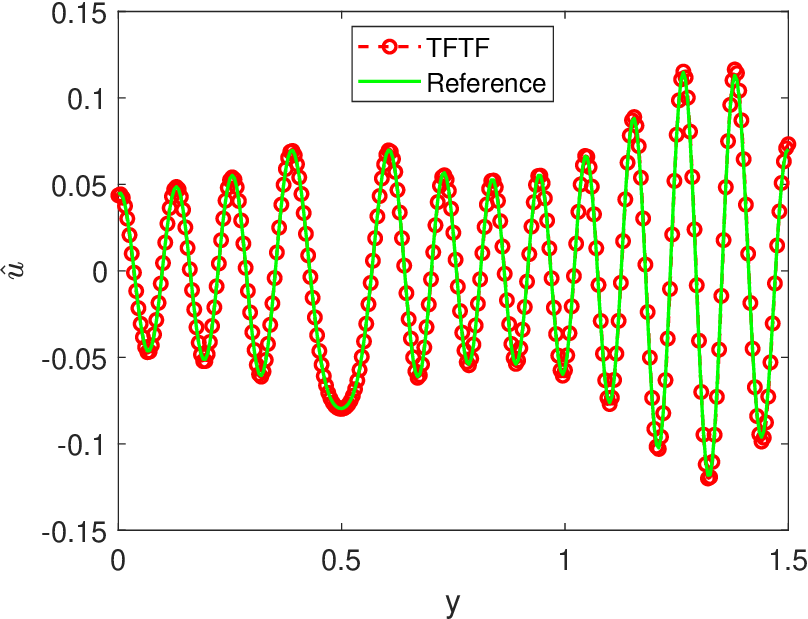}}
                    \subfigure[]{
     \includegraphics[scale=0.35]{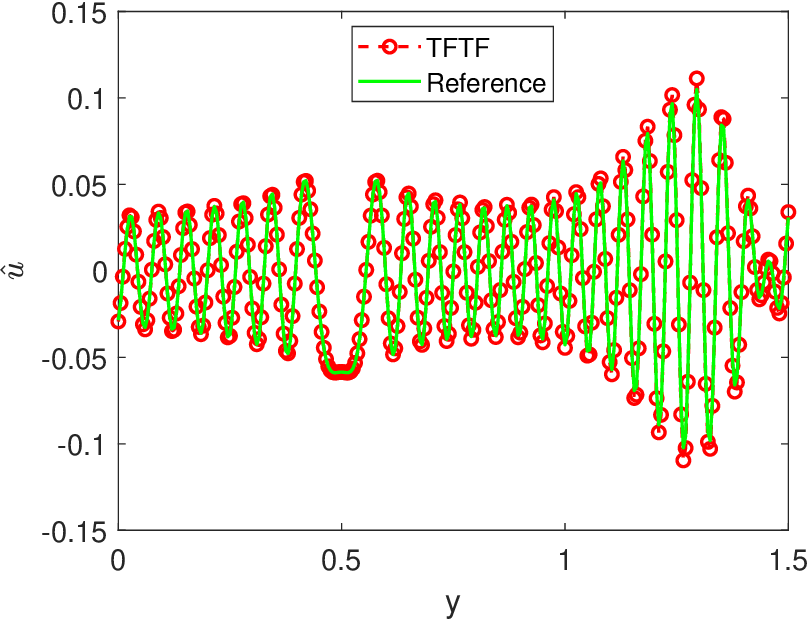}}\\
                \subfigure[]{
     \includegraphics[scale=0.35]{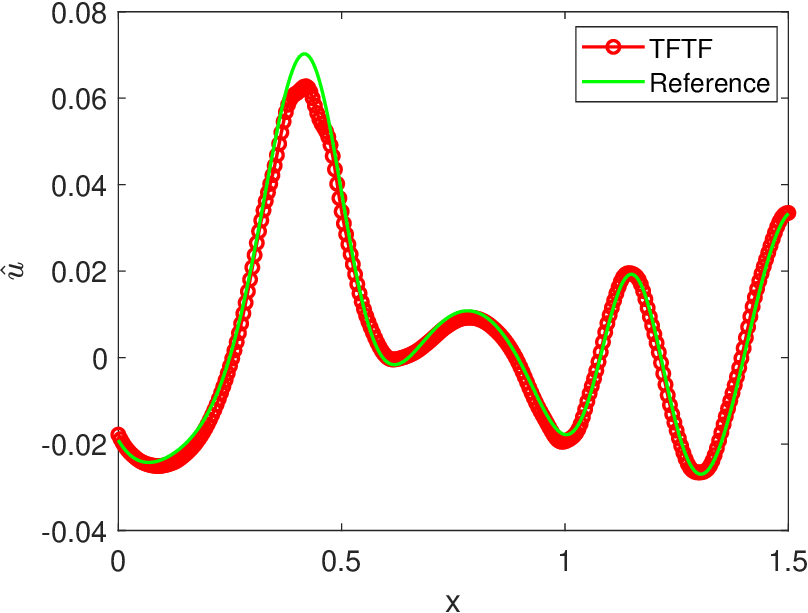}}
                \subfigure[]{
     \includegraphics[scale=0.35]{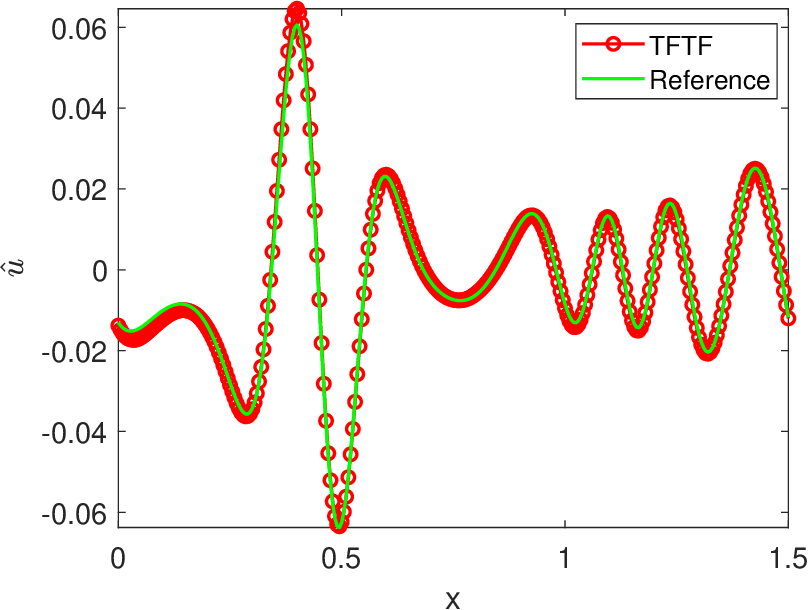}}
                \subfigure[]{
     \includegraphics[scale=0.35]{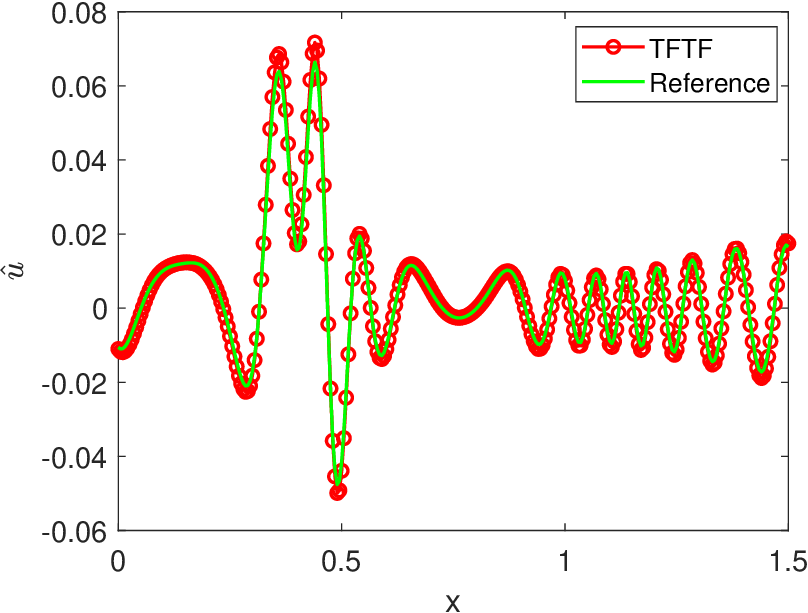}}\\
     \caption{Sinusoidal model. Slices of frequency-domain point source wave fields $\hat{u}(\omega,\bx)$. (a) a slice at $x=0.4$ with $\omega=8\pi$; (b) a slice at $x=0.4$ with $\omega=16\pi$; (c) a slice at $x=0.4$ with $\omega=32\pi$; (d) a slice at $y=1.4$ with $\omega=8\pi$; (e) a slice at $y=1.4$ with $\omega=16\pi$; (f) a slice at $y=1.4$ with $\omega=32\pi$}\label{figure4}
     \vspace{-3.75mm}
     \end{figure}
     
\subsection{Smoothed Heaviside model}
We introduce a $y$-dependent velocity model $c$, which is analogous to a scaled, smoothed, and shifted Heaviside function.
\begin{itemize}
\item $\rho=\left(\fr{1+e^{-20(y-1)}}{1.25+0.8 e^{-20(y-1)}}\right)^2$ and $\nu=1$ such that $c=0.8+\fr{0.45}{1+e^{-20(y-1)}}$.
\item $B=48\pi$, $h=\frac{1}{160}$, $\Delta\omega=\frac{\pi}{2}$, and $\Delta t=\frac{1}{144}$.
\item $\Omega=[-3,3]\times[-1,4]$; the sizes of $\Omega_S^{\ell}$ and $\Omega_R^{\ell}$ are $0.2\times 0.2$ and $0.6\times 0.6$,  respectively; $\Delta T=\frac{1}{8}$.
\item Orders of Chebyshev interpolations of Hadamard ingredients in the four variables $x_0$, $y_0$, $x$, and $y$ are $11$, $11$, $13$, and $13$, respectively.
\item Tolerance used in interpolative decomposition: ${\rm tol} = 10^{-9}$.
\end{itemize}
\color{b}
We show the velocity model, some rays, and wavefronts in Fig.~\ref{velocity2}. The rapid change in velocity near $y=1$ results in overturning rays and caustics. The FHS method, which propagates the wave field along a specific spatial direction, is no longer valid due to these overturning rays. However, the proposed Hadamard integrator, which propagates the wave field along the time direction, can naturally handle such caustics.
\color{black}
     \begin{figure}[htbp]
     \centering
     \subfigure[]{
     \includegraphics[scale=0.4]{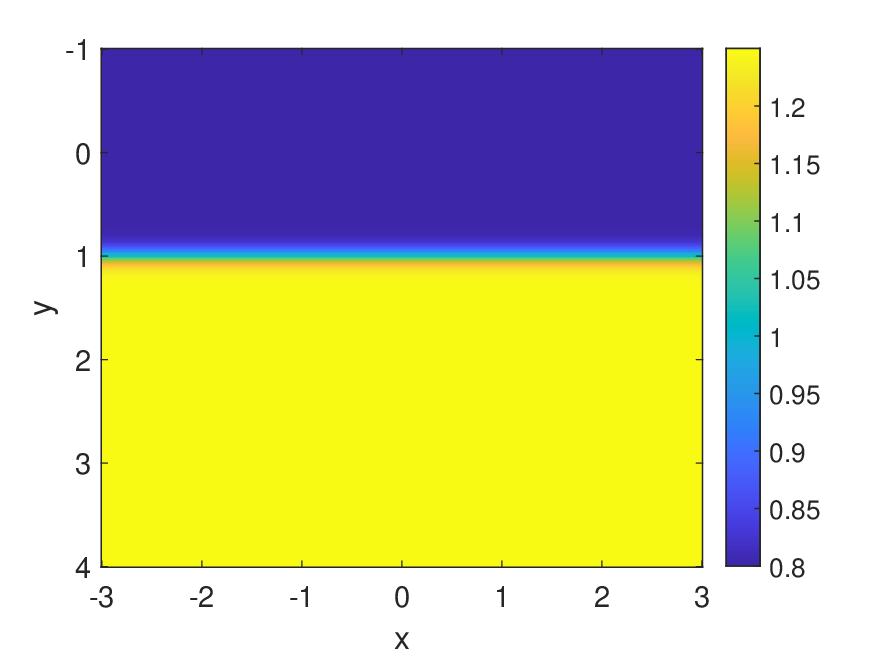}}
     \subfigure[]{
     \includegraphics[scale=0.4]{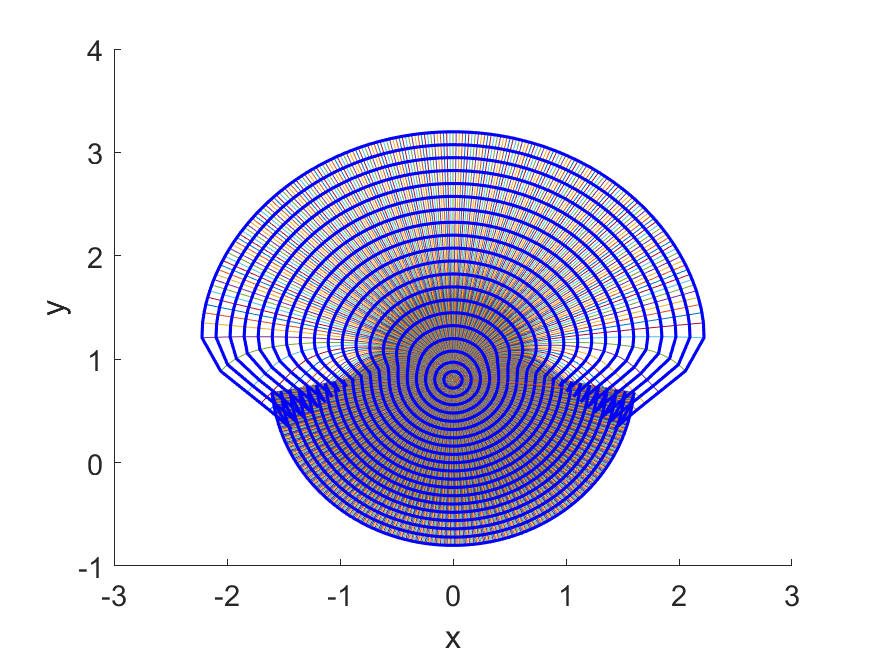}}
     \caption{\color{b} Heaviside model. (a) The velocity; (b) Rays and wavefronts with source $\bx_0=[0,0.8]$. The thick blue lines represent equal-time wavefronts (traveltime contours) with the contour interval equal to $0.1$, and thin colored lines represent rays with different take-off angles}
     \label{velocity2}
     \end{figure}
\color{black}
\subsubsection{TFT method for Case 1}
For Case 1, we take the initial conditions as
\begin{equation}
    u(0,x,y)=\sin\left(\beta \pi \frac{x+y-1}{\sqrt{2}}\right)\exp(-100(x^2+(y-0.8)^2)),\;\;\; u_t(0,x,y)=0,
\end{equation}
where $\beta=8$, $16$, $24$, respectively, and we utilize the Hadamard integrator to simultaneously solve these three Cauchy problems. Fig.~\ref{figure5} shows the comparisons between the TFT solutions and the \color{b}reference \color{black} solutions at $t=0.5$. Fig. \ref{figure6.0} further compares the solutions along some lines. \color{b} Table \ref{table3} shows the relative $L^2$ and $L^{\infty}$ errors between the TFT solutions and the reference solutions.  
\color{black}

 \begin{figure}[htbp]
     \centering
               \subfigure[]{
     \includegraphics[scale=0.35]{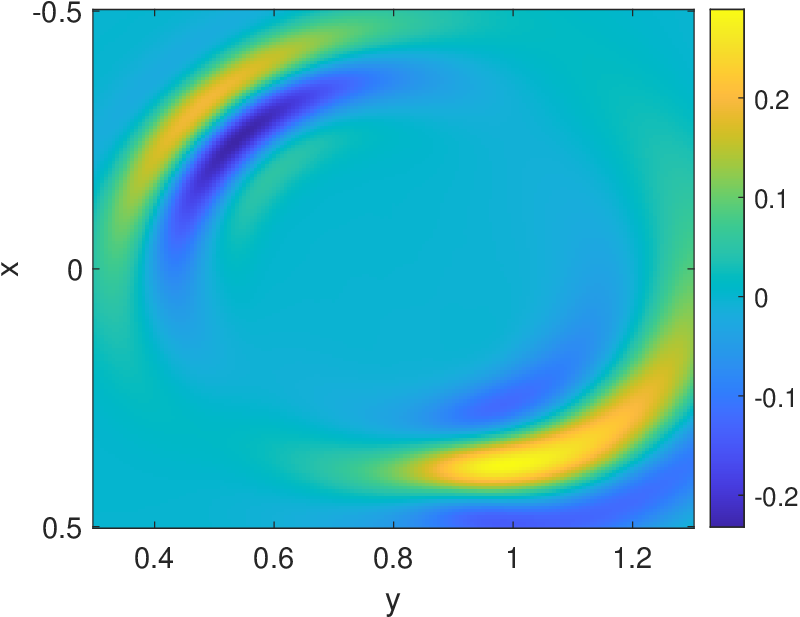}}
                   \subfigure[]{
     \includegraphics[scale=0.35]{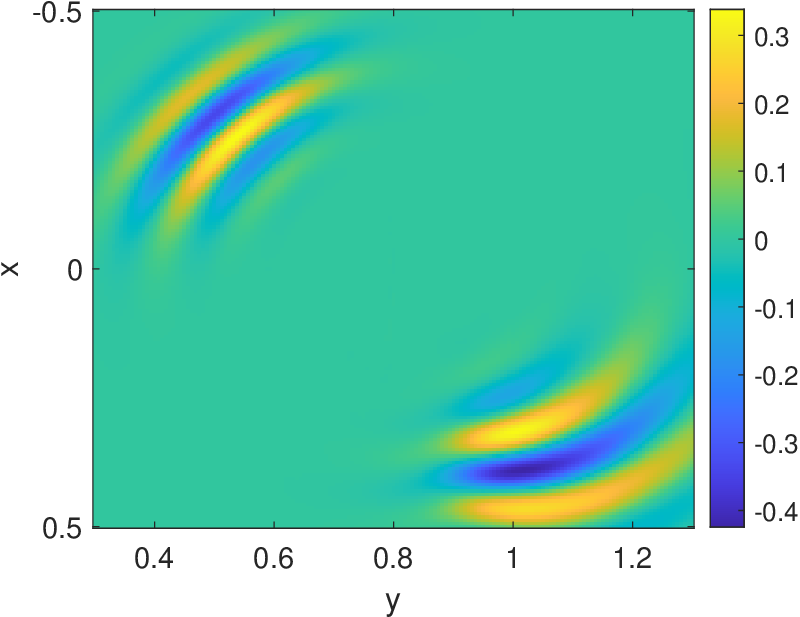}}
                    \subfigure[]{
     \includegraphics[scale=0.35]{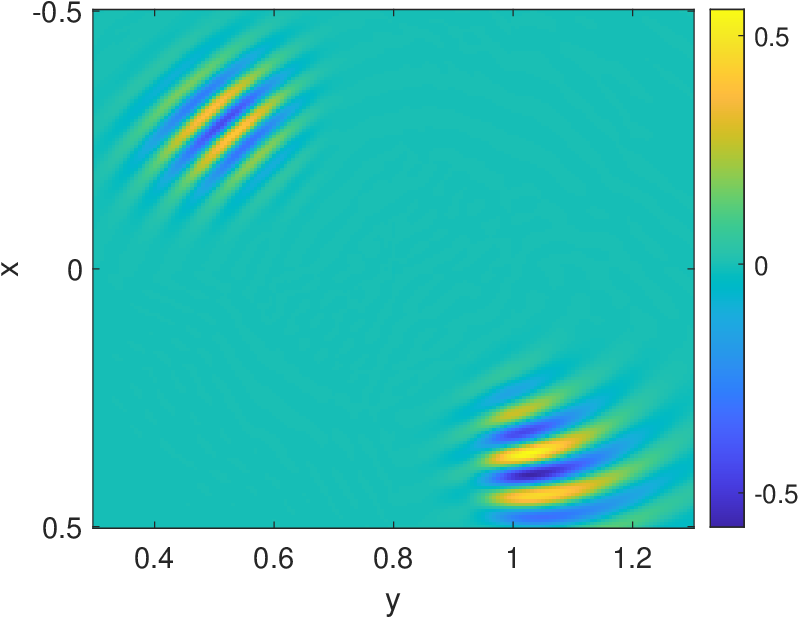}}\\
                \subfigure[]{
     \includegraphics[scale=0.35]{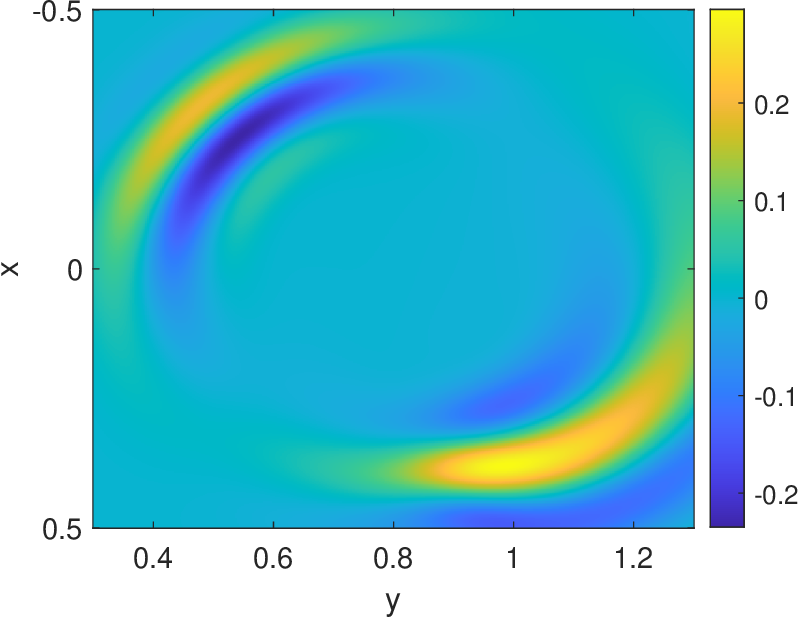}}
                \subfigure[]{
     \includegraphics[scale=0.35]{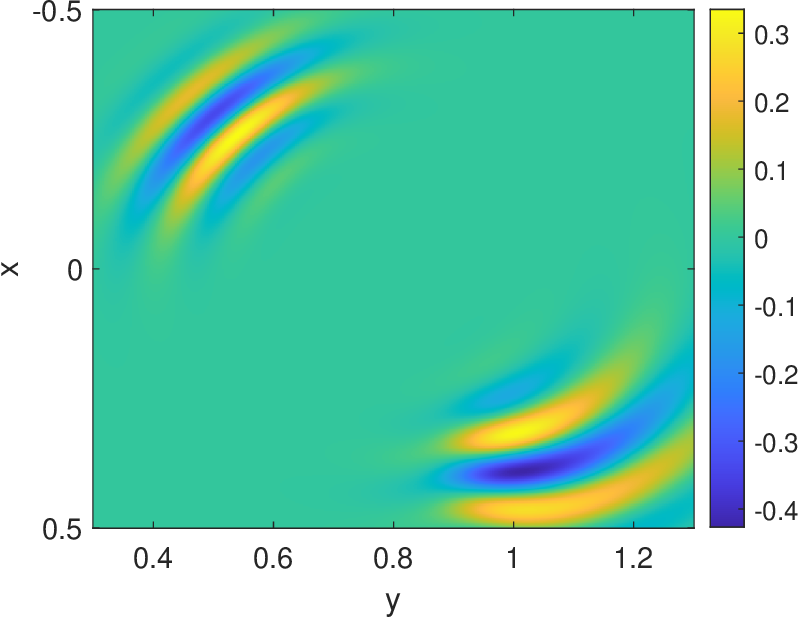}}
                \subfigure[]{
     \includegraphics[scale=0.35]{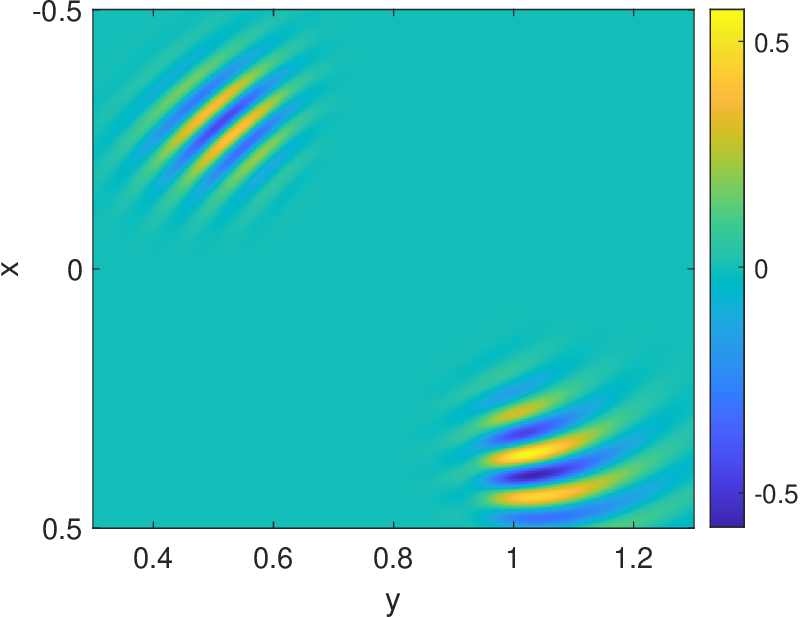}}\\
     \caption{Heaviside model. $u(t,\bx)$, $t=0.5$: (a) TFT solution with $\beta=8$;  (b) TFT solution with $\beta=16$;  (c) TFT solution with $\beta=32$; (d) \color{b}reference \color{black} solution with $\beta=8$;  (e) \color{b}reference \color{black} solution with $\beta=16$; (f) \color{b}reference \color{black} solution with $\beta=32$ }
     \label{figure5}
     \end{figure}
\begin{table}[ht]
\centering
\caption{The relative $L^2$ and $L^{\infty}$ errors of TFT solutions for the Heaviside model at $t=0.5$}\label{table3}
\begin{tabular}{ccccc}
\toprule
$\beta$ & 8 & 16  & 32\\
\midrule
Relative $L^2$ error & $2.76e-2$ & $2.38e-2$  & $4.62e-2$\\
Relative $L^{\infty}$ error & $4.05e-2$ & $2.81e-2$ & $4.92e-2$\\
\bottomrule
\end{tabular}
\end{table}     
\begin{figure}[htbp]
     \centering
                    \subfigure[]{
     \includegraphics[scale=0.35]{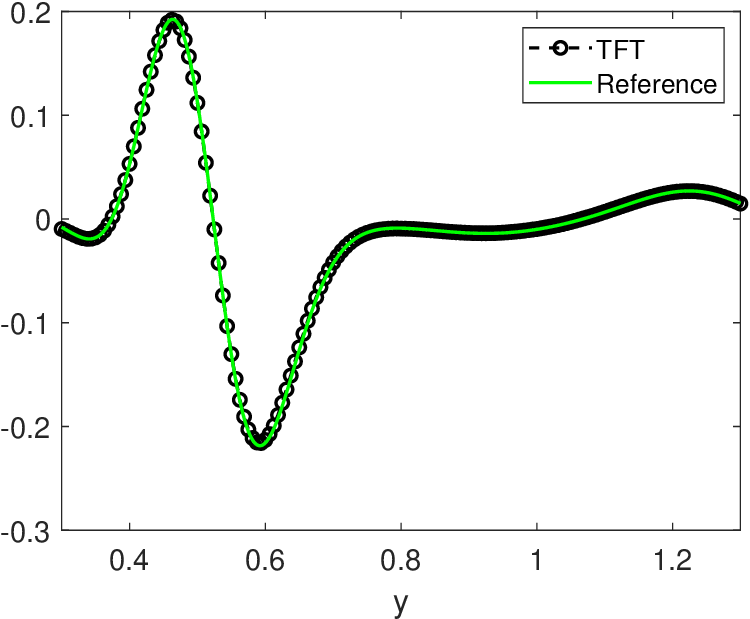}}
                  \subfigure[]{
     \includegraphics[scale=0.35]{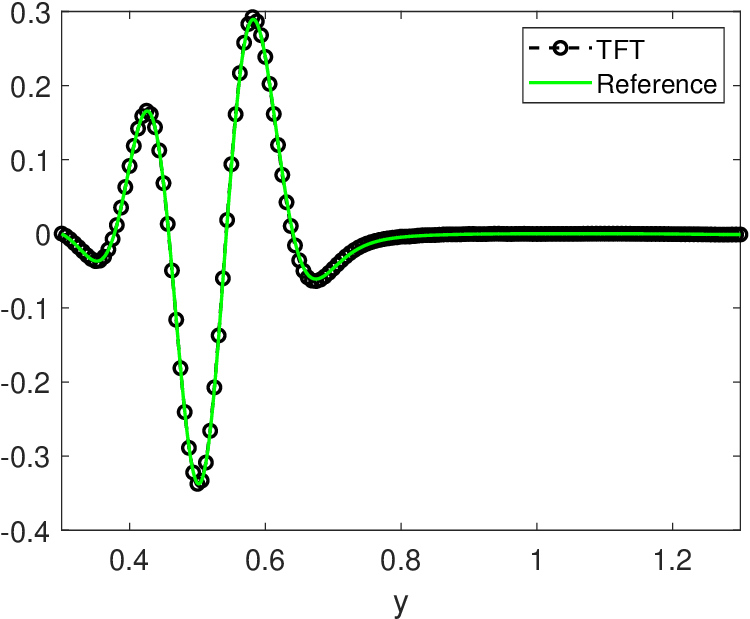}}
                    \subfigure[]{
     \includegraphics[scale=0.35]{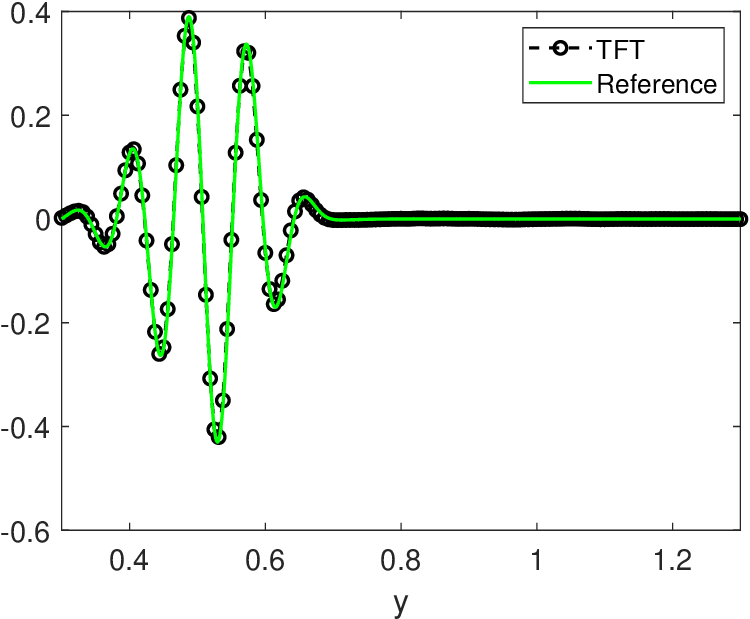}}\\
                \subfigure[]{
     \includegraphics[scale=0.35]{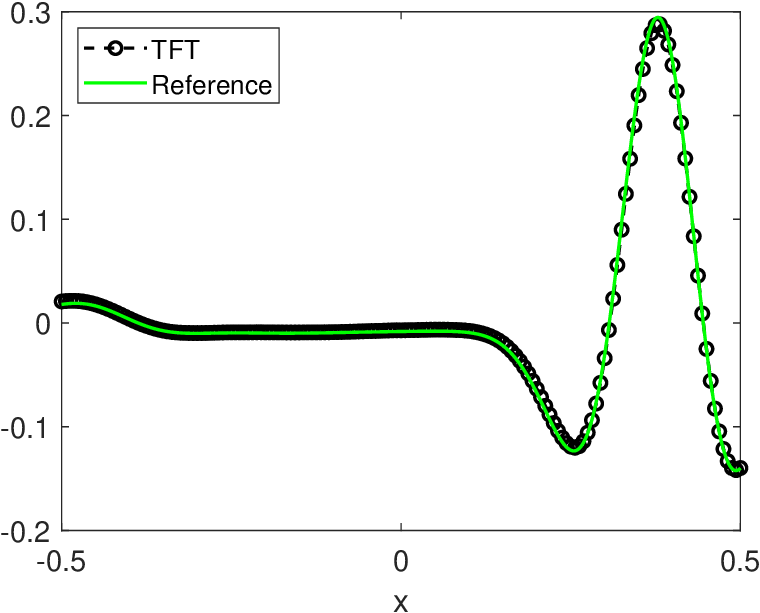}}
                \subfigure[]{
     \includegraphics[scale=0.35]{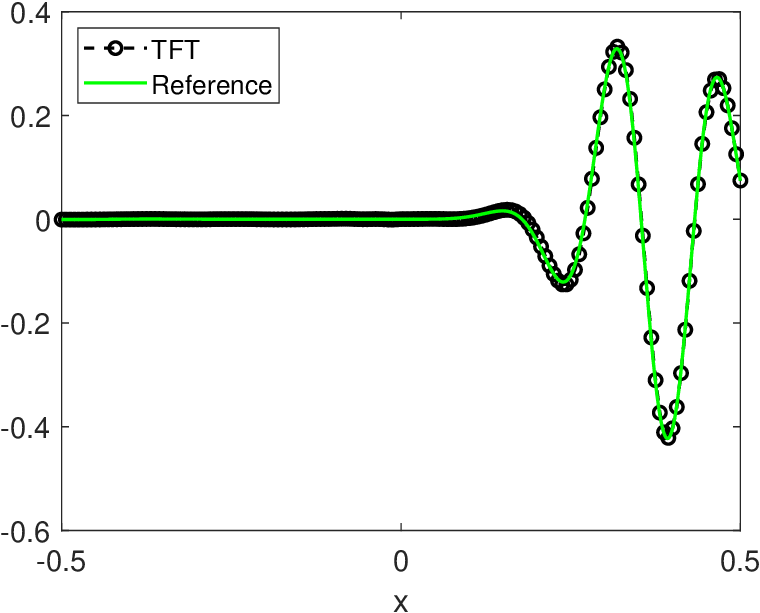}}
                \subfigure[]{
     \includegraphics[scale=0.35]{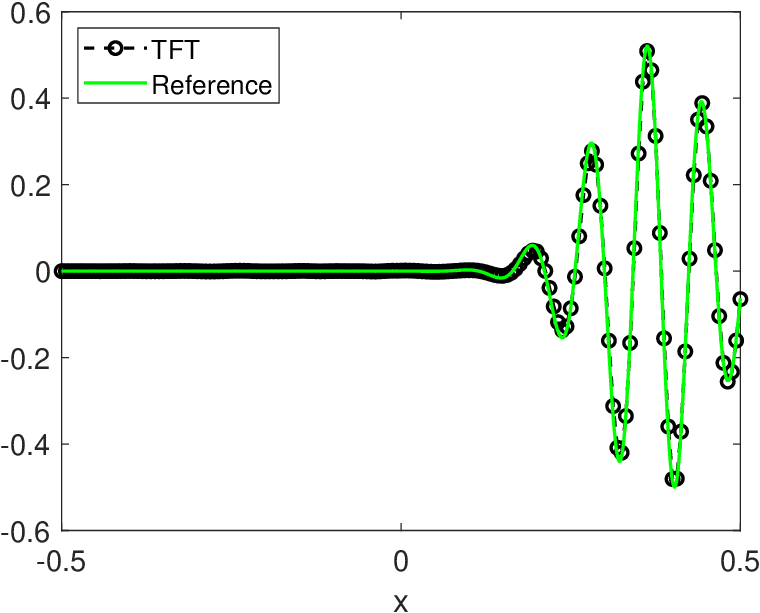}}\\
     \caption{Heaviside model. Slices of time-domain wave fields $u(t,\bx)$ when $t=0.5$: (a) a slice at $x=-0.3$ with $\beta=8$; (b) a slice at $x=-0.3$ with $\beta=16$; (c) a slice at $x=-0.3$ with $\beta=24$; (d) a slice at $y=1$ with $\beta=8$; (e) a slice at $y=1$ with $\beta=16$; (f) a slice at $y=1$ with $\beta=24$
     }
      \label{figure6.0}
     \end{figure}
     
\subsubsection{TFTF method for Case 2}
For Case 2, taking $\bz=[0,0.8]$ and $T_{end}=2$, we obtain the time-domain point-source wave fields $u(t,\bx)$. As shown in Fig. \ref{example3.0}, we can clearly observe the caustics induced by the overturning rays.

Further applying the Fourier transform, we obtain the frequency-domain point-source wave fields $\hat{u}(\omega,\bx)$ with different angular frequencies $\omega$. In Fig. \ref{figure303}, we compare the TFTF solution with the \color{b}reference \color{black} solution for $\omega=8\pi$, $16\pi$, and $24\pi$, respectively, where $\bx\in [-0.6,1.4]\times[0.5,2.5]$. \color{b} We then show in Table \ref{table4} the relative $L^2$ and $L^{\infty}$ errors of the TFTF solutions compared with the reference solutions.  
\color{black}

Fig. \ref{figure404} shows the comparisons of the two solutions along lines that traverse through the caustics region for various frequencies $\omega$. The TFTF solutions match with the \color{b}reference \color{black} solutions well.

    \begin{figure}[htbp]
     \centering
                         \subfigure[]{
     \includegraphics[scale=0.35]{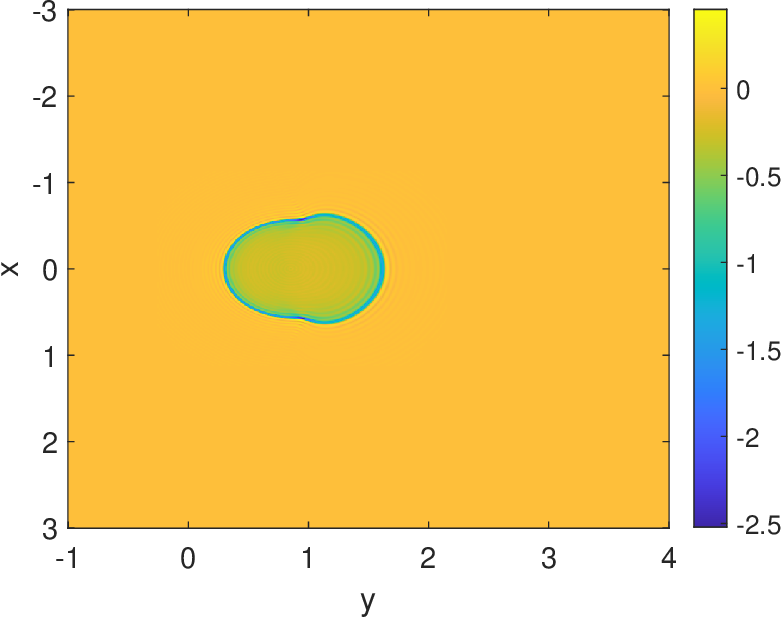}}
                              \subfigure[]{
     \includegraphics[scale=0.35]{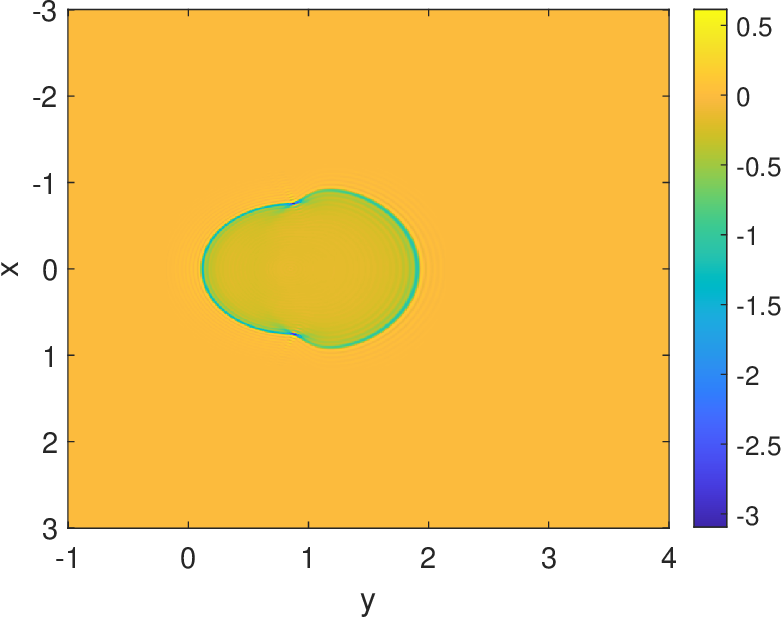}}
                              \subfigure[]{
     \includegraphics[scale=0.35]{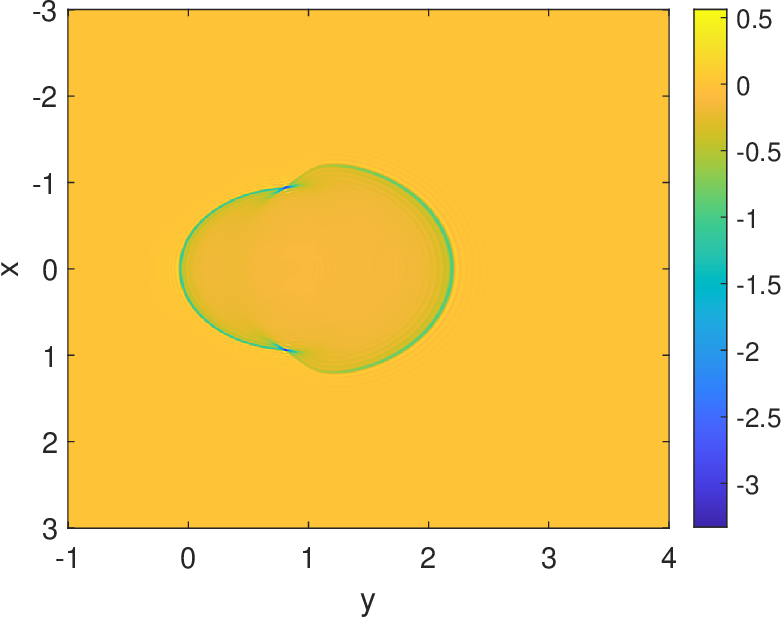}}\\
                              \subfigure[]{
     \includegraphics[scale=0.35]{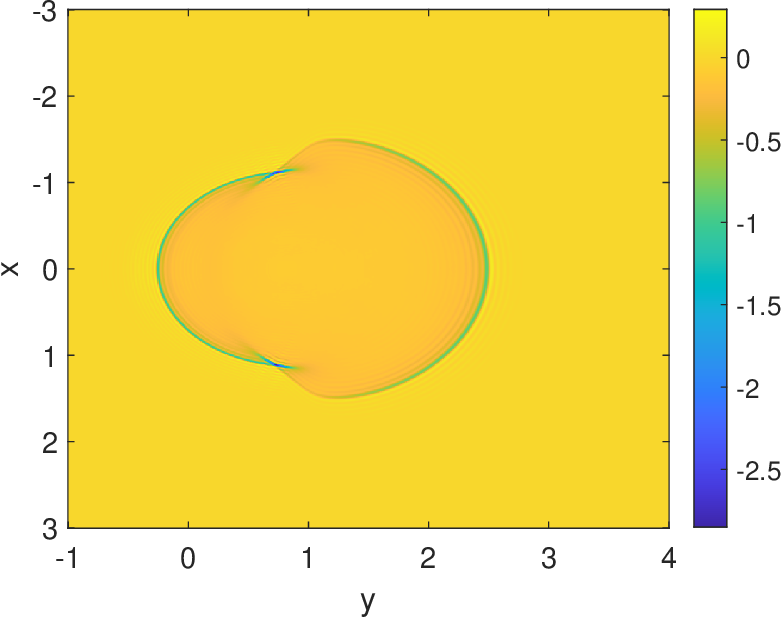}}
                              \subfigure[]{
     \includegraphics[scale=0.35]{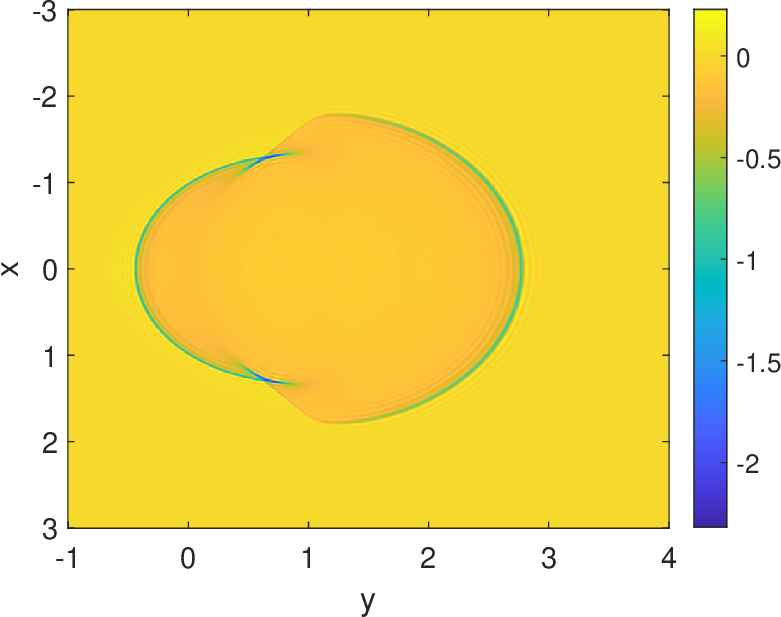}}
                              \subfigure[]{
     \includegraphics[scale=0.35]{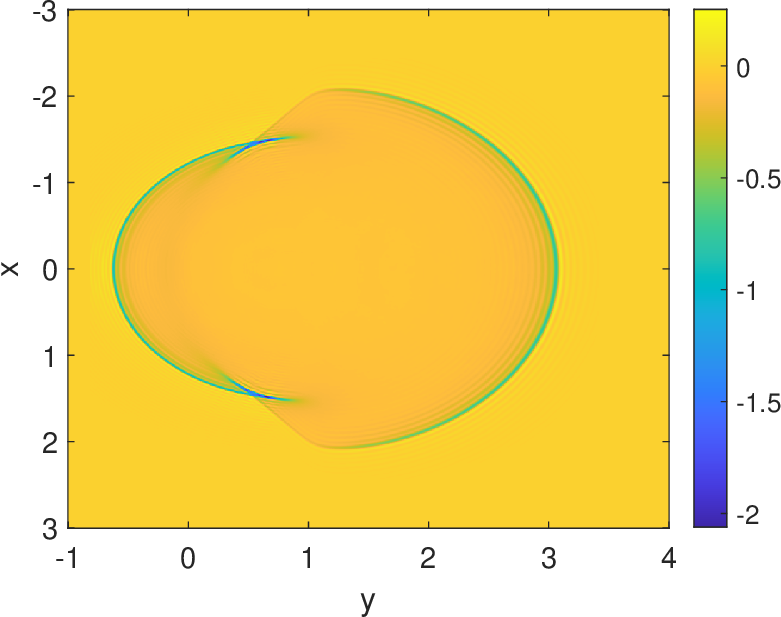}}\\
     \caption{Heaviside model. Time-domain point source wave fields $u(t,\bx)$. (a) $T=0.75$; (b) $T=1$; (c) $T=1.25$; (d) $T=1.5$; (e) $T=1.75$; (f) $T=2$}
     \label{example3.0}
     \end{figure}
     
      \begin{figure}[htbp]
     \centering
          \subfigure[]{
     \includegraphics[scale=0.35]{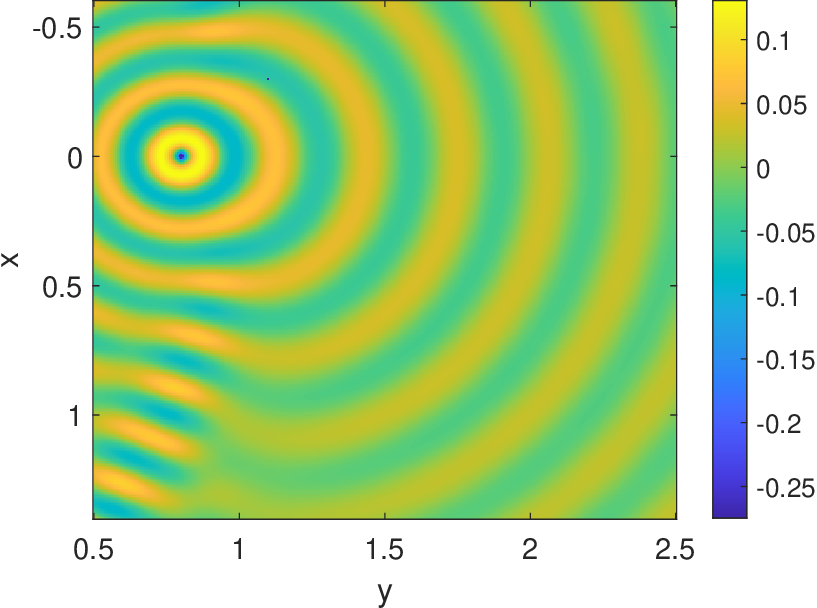}}
                     \subfigure[]{
     \includegraphics[scale=0.35]{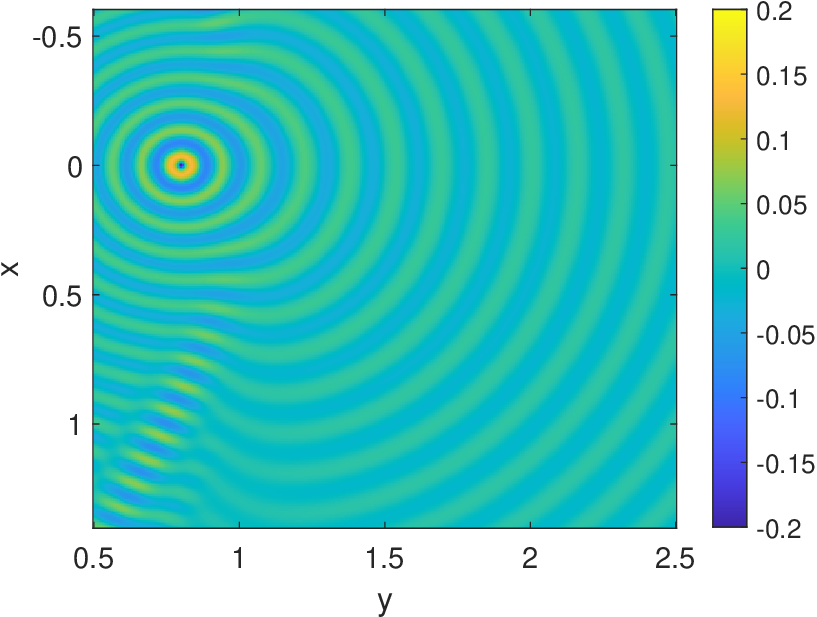}}
          \subfigure[]{
     \includegraphics[scale=0.35]{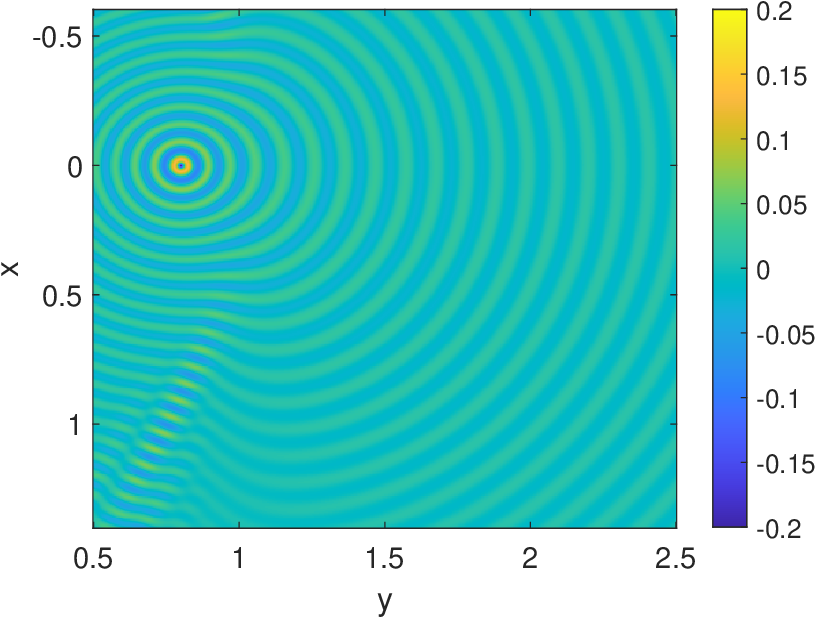}}\\
                    \subfigure[]{
     \includegraphics[scale=0.35]{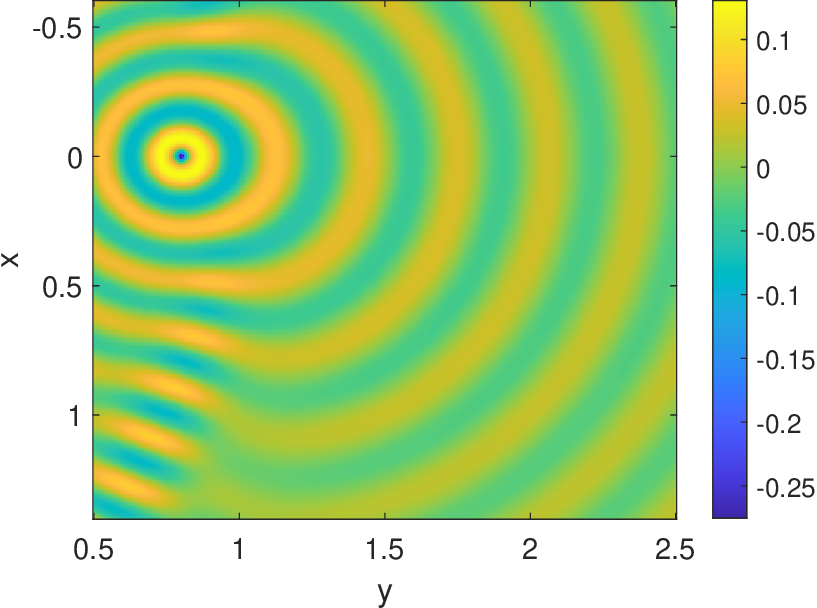}}
                     \subfigure[]{
     \includegraphics[scale=0.35]{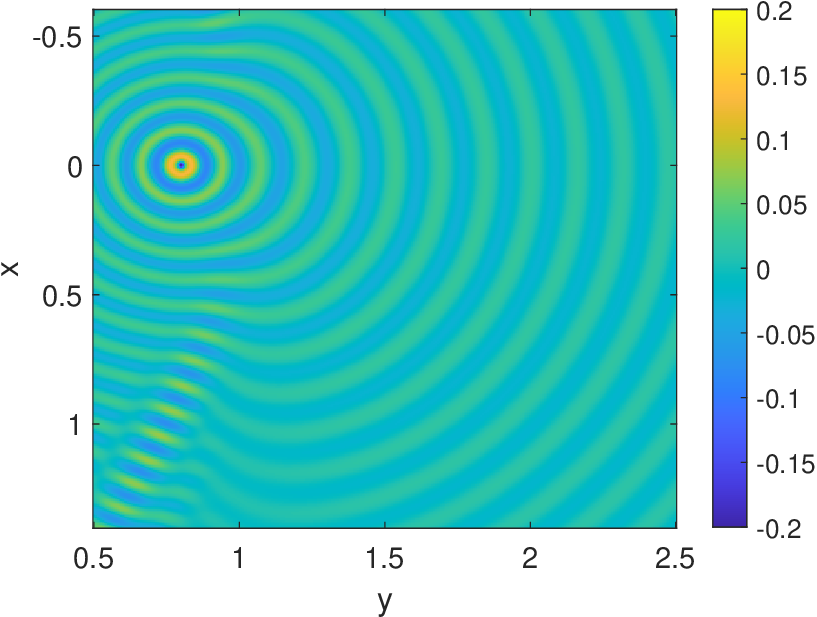}}
          \subfigure[]{
     \includegraphics[scale=0.35]{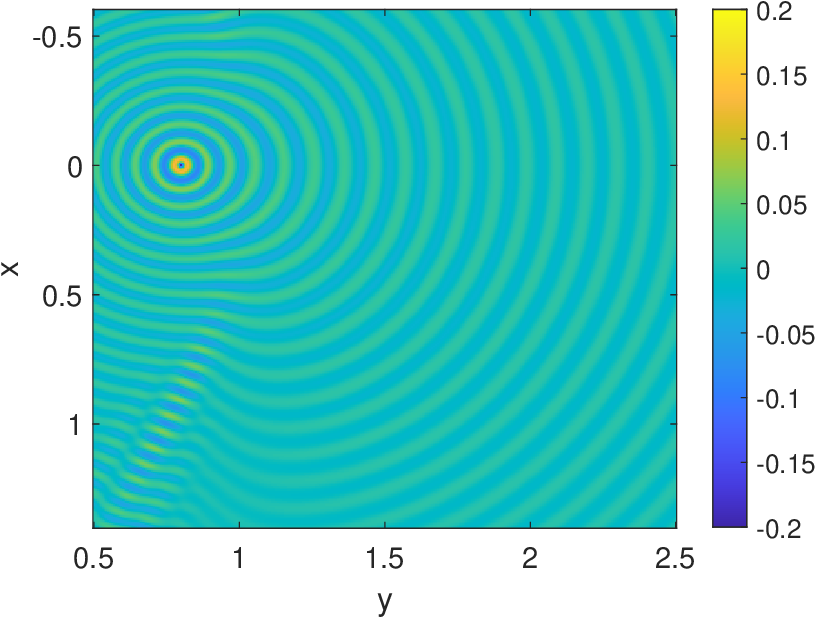}}\\
     \caption{Heaviside model. $\hat{u}(\omega,\bx)$. (a) TFTF solution with $\omega=8\pi$; (b) TFTF solution with $\omega=16\pi$; (c) TFTF solution with $\omega=24\pi$; (d) \color{b}reference \color{black} solution with $\omega=8\pi$;   (e) \color{b}reference \color{black} solution with $\omega=16\pi$; (f) \color{b}reference \color{black} solution with $\omega=24\pi$}\label{figure303}
     \vspace{-3.75mm}
     \end{figure}
     
\begin{table}[ht]
\centering
\caption{The $L^2$ and $L^{\infty}$ errors of TFTF solutions for Heaviside model}\label{table4}
\begin{tabular}{cccc}
\toprule
$\omega$ & $8\pi$ & $16\pi$ & $24\pi$ \\
\midrule
Relative $L^2$ error & $4.33e-2$ & $2.98e-2$ & $ 5.34e-2$   \\
Relative $L^{\infty}$ error & $3.47e-2$ & $3.07e-2$ & $4.91e-2$  \\
\bottomrule
\end{tabular}
\end{table}     
\begin{figure}[htbp]
     \centering
          \subfigure[]{
     \includegraphics[scale=0.35]{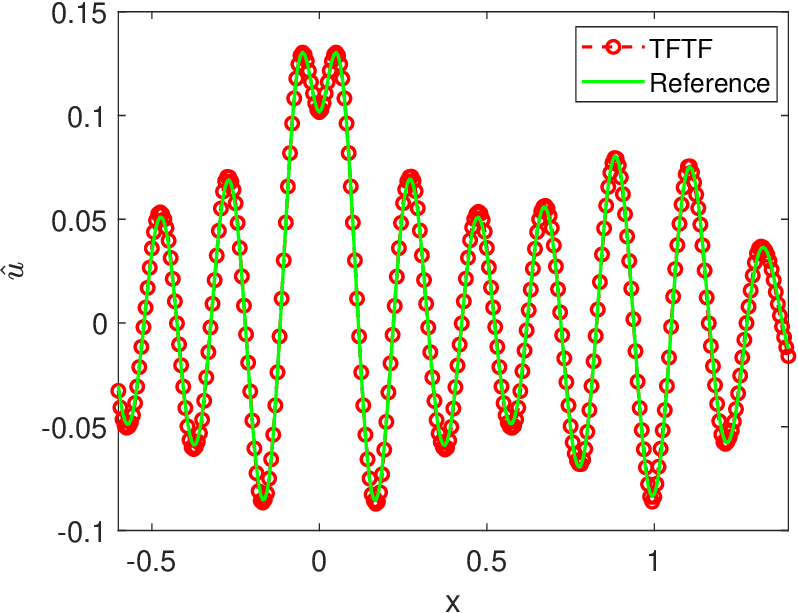}}
               \subfigure[]{
     \includegraphics[scale=0.35]{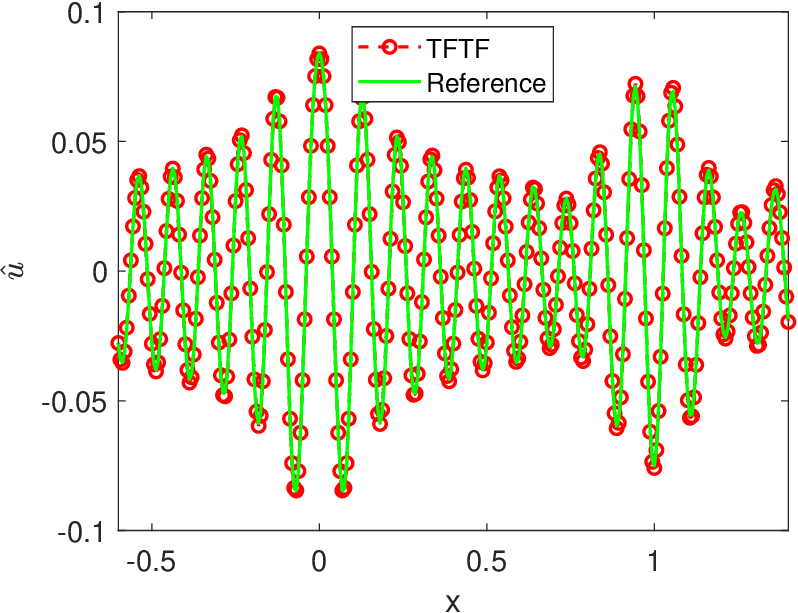}}
                    \subfigure[]{
     \includegraphics[scale=0.35]{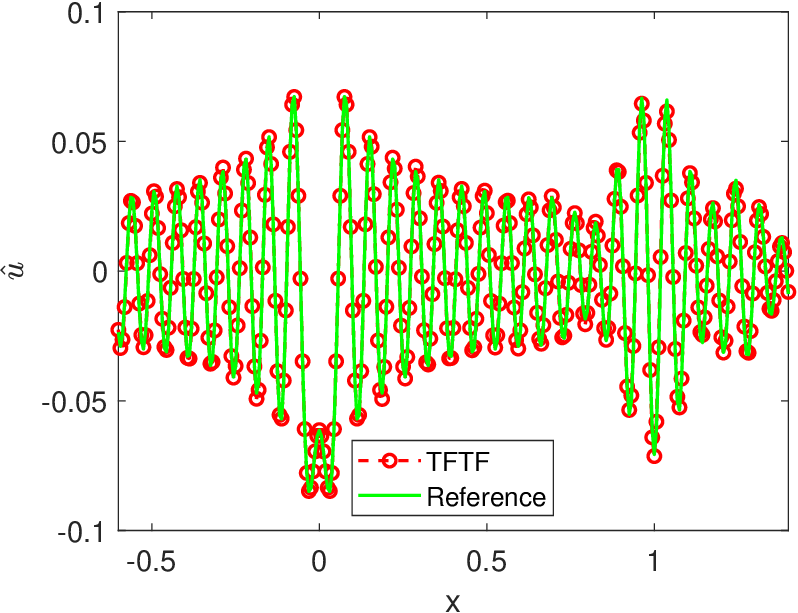}}\\
                \subfigure[]{
     \includegraphics[scale=0.35]{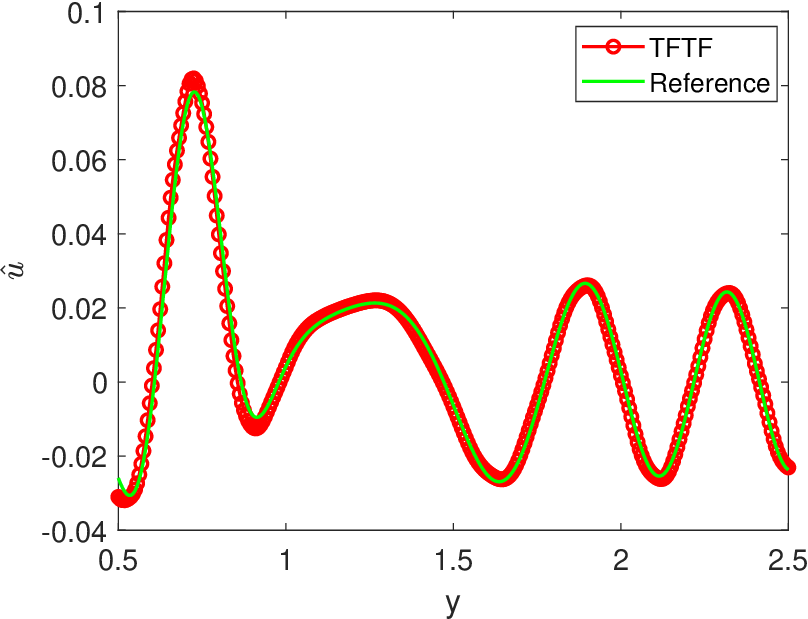}}
                \subfigure[]{
     \includegraphics[scale=0.35]{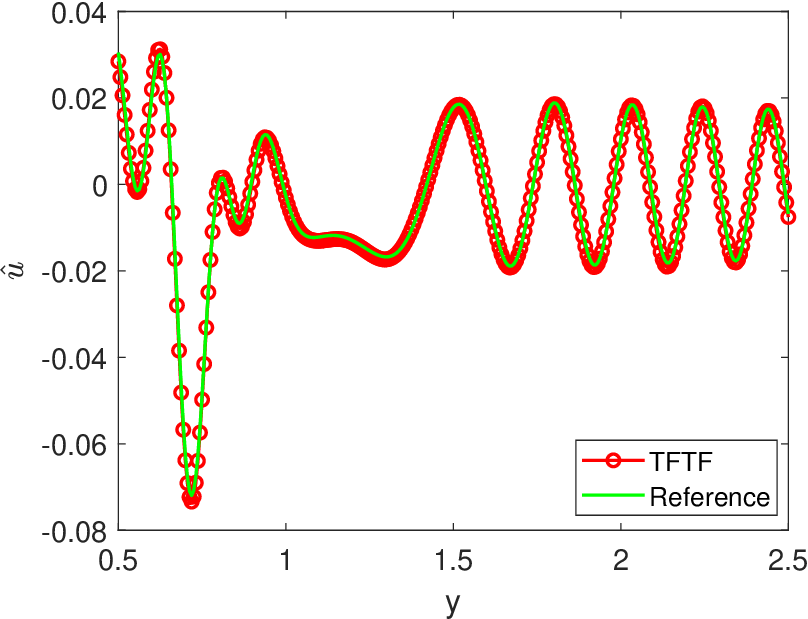}}
                \subfigure[]{
     \includegraphics[scale=0.35]{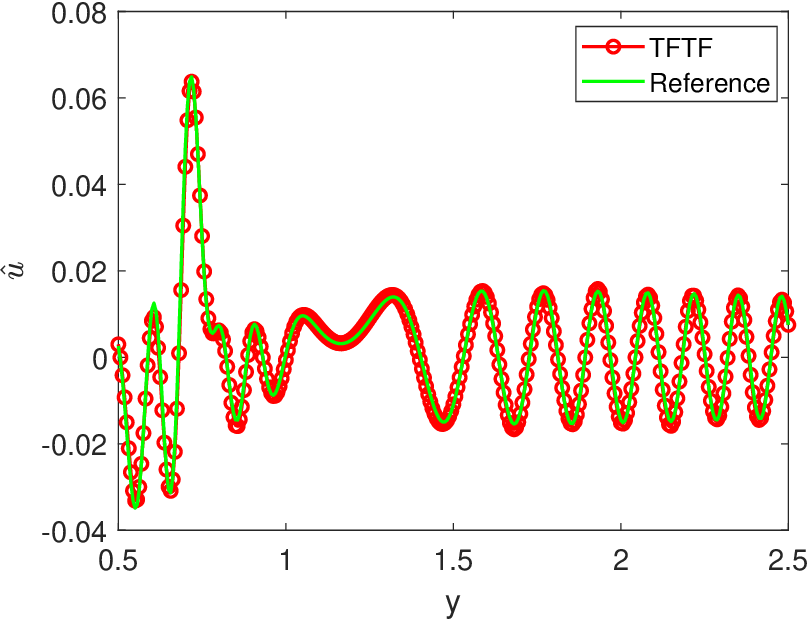}}\\
     \caption{Heaviside model. Slices of frequency-domain point source wave fields $\hat{u}(\omega,\bx)$. (a) a slice at $y=0.75$ with $\omega=8\pi$; (b) a slice at $y=0.75$ with $\omega=16\pi$; (c) a slice at $y=0.75$ with $\omega=24\pi$;
     (d) a slice at $x=1.3$ with $\omega=8\pi$;  (e) a slice at $x=1.3$ with $\omega=16\pi$;  (f) a slice at $x=1.3$ with $\omega=24\pi$}\label{figure404}
     \vspace{-3.75mm}
     \end{figure}
     
\subsection{Complexity validation}
In this subsection, we validate the CPU and memory complexities of the proposed Eulerian Hadamard integrators. Here we use the Sinusoidal model, and we set the source region $\Omega_1=[-0.1,0.1]\times[-0.1,0.1]$ and one of the corresponding adjacent receiver subregions $\Omega_2=[-0.3,-0.1]\times [-0.1,0.1]$. 

Since there exist singularities in our kernels, the compression error will increase as the grid becomes finer; however, we can choose proper tolerance ${\rm tol}$ and oversampling parameters $\chi_1$ and $\chi_2$ so that, even for an extremely fine mesh, such as $h=\frac{1}{3200}$, the maximum relative error is controlled to be smaller than $10^{-4}$, which can be neglected in comparison to other errors. 

We set ${\rm tol} = 10^{-9}$ and $\chi_1=\chi_2=3$, and we vary the frequency and cell count from $\omega=10\pi$ and $N=81^2$ to $\omega=160\pi$ and $N=1281^2$, respectively. We sequentially construct the IDBF matrices of the kernels $U_1^{\omega}(\Omega_1,\Omega_2)$, $U_2^{\omega}(\Omega_1,\Omega_2)$, $U_3^{\omega}(\Omega_1,\Omega_2)$ and $U_4^{\omega}(\Omega_1,\Omega_2)$ and HODBF matrices of the kernels $U_1^{\omega}(\Omega_1,\Omega_1)$, $U_2^{\omega}(\Omega_1,\Omega_1)$, $U_3^{\omega}(\Omega_1,\Omega_1)$ and $U_4^{\omega}(\Omega_1,\Omega_1)$. And we record the total CPU time and total memory usage of IDBF and HODBF, respectively. The CPU time and memory usage for constructing the IDBF and HODBF matrices are plotted, respectively, in Fig. \ref{12} and Fig. \ref{13}, which validate our complexity estimates. Thus, we can state that for given tolerance ${\rm tol}$ and oversampling parameters $\chi_1$ and $\chi_2$, the computational complexity and the CPU memory usage of the butterfly compression grow as $O(N\log(N))$ for IDBF and $O(N\log^2(N))$ for HODBF, respectively. 

     \begin{figure}[htbp]
     \centering
               \subfigure[]{
     \includegraphics[scale=0.35]{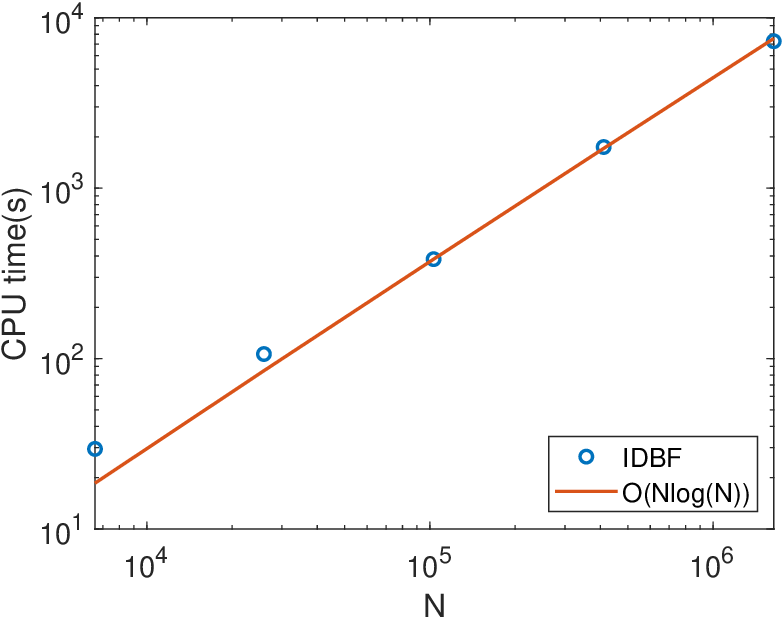}}
                \subfigure[]{
     \includegraphics[scale=0.35]{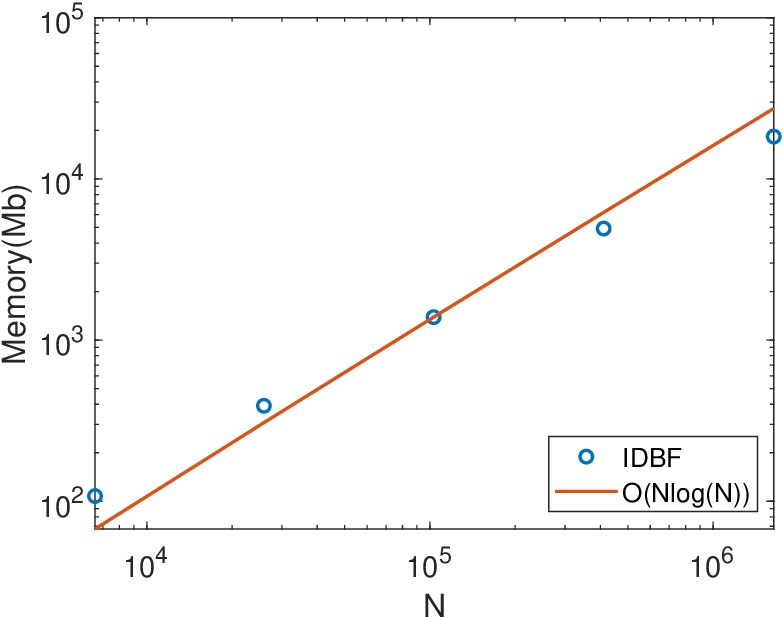}}\\
     \caption{CPU time and memory usage of IDBF without parallel implementation}
    \label{12}
     \end{figure}
     
          \begin{figure}[htbp]
     \centering
              \subfigure[]{
     \includegraphics[scale=0.35]{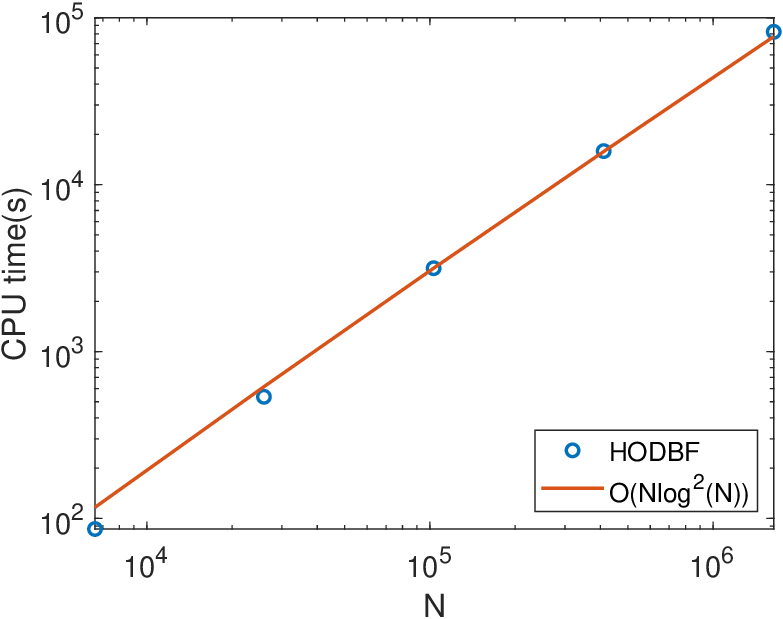}}
                \subfigure[]{
     \includegraphics[scale=0.35]{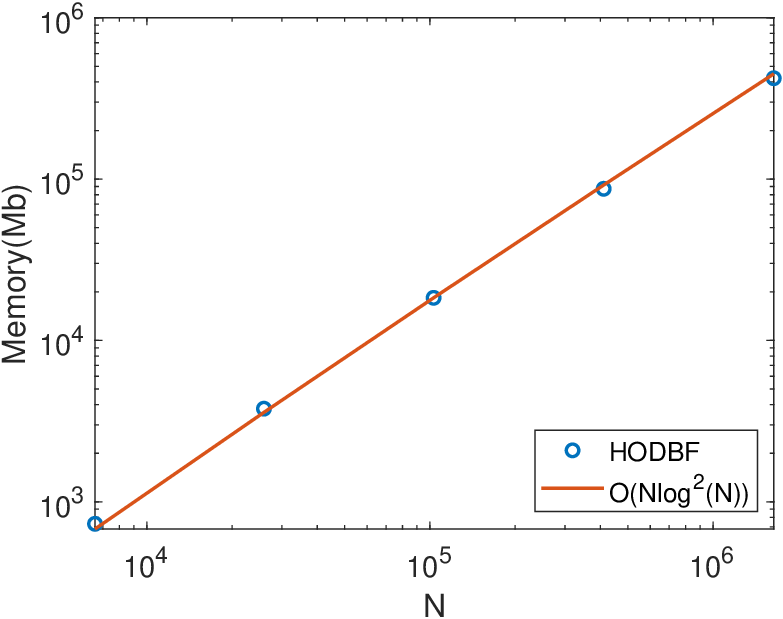}}\\
     \caption{CPU time and memory usage of HODBF without parallel implementation}
     \label{13}
     \end{figure}
     
Due to its quasi-linear memory usage, the computations of IDBF and HODBF can be parallelized quite directly. Although there are $O(B)$ kernels to compress, we will not embed parallelism into the algorithm so that we will only show the computational time of a sequential implementation of the algorithm. Nevertheless, simply using ``parfor'' with respect to $\omega_k$ and $t_\ell$ reduces the computational time of each kernel by roughly $O(N_c)$ times, where $N_c$ is the number of cores of the computer, and parallel implementation on the 56-core computer that we used here can improve the efficiency of the IDBF and HODBF in our numerical examples by roughly about $25$ times.

\subsection{Convergence test}
\label{sec6.4}
As mentioned earlier, $B$ is an artificial frequency bandwidth, and $\omega$ or $\beta$ is the frequency parameter from the problem, both of which affect the accuracy of the Hadamard integrator. 

In the TFT method for Case 1, we have assumed that the time-domain solution $u(t,\bx)$ decays rapidly in the frequency domain, which means that increasing $B$ will not improve the accuracy of the TFT solution. On the other hand, we regard $\beta$ as a frequency parameter, but as $\beta \rightarrow 0$, the frequency-domain bandwidth of the solution does not tend to zero, due to the influence of components in the initial conditions ensuring compact support. Consequently, it is non-trivial to discuss the asymptotic behavior of the TFT method in terms of $B$ or $\beta$, and thus we will not do that here. However, we need to choose an appropriately large $B$ in conjunction with $\beta$ to ensure the rapid decay beyond the bandwidth. 

In the TFTF method for Case 2, generating the time-domain Green's function from the Fourier summation of frequency-domain Green's functions yields an $O(\frac{1}{B^{0.5}})$ asymptotic error according to \rf{3.16.0}. However, when using these time-domain solutions to generate the final frequency-domain solution, we observe an asymptotic behavior of $O(\frac{1}{B})$, which we believe is due to the fact that the TFTF method is somewhat analogous to reconstructing a signal using the Fourier transform. Additionally, the angular frequency also affects the accuracy. According to \cite{liusonburqia23}, taking the one-term Babich's ansatz $\hat{G}_0(\omega,\bx_0;\bx)$ as the frequency-domain Green's function introduces an $O(\frac{1}{\omega^{1.5}})$ asymptotic error. 

To verify the asymptotic behavior with respect to $\omega$, we take the angular frequency $\omega=2\pi$, $4\pi$, $8\pi$, and $16\pi$, respectively, and evaluate the $L^2$-error with respect to the \color{b}reference \color{black} solution. As shown in Fig. \ref{figure12}(a), when $\omega$ is small, the asymptotic error of the one-term Babich's ansatz dominates, resulting in $O(\frac{1}{\omega^{1.5}})$ convergence.

\begin{figure}[htbp]
    \centering
     \subfigure[]{
     \includegraphics[scale=0.35]{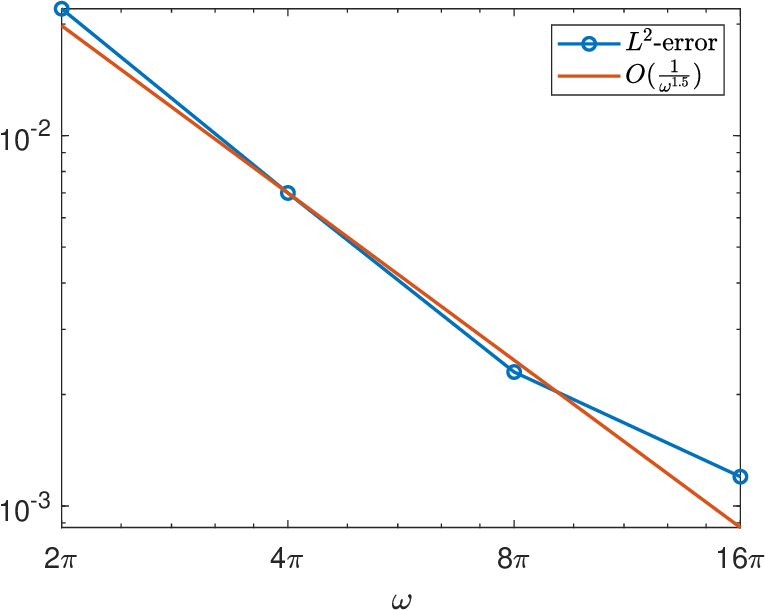}}
     \subfigure[]{
     \includegraphics[scale=0.35]{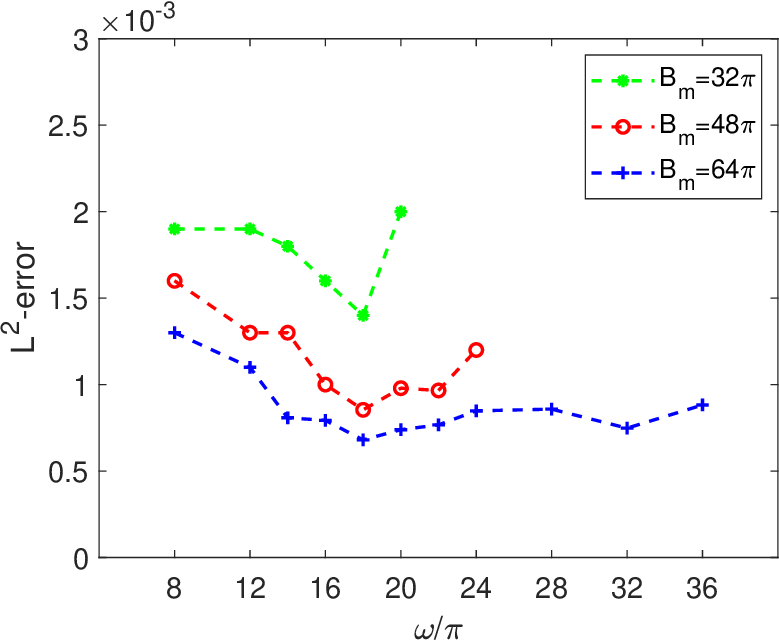}}
\caption{Convergence. (a): $L^2$-error with $B=64\pi$ and different $\omega$; (b): $L^2$-error with different $B$ and $\omega$}
\label{figure12}
\end{figure}

To verify the asymptotic behavior with respect to $B$, we consider the following setups in the Sinusoidal model:
\begin{enumerate}
    \item $B=32\pi$, $h=\frac{1}{100}$, and $\Delta t=\frac{1}{96}$;
    \item $B=48\pi$, $h=\frac{1}{160}$, and $\Delta t=\frac{1}{144}$;
    \item $B=64\pi$, $h=\frac{1}{200}$, and $\Delta t=\frac{1}{192}$.
\end{enumerate}

We further set $T_{end}=1$ and compute the frequency-domain solution in $[0.2,1.2]\times[0.2,1.2]$. 
The $L^2$ errors for different $B$ and $\omega$ are shown in Fig. \ref{figure12}(b). 


For a fixed $B$, as $\omega$ increases, the error initially decreases due to the dominance of the asymptotic error from Babich's ansatz. Subsequently, the error saturates due to the dominance of the asymptotic error caused by truncation of the frequency domain and the fixed $B$. Finally, after $\omega > \frac{B}{2}$, the error increases as the artificial bandwidth $B$ can no longer guarantee a sufficient sampling of the frequency domain. According to the Shannon-Nyquist Sampling theorem, we set $B> 2\omega$. By comparing the errors for different bandwidths $B$ in the interval where the errors do not vary much, we observe an asymptotic behavior of $\frac{1}{B}$. Therefore, we can conclude that when solving the Helmholtz equation, the Hadamard integrator enjoys a frequency-dependent asymptotic convergence in the form of $O\left(\frac{1}{\omega^{1.5}}\right) + O\left(\frac{1}{B^{0.5}}\right)$.

\section{Conclusion}
\label{sec7}
Based on the Kirchhoff-Huygens representation, Hadamard’s ansatz and Babich's ansatz, we developed novel butterfly-compressed Eulerian Hadamard integrators for solving the high-frequency wave equation in time and frequency domains. We derive the Eulerian formulations via the Fourier transform and utilize the butterfly algorithm to accelerate the resulting matrix-vector multiplication. The proposed integrators propagate wave fields beyond caustics implicitly and advance spatially overturning waves naturally, with quasi-linear computational complexity and memory usage. Once constructed, the integrators can simultaneously solve the time-domain wave equations with different initial conditions or the frequency-domain wave equations with different point sources in the computational domain. Numerical examples illustrate the accuracy and efficiency of the new integrators.

\section*{Compliance with Ethical Standards} 
\subsection*{Conflict of interest}
On behalf of all authors, the corresponding author states that there is no conflict of interest.

\section*{Data Availability Statement}
Data sets generated during the current study are available from the corresponding author on reasonable request.

 \setcounter{equation}{0}
\renewcommand\theequation{8.\arabic{equation}}
\begin{appendices}
\section{Computation of Singular Terms}
\label{appendices}
\subsection{Frequency-space integral}
\label{A1}
Consider
\begin{equation}\label{a1}
    I(t,\bx_0;\bx)= \int_0^{\Delta \omega} e^{-i\omega t} \hat{G}_0(\omega,\bx_0;\bx) d\omega,
\end{equation}
which is a frequency-space integral. Since the HKH propagator is short-time valid, both $t$ and $\tau$ are bounded by a sufficiently small $\Delta T$. In the numerical examples, we take $\Delta \omega=\frac{\pi}{2}, \Delta T=\frac{1}{8}$.
The small arguments $\omega t$ and $\omega\tau$ allow us to use the Taylor expansion of $e^{-i \omega t}$ and polynomial approximation of Bessel functions. According to \cite{abrste65},  we have the following polynomial approximations of Bessel functions $J_0(x)$ and $Y_0(x)$ when $0<x<3$,
\begin{equation}
    J_{0}(x)=1-a_1(\frac{x}{3})^{2}+a_2(\frac{x}{3})^{4} +\epsilon,\quad
|\epsilon|=O\left((\frac{x}{3})^6\right),
\end{equation}
and 
\begin{equation}
    Y_{0}(x)=\frac{2}{\pi} \ln (\frac{x}{2}) J_{0}(x)+b_1+b_2(\frac{x}{3})^{2}-b_3(\frac{x}{3})^{4} +\epsilon, 
\quad|\epsilon|=O\left((\frac{x}{3})^6\right),
\end{equation}
where $$a_1=2.2499997, \;a_2=1.2656208, \;b_1=0.36746691,\; b_2=0.60559366, \; b_3=0.74350384.$$
 Using the Taylor expansion of $e^{-i\omega t}$, we obtain 
\begin{equation}\label{3.21}
    \begin{aligned}
ie^{-i\omega t}H_0^1(\omega \tau)\approx& i(1+i\omega t -\frac{\omega^2t^2}{2}-\frac{i \omega^3t^3}{6}+\frac{\omega^4t^4}{24})[
(1+\frac{2i}{\pi}\ln(\frac{\omega \tau}{2}))(1-a_1(\frac{\omega\tau}{3})^2\\
&+a_2(\frac{\omega\tau}{3})^4)+ib_1+ib_2(\frac{\omega\tau}{3})^2+ib_3(\frac{\omega\tau}{3})^4].\\
\end{aligned}
\end{equation}

Dropping fifth- and higher-order terms as well as the imaginary parts, and calculating the definite integral over $[0,\Delta \omega]$, we obtain 
\begin{equation}
   \label{3.22}
   \int_{0}^{\Delta\omega} e^{-i\omega t} \hat{G}_0(\omega,\boldsymbol{x}_0;\bx) d \omega\approx \frac{\sqrt{\pi}}{2}v_0(L_1-L_2),
\end{equation}
where 
\begin{equation}
    \label{3.23}
      \begin{aligned}
L_1  \approx&\frac{1}{(\Delta \omega) \pi}\{- 2 (\Delta \omega)^2+ 2\frac{[(\Delta \omega)^4 (72 t^2 + 16 a_1 \tau^2)]}{648} - \frac{(\Delta \omega)^6 (54 t^4 + 72 a_1 t^2 \tau^2 + 16 a_2 \tau^4)}{16200}\\
& +(\Delta \omega)^2 \ln(\frac{(\Delta \omega) \tau}{2}) + \frac{(\Delta \omega)^6 \ln(\frac{(\Delta \omega) \tau}{2}) (54 t^4+ 72 a_1 t^2 \tau^2 + 16 a_2 \tau^4)}{3240} \\
&- \frac{(\Delta \omega)^4 \ln(\frac{(\Delta \omega)\tau}{2}) (216 t^2 + 48 a_1 \tau^2)}{648}\}-\frac{1}{(\Delta \omega)}[- b_1 (\Delta \omega)^2\\
&+\frac{(\Delta \omega)^4 (108 b_1 t^2 - 24 b_2 \tau^2)}{648} 
  + \frac{(\Delta \omega)^6 (- 27 b_1 t^4 + 36 b_2 t^2 \tau^2 + 8 b_3 \tau^4)}{3240}],
\end{aligned}
\end{equation}
and 
\begin{equation}\label{3.24}
\begin{aligned}
 L_2 & \approx \int_{0}^{\Delta \omega}(\omega t-\frac{\omega^3 t^3}{6}+\frac{\omega^5 t^5}{120})(1-a_1(\frac{\omega\tau}{3})^2+a_2(\frac{\omega\tau}{3})^4)d\omega \\
&\approx\int_{0}^{\Delta \omega} \omega t(1-a_1(\frac{\omega\tau}{3})^2+a_2(\frac{\omega\tau}{3})^4)-\frac{\omega^3 t^3}{6}(1-a_1(\frac{\omega\tau}{3})^2)+\frac{\omega^5 t^5}{120} d\omega\\
&=\int_{0}^{\Delta \omega}\omega t-(a_1\frac{t \tau^2}{9}+\frac{t^3}{6})\omega^3+(a_2\frac{t\tau^4}{81}+a_1\frac{t^3\tau^2}{54}+\frac{t^5}{120})\omega^5 d\omega\\
&=\frac{(\Delta \omega)^2t}{2}-(a_1\frac{t \tau^2}{9}+\frac{t^3}{6})\frac{(\Delta\omega)^4}{4}+(a_2\frac{t\tau^4}{81}+a_1\frac{t^3\tau^2}{54}+\frac{t^5}{120})\frac{(\Delta \omega)^6}{6}.
\end{aligned}
\end{equation}

\subsection{Diagonal term of $I(t,\bx_0;\bx)$}
\label{A2}
 As shown in the appendix of \cite{liusonburqia23}, the integral of self-intersection term inside source cell $c_j$
 \begin{equation}\label{3.20}
 \begin{aligned}
\int_{c_j}\hat{G}_0(\omega_k,\bx_0;\bx)d\bx\approx&\frac{i\sqrt{\pi}}{2}v_0(\bx_0;\bx_0)\frac{1}{\left(n_0 \omega_k\right)^2}\left[8 \int_0^{\frac{\pi}{4}} \frac{h n_0 \omega_k}{2 \cos \theta} H_1^{(1)}\left(\frac{h n_0 \omega_k}{2 \cos \theta}\right) d \theta+4 i\right]\\
=&\frac{i}{4\nu_0\left(n_0 \omega_k\right)^2}\left[8 \int_0^{\frac{\pi}{4}} \frac{h n_0 \omega_k}{2 \cos \theta} H_1^{(1)}\left(\frac{h n_0 \omega_k}{2 \cos \theta}\right) d \theta+4 i\right],\\
 \end{aligned}
 \end{equation}
where $n_0=n(\bx_0)$, $\nu_0=\nu(\bx_0)$, and we use the initial condition \rf{3.9.0}.

Now we consider a self-intersection term of $I(t,\bx_0;\bx)$. Starting from \rf{3.20}, we have 
\begin{equation}\label{3.25}
\begin{aligned}
 \int_{c_j} I(t,\bx_0;\bx) d\boldsymbol{x}&=\int_{0}^{\Delta \omega} \frac{ie^{-i\omega t}}{4\nu_0\left(n_0 \omega\right)^2}\left[8 \int_0^{\frac{\pi}{4}} \frac{h n_0 \omega}{2 \cos \theta} H_1^{(1)}\left(\frac{h n_0 \omega}{2 \cos \theta}\right) d \theta+4 i\right]d \omega\\
 &=\int_{0}^{\frac{\pi}{4}}\int_{0}^{\Delta \omega} \frac{2ie^{-i\omega t}}{\nu_0\left(n_0 \omega\right)^{2}}\left[  \frac{h n_0 \omega}{2 \cos \theta} H_{1}^{(1)}\left(\frac{h n_0 \omega}{2 \cos \theta}\right) +\frac{2i}{\pi} \right] d\omega d \theta\\
 &=\int_{0}^{\frac{\pi}{4}}\int_{0}^{\Delta \omega} \frac{ie^{-i\omega t}h^2}{2\nu_0\cos^2\theta }\left(\frac{2\cos\theta}{h n_0\omega}\right)^2\left[  \frac{h n_0 \omega}{2 \cos \theta} H_{1}^{(1)}\left(\frac{h n_0 \omega}{2 \cos \theta}\right) +\frac{2i}{\pi} \right] d\omega d \theta\\
 &\doteq S_1.
\end{aligned}
\end{equation}

To handle $H_1^{(1)}(s)=J_1(s)+iY_1(s)$, we introduce the polynomial approximation of Bessel functions $J_1$ and $Y_1$ \cite{abrste65} when $0<s<3$,
\begin{equation}
    \label{3.26}
\begin{aligned}
s^{-1} J_{1}(s)=&\frac{1}{2}-0.56249985(\frac{s}{3})^{2}+0.21093573(\frac{s}{3})^{4}-0.03954289(\frac{s}{3})^{6}+0.00443319(\frac{s}{3})^{8} \\
&-0.00031761(\frac{s}{3})^{10}+0.00001109(\frac{s}{3})^{12}+\epsilon, \quad|\epsilon|<1.3 \times 10^{-8},
\end{aligned}
\end{equation}
and 
\begin{equation}\label{3.27}
\begin{aligned}
s Y_{1}(s)=&\frac{2}{\pi} s  \ln (\frac{s}{2})  J_{1}(s)-\frac{2}{\pi} +0.2212091(\frac{s}{3})^{2}+2.1682709(\frac{s}{3})^{4}-1.3164827(\frac{s}{3})^{6}\\
&+0.3123951(\frac{s}{3})^{8} -0.0400976(\frac{s}{3})^{10}+0.0027873(\frac{s}{3})^{12}+\epsilon, \quad|\epsilon|<1.1 \times 10^{-7}.
\end{aligned}
\end{equation}
Thus as $s=\frac{hn_0\omega}{2\cos\theta}\rightarrow 0$, we obtain
\begin{equation}\label{3.28}
\begin{aligned}
\frac{1}{s^2}\left(sH_1^1(s)+\frac{2i}{\pi}\right)&=\frac{1}{s^2}\left[s\left(J_1(s)+iY_1(s)\right)+\frac{2i}{\pi}\right]\\
&=\frac{1}{s}J_1(s)+i(\frac{1}{s^2})\left(sY_1(s)+\frac{2}{\pi}\right)\\
&\sim \frac{1}{2}+i\frac{1}{\pi}\ln(\frac{1}{2}s).
\end{aligned}
\end{equation}

Dropping higher order terms, we have 
\begin{equation}
    \label{3.29}
    s^{-1} J_{1}(s)=\frac{1}{2}-0.56249985(\frac{s}{3})^{2}+0.21093573(\frac{s}{3})^{4}
\end{equation}
for the real part and 
\begin{equation}
    \label{3.30}
    \begin{aligned}
\frac{1}{s^2}(s Y_{1}(s)+\frac{2}{\pi})=&\frac{2}{\pi}\ln (\frac{s}{2})  s^{-1}J_{1}(x)+\frac{1}{s^2}\left(0.2212091(\frac{s}{3})^{2}\right.
\left.+2.1682709(\frac{s}{3})^{4}-0.03954289(\frac{s}{3})^{6}\right)
\end{aligned}
\end{equation}
for the imaginary part. The singularity appearing in the first term of \rf{3.30} can be removed in the following way, 
\begin{equation}
    \label{3.31}
    \begin{aligned}
S_1=&\int_{0}^{\frac{\pi}{4}} \frac{ih^2}{2\nu_0\cos^2\theta}\int_{0}^{\Delta \omega}e^{-i\omega t} \left(\frac{2\cos\theta}{h n_0 \omega}\right)^{2}\left[  \frac{h n_0 \omega}{2 \cos \theta} H_{1}^{(1)}\left(\frac{h n_0 \omega}{2 \cos \theta}\right) +\frac{2i}{\pi} \right] d\omega d \theta\\
=&\int_{0}^{\frac{\pi}{4}} \frac{i h^2}{2\nu_0\cos^2\theta}\int_{0}^{\Delta \omega}e^{-i\omega t} \left(\frac{2\cos\theta}{h n_0 \omega}\right)^{2}\left[  \frac{h n_0 \omega}{2 \cos \theta} H_{1}^{(1)}\left(\frac{h n_0 \omega}{2 \cos \theta}\right) +\frac{2i}{\pi} \right]\\
&-\frac{i}{\pi}\ln(\frac{h n_0 \omega}{4 \cos \theta}) d\omega d \theta
+\int_{0}^{\frac{\pi}{4}}\int_{0}^{\Delta \omega}  \frac{i^2h^2}{2\nu_0\cos^2\theta}\frac{1}{\pi}\ln(\frac{h n_0 \omega}{4 \cos \theta}) d\omega d \theta,
\end{aligned}
\end{equation}
where the last singular term in \rf{3.31} can be handled as 
\begin{equation}
    \label{3.32}
    \begin{aligned}
&\int_{0}^{\frac{\pi}{4}}\int_{0}^{\Delta \omega}  \frac{i^2h^2}{2\nu_0\cos^2\theta}\frac{1}{\pi}\ln(\frac{h n_0 \omega}{4 \cos \theta}) d\omega d \theta\\
=&-\frac{2h^2}{\pi\nu_0}\int_{0}^{\frac{\pi}{4}}\int_{0}^{\Delta \omega}  \frac{1}{\cos^2\theta}\ln(\frac{h n_0 \omega}{4 \cos \theta}) d\omega d \theta\\
=&-\frac{2h^2}{\pi\nu_0}\int_{0}^{\frac{\pi}{4}}\int_{0}^{\Delta \omega}  \frac{1}{\cos^2\theta}\left[\ln(h n_0\omega)-\ln(4 \cos \theta)\right] d\omega d \theta\\
=&-\frac{2h^2}{\pi\nu_0}\int_{0}^{\frac{\pi}{4}}  \frac{1}{\cos^2\theta}[\Delta\omega(\ln(\frac{ h n_0 \Delta\omega}{4\cos\theta})-1)] d \theta.\\
\end{aligned}
\end{equation}

\subsection{Diagonal term of $\hat{G}_0(\omega,\bx_0;\bx)\left(\fr{\partial}
{\tau \partial \tau}\right)\left( \fr{u^k}{v_0} \right)$}
\label{A3}
Consider the self-intersection terms \rf{dig3} with non-zero $\omega$ in cell $c_j$, where $n$ and non-singular ingredients are taken to be  constant. We  state that
\begin{equation}\label{A.14}
    S_2\doteq \int_{c_j}\hat{G}_0(\omega,\bx_0;\bx)\left(\fr{\partial}{\tau \partial \tau}\right)\left(\frac{u^k}{v_0}\right) d\boldsymbol{x}=0, \quad k=1,2,
\end{equation}
which implies that the self-intersection term is zero.

To see this, we note that, when $n=n_0$,
\begin{equation}
    \frac{\partial }{\partial \tau}=c^2\nab \tau \cdot \nab=\frac{1}{n_0}(\cos\theta,\sin\theta)\cdot \nab,
\end{equation}
where $\theta$ is the take-off angle, which is consistent with that in the polar coordinate transformation when $n=n_0$. And
\begin{equation}
\hat{G}_0(\omega,\bx_0;\bx)=\frac{i}{4\nu_0}H_0^1(\omega n_0|\bx-\bx_0|),
\end{equation}
which is a radial function. Then we can obtain \eqref{A.14} by symmetry. Actually, we can divide the cell $c_j$ into four parts
\begin{equation}
    \left\{
\begin{aligned}
&c_j^1 :\theta\in [-\frac{\pi}{4},\frac{\pi}{4}], r\in[0,\frac{h}{2\cos(\theta)}],\\
&c_j^2 :\theta\in [\frac{\pi}{4},\frac{3\pi}{4}], r\in[0,\frac{h}{2\sin(\theta)}],\\
&c_j^3 :\theta\in [\frac{3\pi}{4},\frac{5\pi}{4}], r\in[0,-\frac{h}{2\cos(\theta)}],\\
&c_j^4 :\theta\in [\frac{5\pi}{4},\frac{7\pi}{4}], r\in[0,-\frac{h}{2\sin(\theta)}].
\end{aligned}
\right.
\end{equation}
Now we use the polar coordinate transformation in \rf{A.14}, yielding  
\begin{equation}
\begin{aligned}
    S_2=&\int\int \frac{i}{4\nu_0}H_0^1(\omega n_0|\bx-\bx_0|) \frac{1}{n_0}(\cos\theta,\sin\theta)\cdot \nab\left(\frac{u^k}{v_0}\right) \frac{1}{n_0 |\bx-\bx_0|} |\bx-\bx_0| d r d\theta \\
    =& \frac{i}{4\rho_j }(\int_{c_j^1}+\int_{c_j^2}+\int_{c_j^3}+\int_{c_j^4}) H_0^1(\omega n_0|\bx-\bx_0|) (\cos\theta,\sin\theta)\cdot \nab\left(\frac{u^k}{v_0}\right)  d r d\theta \\
    \doteq&S_{2,1}+S_{2,2}+S_{2,3}+S_{2,4}.
\end{aligned}
    \end{equation}
    Substituting $\xi=\theta+\pi$, it is easy to verify that 
$$
S_{2,1}+S_{2,3}=0,\quad S_{2,2}+S_{2,4}=0.
$$
That is, $S_2=0$.
\end{appendices}


%
%

\bibliographystyle{spmpsci}      

\bibliography{myref}   

%
%

\end{document}